\documentclass[review]{elsarticle}
\usepackage{latexsym,amsfonts,amssymb,amsthm,amsmath,mathrsfs,color,amscd,graphicx,fullpage,parskip,hyperref,bm,bbm}

\usepackage{lineno,hyperref}
\modulolinenumbers[5]
\journal{arXiv:2002.11055}

\usepackage{nicefrac}
\usepackage{comment}
\usepackage{commath,mathtools,setspace}
\usepackage{paralist}
\usepackage{extarrows}

\includecomment{uncomment} 
\excludecomment{comment} 

\newcommand{\s}[1]{{\mathcal #1}}

\newcommand{\bb}[1]{{\mathbb #1}}

\DeclareMathOperator*{\essinf}{ess\, inf}

\newtheorem{theorem}{Theorem} 
\newtheorem{corollary}[theorem]{Corollary}

\newtheorem{lemma}[theorem]{Lemma}

\newtheorem{definition}[theorem]{Definition}

\newtheorem{remark}[theorem]{Remark}
\newtheorem{assumption}[theorem]{Assumption}

\numberwithin{equation}{section}
\numberwithin{theorem}{section}
\renewcommand{\color}[2]{#2}

\begin{document}

\begin{frontmatter}

\title{Nonlocal Bertrand and Cournot Mean Field Games \\with General Nonlinear Demand Schedule} 

\author{P. Jameson Graber\fnref{fn1}}
\address{P. J. Graber: Baylor University, Department of Mathematics;
	One Bear Place \#97328;
	Waco, TX 76798-7328 
}
\ead{Jameson\_Graber@baylor.edu}

\fntext[fn1]{P. Jameson Graber is grateful to be supported by the National Science Foundation under NSF Grants DMS-1612880 and DMS-1905449.}

\author{Vincenzo Ignazio\fnref{fn2}}
\address{Vincenzo Ignazio: ETH Zurich, Department of Mathematics;
	R\"amistrasse 101,
	8092 Zurich, Switzerland
}
\ead{vincenzo.ignazio@math.ethz.ch}

\fntext[fn2]{Vincenzo Ignazio is grateful to be partly supported by the Swiss National Science Foundation grant SNF~200020$\_$172815.}

\author{Ariel Neufeld\fnref{fn3}\corref{mycorrespondingauthor}}
\address{Ariel Neufeld: NTU Singapore, Division of Mathematical Sciences;
		21 Nanyang Link 
		Singapore 637371
}
\ead{ariel.neufeld@ntu.edu.sg}

\fntext[fn3]{Ariel Neufeld is grateful to be  supported by his NAP Grant \textit{Machine Learning based Algorithms in Finance and Insurance.}}

\cortext[mycorrespondingauthor]{Corresponding author}


\begin{abstract}
		In this article we prove the existence of classical solutions to a system of mean field games arising in the study of exhaustible resource production under market competition.
Individual trajectories are modeled by a controlled diffusion process with jumps, which adds a nonlocal term to the PDE system.
The assumptions on the Hamiltonian are sufficiently general to cover a large class of examples proposed in the literature on Bertrand and Cournot mean field games.
Uniqueness also holds under a sufficient restriction on the structure of the Hamiltonian, which in practice amounts to a small upper bound on the substitutability of goods.
\end{abstract}

\begin{keyword}
	mean field games of controls
	\sep  
	extended mean field games
	\sep 
	 economic models
	 \sep
	   Bertrand competition\sep 
	    Cournot competition
\MSC[2010] 35Q91\sep  	35Q93  \sep  	91A13 
\end{keyword}

\end{frontmatter}


\section{Introduction}
The intent of this paper is to study a {\color{red} system of} partial differential equations arising from the theory of mean field games, of the form
\begin{equation}\label{eq:Main}
\begin{array}{lc}
u_t+\frac{1}{2}\sigma^2 u_{xx}-ru+H(t,u_x,{\color{red}\pi}[u_x,m],\eta)+{\s{I}}[u]=0,&  0<t<T,\ 0<x<\infty\\
m_t-\frac{1}{2}\sigma^2 m_{xx}-(D_\xi H(t,u_x,{\color{red}\pi}[u_x,m],\eta)m)_x-\s{I}^*[m]=0,& 0<t<T,\ 0<x<\infty\\
\eta(t) = \int_0^\infty m(t,x)\dif x, & 0 < t < T,\\
m(0,x)=m_0(x),\ u(T,x)=u_T(x), & 0\leq x< \infty\\
u(t,x)=m(t,x)=0, & 0\leq t\leq T,  -\infty< x\leq 0.
\end{array}
\end{equation}
where $T>0$ is a known terminal time, $m_0$ and $u_T$ are given smooth functions. The operator $H$, {\color{red} which among others has a dependency on a non-local coupling {\color{blue} function} ${\color{red}\pi}[u_x,m]$ (defined in Assumption~\ref{as:coupling}),}  is determined by a given equilibrium price-demand function, which also characterizes the market on which the problem is modeled on. Here ${\s{I}}$ is defined as a non-local operator of the form 
\begin{equation}
\mathcal{I}[u](t,x):=\int_{-\infty}^\infty \big[u(t,x+z)-u(t,x)-u_x(t,x)z\mathbf{1}_{\{|z|\leq 1\}}\big]\,F(dz),
\end{equation}
where $F(dz)$ denotes a L\'evy measure. 
{\color{red}The adjoint operator  
 $\s{I}^*$ of ${\s{I}}$ is given by
\begin{equation*}
	\s{I}^*[m](t,x):=
	\int_{-\infty}^\infty \big[m(t,x-z)-m(t,x)-m_x(t,x)z\mathbf{1}_{\{|z|\leq 1\}}\big]\,F(dz).
\end{equation*}}
These equations represent a mean field game based on a competition between producers with different amount of capacities, which are based on ideas 
first introduced by Gu\'eant, Lasry, and Lions \cite{gueant2010paris}
and then further explored by Chan and Sircar in \cite{chan2015bertrand}, where the authors studied a model with an exhaustible source. Further results on the existence and uniqueness for this kind of equations were proven, under suitable conditions, by Bensoussan and Graber in \cite{graber2018existence}. Inspired by their work, we provide a generalization of their equations by allowing them to incorporate additional abrupt changes in the available resource quantities.\\\\  
Mean field games were first introduced in \cite{huang2006large} and \cite{lasry07} in order to analyze the behavior of a large crowd of players in a differential game.
The entire ``mass'' of players is approximated through a density function $m(x,t)$, by viewing them as a homogeneous continuum. The arising equations are given by the coupling of a backward Hamilton--Jacobi--Bellman equation, describing the optimal control aspect of individuals in the game, with a forward Fokker--Planck equation, characterizing the change in the density of the players (or sub-density, in the case that players can leave the game). 
\\\\
To date, most results studying the existence and uniqueness of this equation either couple the HJB and Fokker--Planck equations through a local coupling, as in 
\cite{porretta2015weak, 
	cardaliaguet2015second,
	cardaliaguet2012long, 
	gomez2014time,
	cardaliaguet2013weak,
	cardaliaguet2015mean}, or through a non-local coupling, as in \cite{cardaliaguet2015long}. 
It is also worth mentioning that different authors have approximated solutions numerically both for the local case \cite{briceno2018long} and for the non-local case \cite{nurbekyan2018long}. Nonetheless, there is currently no available literature which simultaneously considers a non-local coupling inside the Hamiltonian, the addition of a non-local linear operator describing sudden changes (``jumps") in the given data and an unbounded domain for the {\color{red} state} variable~$x$.
(See, however, \cite{cesaroni2019stationaryfractional,cirant2019fractionalmfg} for some mean field games with a fractional Laplacian.)
This additional nonlocal coupling inside the Hamiltonian derives from the feature that players in the mean field game maximize an objective that depends on the distribution of controls.
Such games are now called in the literature ``mean field games of controls" or else ``extended mean field games," which have been studied in some generality from a PDE point of view in the following references: \cite{gomes2014extended,gomes2016extended,cardaliaguet2018mfgcontrols,kobeissi2019mfgcontrols,bonnans2019schauder}.
Although the results of this article are tailored to a specific application, the general techniques extend many of these results to a broader class of Hamiltonians and nonlocal interactions.\\\\
The main purpose of this article is to prove, under general conditions, that a (classical) solution for the system \eqref{eq:Main} exists, and that uniqueness can be obtained under certain restrictions on the non-linear Hamiltonian.
In particular, we do not require the Hamiltonian to have any growth restriction as the variable $\abs{u_x}$ tends to infinity; on the other hand, it possesses certain monotonicity properties which have a regularizing effect on the solution.
Furthermore, we show that our assumptions are satisfied for a general class of price-demand functions appearing in economic applications. Our results, for example, apply to the entire collection of examples proposed by Chan and Sircar in \cite{chan2017fracking} (cf.~\cite{harris2010games}).  We consider only the case of a non-vanishing positive diffusion constant~$\sigma$, in order to avoid degenerate parabolic equations.\\\\
In Subsection~\ref{sec:Multi-Per} we provide motivation for the model in terms of Bertrand and Cournot competition in the mean field limit.
Then in Section~\ref{sec:Main-Res} we introduce the model, the required notation, the spaces in which our analysis will take place, and the required assumptions on the data. Additionally, we state in Theorem~\ref{thm:main-result} the two main results of the paper (existence and uniqueness), which are proven in Section~\ref{sec:existenceP} and Section~\ref{sec:uniqueness}.\\\\
In Section~\ref{sec:auxiliary} we prove analogues of standard results, such as the maximum principle or bootstrapping estimates, in the case where there is an additional linear non-local jump term $\mathcal{I}[u]$ present. In fact, due to the linearity of the operator, it behaves like a ``perturbation" which does not interact ``too strongly" with the non-local coupling of the gradient $u_x$.\\\\
In Section~\ref{sec:apriori} we establish a priori estimates in order to obtain bounds on $u$, $u_x$ and $m$, taking into consideration the final data $u_T$. These estimates lead to an upper bound for the non-local coupling ${\color{red}\pi}[u_x,m]$ {\color{red}defined in Assumption~\ref{as:coupling}}, which we exploit for the results in the following sections. \\\\
Section~\ref{sec:existenceP} is centered around establishing the existence of a (classical) solution $(u,m)$ for the system of equations \eqref{eq:Main}.  These results are based on the Leray--Schauder fixed point theorem and rely on the assumptions of local Lipschitz continuity and boundedness of the Hamiltonian $H$. \\\\ 
In  Section~\ref{sec:uniqueness} the solution is shown to be unique, granted that the Hamiltonian satisfies certain \textit{smallness} and \textit{uniform convexity} conditions, defined below in Subsection~\ref{subsec:ass-unique}.
In economic applications, the conditions under which uniqueness is proved amount to a small upper bound on the degree to which competing firms can be substituted for one another.

\subsection{Bertrand and Cournot Type Mean Field Games with sensitive market changes}\label{sec:Multi-Per}
{\color{red}The mean field game we consider is one in which producers sell their stock in an optimal way, given the market is in equilibrium; see {\color{blue} the derivation found in} \cite{chan2017fracking} and \cite{gueant2011mean}.}
Let $t$ be a time and  {\color{red} let the state variable} $x$ be the producer's capacity. We assume to have already such a big number of players that the mean field games ``limit" has already been obtained and is therefore described by the density of the players, namely the variable $m(t,x)$.\\
We introduce the variable 

\begin{equation}
\eta(t)=\int_0^{\infty} m(t,x)dx,\ \ 0\leq\eta(t)\leq 1,
\end{equation}
which represents the total mass of producers remaining into the ``game", i.e., those with positive stock. The function $\eta(t)$ is decreasing in time.\\\\
The producer's capacity is driven by a stochastic differential equation with L\'evy jump term $Z$ of the form

\begin{equation}
\begin{cases}
dX(s)=-q(s,X(s))ds+\sigma dW(s)+dZ(s)\\
X(0)=X_0
\end{cases}
\end{equation}
where $m(0,x)=m_0(x)$ describes the density of the initial random variable $X_0$ and $Z$ is a L\'evy process of the form $Z=\int_0^\cdot \int_{\bb R} z\,\mathbf{1}_{\{|z|\leq 1\}}[\mu^Z(dz,dt)-F(dz)dt]+\int_0^\cdot \int_{\bb R}z\,\mathbf{1}_{\{|z|> 1\}}\,\mu^Z(dz,dt)$, where $\mu^Z(dz,dt)$ denotes its jump measure  with compensator $F(dz)\,dt$ for a L\'evy measure $F(dz)$. The corresponding integro operator of $Z$ is then of the form
\begin{equation*}
\mathcal{I}[v](t,x):=\int_{-\infty}^\infty \big[v(t,x+z)-v(t,x)-v_x(t,x)z\mathbf{1}_{\{|z|\leq 1\}}\big]\,F(dz).
\end{equation*}
Moreover, we introduce two different time variables, $T$ and $\tau$. The former variable describes the time horizon, while the latter is a so called stopping time (i.e.\ a random variable) that describes the time when our process reaches zero, which means that the producer will have no more capacity and leaves the market.\\\\
The {\color{red}factor $q\equiv q(t,x)$ determines} the quantity of produced material, 
{\color{red} whereas the factor $p\equiv p(t,x)$ describes the price demanded for the goods.}
	The quantity $q$ is in a natural equilibrium with the demand of the good, namely 
\begin{equation} \label{eq:demand schedule}
q=\s{D}^{\eta}(p,\bar{p}),
\end{equation} 
which depends on the mass of players $\eta$ remaining in game, the price $p$ the consumers are willing to pay for buying the good and $\bar{p}$ the average price offered by all producers.

Alternatively, $q$ can also be related to $p$ via an \emph{inverse demand schedule}, namely
\begin{equation} \label{eq:inverse demand schedule}
p=P(q,Q),
\end{equation}
where $Q$ is the total market quantity produced.
In some cases, such as when the demand schedule is linear, it can be shown that \eqref{eq:demand schedule} and \eqref{eq:inverse demand schedule} are equivalent \cite{chan2015bertrand}.

\subsubsection{Bertrand competition}

Bertrand competition refers to the situation in which producers set prices in order to maximize profit.
Thus, the variable $p$ is the control.
The optimal utility is given by 
\begin{equation}\label{eq:maxproblem-Bertrand}
u(t,x):=\sup_{p\geq 0}\bb{E} \left\{\int^{T\wedge\tau}_t e^{-r(s-t)}{\color{red}p(s,X(s))q(s,X(s))}\,ds+e^{-r(T-t)}u_T(X(T))\,\bigg|\,\ X(t)=x\right\}
\end{equation}  
where  the supremum is taken over {\color{red}Markov controls given by the price random fields 
$p\equiv \{p(s,x)\}_{s\in [t,T],x\in[0,\infty)}$, 
	and where $q\equiv \{q(s,x)\}_{s\in [t,T],x\in[0,\infty)}$ denotes the quantity random fields of produced material characterized by the equilibrium relation given in \eqref{eq:demand schedule}.}

The corresponding Hamilton-Jacobi-Bellman equation has the following form 
\begin{equation}
u_t+\frac{1}{2}\sigma^2u_{xx}-ru+\max_{p\geq 0}[\s{D}^{\eta}(p,\bar{p})(p-u_x)]+{\s{I}}[u]=0.
\end{equation}  

The average price is defined as
\begin{equation}
\bar{p}=\frac{1}{\eta(t)}\int_0^{\infty}p^{*}(t,x)m(t,x)dx,
\end{equation}
where $p^*(t,x)$ is the Nash equilibrium price.
We then set the Hamiltonian to be
\begin{equation} \label{eq:hamiltonian Bertrand}
H = \s{D}^{\eta}(p^*,\bar{p})(p^*-u_x)
\end{equation}
and the equilibrium production rate is
\begin{equation} \label{eq:production Bertrand}
\s{D}^{\eta}(p^*,\bar{p}). 
\end{equation}

\subsubsection{Cournot competition}

Cournot competition refers to the situation in which producers set their level of production in order to maximize profit.
Thus, the variable $q$ is the control.
The optimal utility is given by 
\begin{equation}\label{eq:maxproblem-Cournot}
u(t,x):=\sup_{q\geq 0}\bb{E} \left\{\int^{T\wedge\tau}_t e^{-r(s-t)}{\color{red}p(s,X(s))q(s,X(s))}\,ds+e^{-r(T-t)}u_T(X(T))\,\bigg|\, \ X(t)=x\right\}
\end{equation}  
where the supremum is taken over {\color{red}Markov controls given by  the quantity random fields of produced material
$q\equiv \{q(s,x)\}_{s\in [t,T],x\in[0,\infty)}$, 
and where $p\equiv \{p(s,x)\}_{s\in [t,T],x\in[0,\infty)}$ denotes the price random fields  characterized by the equilibrium relation given in \eqref{eq:inverse demand schedule}.}

The corresponding Hamilton-Jacobi-Bellman equation has the following form
\begin{equation}
u_t+\frac{1}{2}\sigma^2u_{xx}-ru+\max_{q\geq 0}[q(P(q,Q)-u_x)]+{\s{I}}[u]=0.
\end{equation}  

The total market production is defined as
\begin{equation}
Q=\int_0^{\infty}q^{*}(t,x)m(t,x)\dif x,
\end{equation}
where $q^*(t,x)$ is the Nash equilibrium quantity produced.
We then set the Hamiltonian to be
\begin{equation} \label{eq:hamiltonian Cournot}
H = q^*(P(q^*,Q)-u_x).
\end{equation}

\section{System of equations and main result}\label{sec:Main-Res}
\subsection{System of Equations}
Observe that whether the competition is of Bertrand or Cournot type, the final system of equations reads in the general form
\vspace{0.2cm}
\begin{equation}\label{eq:HJBs2}
\begin{array}{lc}
u_t+\frac{1}{2}\sigma^2 u_{xx}-ru+H(t,u_x,{\color{red}\pi}[u_x,m],\eta)+{\s{I}}[u]=0,&  0<t<T,\ 0<x<\infty\\
m_t-\frac{1}{2}\sigma^2 m_{xx}-(D_\xi H(t,u_x,{\color{red}\pi}[u_x,m],\eta)m)_x-\s{I}^*[m]=0,& 0<t<T,\ 0<x<\infty\\
\eta(t) = \int_0^\infty m(t,x)\dif x, & 0 < t < T,\\
m(0,x)=m_0(x),\ u(T,x)=u_T(x), & 0\leq x< \infty\\
u(t,x)=m(t,x)=0, & 0\leq t\leq T,  -\infty< x\leq 0.
\end{array}
\end{equation}
Here $m_0$ represents a given initial distribution of competitors, $u_T$ is the final cost, $H = H(t,\xi,\upsilon,\eta)$ is a given function of three real variables, $\mathcal{I}^*$ denotes the adjoint operator of $\mathcal{I}$ (see also \cite[Section~2.4]{garroni2002second}), and $\upsilon = {\color{red}\pi}[u_x,m]$ in a nonlocal way--thus $H(t,u_x,{\color{red}\pi}[u_x,m],\eta)$ depends on the pointwise value $u_x = u_x(t,x)$ but also on the pair $u_x,m$ in a nonlocal way.
See Assumptions \ref{as:coupling} and \ref{as:coupling continuity}.

The goal of this work is to provide existence and uniqueness of a {\color{magenta}(}classical{\color{magenta})} solution for the  above system of equation \eqref{eq:HJBs2}. For the precise statement of the definition of a solution \eqref{eq:HJBs2} as well as our main result formulated in Theorem~\ref{thm:main-result}, we refer to Subsection~\ref{subsec:main-result}.
\subsection{Notation} \label{sec:notation}
In this subsection, we introduce the notation we will be using in this paper.

$\bullet$ As is customary, we will always use $C$ to denote a generic constant depending on the data appearing in the problem.
Occasionally subscripts or arguments will specify precise quantities on which $C$ depends, in which case we have that $C$ is bounded for bounded values of its parameters.\\
$\bullet$ Moreover, $\alpha \in (0,1)$ will also denote a generic constant, but it will be reserved for exponents; its value may \emph{decrease} from line to line, depending only on the data.\\
$\bullet$ For a domain $\Omega \subset [0,\infty)$ and $p\in[1,\infty)$, we denote by $L^p([0,T] \times \Omega)$ the space of measurable functions $u$ such that $\enVert{u}_p = \enVert{u}_{L^p} = \enVert{u}_{L^p([0,T] \times \Omega)} := \del{\int_0^T \int_{\Omega} \abs{u}^p \dif x \dif t}^{1/p} < \infty$.\\
$\bullet$ We denote by $L^p_{loc}([0,T] \times [0,\infty))$ the space of all measurable functions $u$ such that \\
$\enVert{u}_{L^p_{loc}} := \sup_{M \geq 0} \enVert{u}_{L^p([0,T] \times [M,M+1])} < \infty$.
We point out that this is a strict subspace of the space of all locally $L^p$ functions.\\
$\bullet$ We also use the notation $L^p_{\geq 0}$ to denote the subset of $L^p$ functions which are non-negative a.e.

Let $k$ be a non-negative integer and 
$\alpha \in (0,1)$.

$\bullet$ We denote by $W^k_p(\Omega)$ the space of functions $u$ such that the weak partial derivatives satisfy $\pd[j]{u}{x} \in L^p(\Omega)$ for all $0 \leq j \leq k$ with norm $\enVert{u}_{W^k_p} = \sum_{j=0}^k \enVert{\pd[j]{u}{x}}_p$.\\
$\bullet$ The space $C^k(\Omega)$ denotes all bounded functions $u$ that are $k$ times continuously differentiable having corresponding bounded derivatives, with $\enVert{u}_{C^k} = \sum_{j=0}^k \enVert{\pd[j]{u}{x}}_\infty$; hence $C^k$ is a subspace of $W_\infty^k$ with $\enVert{u}_{C^k} = \enVert{u}_{W_\infty^k}$.\\
$\bullet$ For any non-integer $r>0$ we denote by $W^r_p(\Omega) $ the usual fractional Sobolev space, see also  \cite[Chapter~2]{ladyzhenskaia1968linear}.

$\bullet$ We denote by $C^\alpha(\Omega)$ the space of H\"older continuous functions with norm
\begin{equation}
\enVert{u}_{C^\alpha} = \enVert{u}_\infty + \intcc{u}_{C^\alpha},
\end{equation}
with $\enVert{\cdot}_\infty$ denoting the supremum norm and $\intcc{\cdot}_{C^\alpha}$ denoting the H\"older seminorm
\begin{equation}
\intcc{u}_{C^\alpha} = \sup_{x \neq y} \frac{\abs{u(x)-u(y)}}{\abs{x-y}^\alpha}.
\end{equation}
$\bullet$ The space $C^{k+\alpha}(\Omega)$ is the subset of all functions $u \in C^k(\Omega)$ such that $\pd[k]{u}{x} \in C^\alpha(\Omega)$, with norm given by $\enVert{u}_{C^{k+\alpha}} = \enVert{u}_{C^k} + \intcc{\pd[k]{u}{x}}_{C^\alpha}$.

Finally, we define parabolic spaces used in our work.

$\bullet$ Let $W^{1,2}_p([0,T] \times \Omega)$ be the space of functions $u$ such that $u,u_t,u_x,u_{xx} \in L^p([0,T] \times \Omega)$ with norm given by
\begin{equation}
\enVert{u}_{W^{1,2}_p} = \enVert{u}_{L^p} + \enVert{u_t}_{L^p} + \enVert{u_x}_{L^p} + \enVert{u_{xx}}_{L^p}.
\end{equation}
$\bullet$ We denote by $C^{\alpha/2,\alpha}([0,T] \times \Omega)$ the space of H\"older continuous functions where the corresponding seminorm is given by
\begin{equation}
\intcc{u}_{C^{\alpha/2,\alpha}} = \sup_{(t,x) \neq (\tau,y)} \frac{\abs{u(t,x)-u(\tau,y)}}{\abs{t-\tau}^{\alpha/2} + \abs{x-y}^\alpha}.
\end{equation}
%
$\bullet$ Then we denote by $C^{1+\alpha/2,2+\alpha}([0,T] \times \Omega)$ the space functions $u$ such that $u,u_t,u_x,u_{xx}$ are all contained in $C^{\alpha/2,\alpha}([0,T] \times \Omega)$, with norm
\begin{equation}
\enVert{u}_{C^{1+\alpha/2,2+\alpha}} = \enVert{u}_{C^{\alpha/2,\alpha}} + \enVert{u_t}_{C^{\alpha/2,\alpha}} + \enVert{u_x}_{C^{\alpha/2,\alpha}} + \enVert{u_{xx}}_{C^{\alpha/2,\alpha}}.
\end{equation}

\subsection{Assumptions on data} \label{sec:assms}
In order to obtain classical solutions in the sense described above, we will assume that the data are sufficiently regular. We first introduce the sufficient conditions on the data to obtain existence of a classical solution; see Subsubsection~\ref{subsec:ass-exist}. To obtain also uniqueness, we have to impose additional assumptions on the data, which we present in Subsubsection~\ref{subsec:ass-unique}
\subsubsection{Assumptions on data for existence}
\label{subsec:ass-exist}
In this subsubsection, we present the conditions we impose on the data to obtain existence of a solution for the system of equations \eqref{eq:HJBs2}.
\begin{assumption}[Final cost]
	\label{as:uT} $u_T \in C^{2+\alpha}([0,\infty)) \cap W^3_q([0,\infty))$ for some $\alpha \in (0,1)$, 
	$q \in (3/2,2)$, 
	and $u_T \geq 0$ with $u_T(0)=0$. 
	The constant $c_0 = -\min\{\min u_T',0\} > -\infty$ will appear in several estimates.
\end{assumption}
\begin{assumption}[Initial distribution]
	\label{as:m0} $m_0$ is a smooth probability density: $m_0 \geq 0$, $\int_0^\infty m_0(x) \dif x = 1$, and $m_0 \in W^{2}_2([0,\infty))$.
\end{assumption}
\begin{assumption}[L\'evy measure]
	\label{as:levy} $F(dz)$  is a L\'evy measure (i.e.\ a positive measure with~$\int_{\bb R} |z|^2\wedge1 \, F(dz)<\infty$) 
	which 
	{\color{red} satisfies for $s:=\nicefrac{1}{2}$}
	 that 
	\begin{equation}\label{eq:kernel}
	\int_{|z|\geq 1} |z|^s \,F(dz)<\infty.
	\end{equation}
\end{assumption}
\begin{assumption}[Nonlocal coupling]
	\label{as:coupling}
	The function  $(\phi,m) \mapsto {\color{red}\pi}[\phi,m]$ is a real-valued function on $L^\infty([0,\infty)) \times L_{\geq 0}^1([0,\infty))$ such that either
	\begin{enumerate}
		\item ${\color{red}\pi}[\phi,m]$ is bounded for bounded values of $\int_0^\infty |\phi(x)| m(x)\,dx$, {\color{red}meaning that $\forall M>0$  $\exists N>0$ such that $\int_0^\infty |\phi(x)| m(x)\,dx\leq N$ implies that $\big|\pi[\phi,m]\big| \leq M$,} or
		\item ${\color{red}\pi}[\phi,m]$ is bounded for bounded values of {\color{red}$\Vert \phi \Vert_{L^\infty}:=\essinf \phi$  and of $\enVert{m}_{L^1}:=\int_0^\infty |m(x)|\, dx$, meaning that 
		$\forall M>0$  $\exists N>0$ such that $\Vert \phi \Vert_{L^\infty} + \enVert{m}_{L^1}\leq N$ implies that $\big|\pi[\phi,m]\big| \leq M$}.
	\end{enumerate}
\end{assumption}
If $\phi,m$ are continuous functions on $[0,T] \times [0,\infty)$ such that $\phi(t,\cdot) \in L^\infty(0,\infty)$ and $m(t,\cdot) \in L_{\geq 0}^1(0,\infty)$ for every $t \in [0,T]$, then we can define ${\color{red}\pi}[\phi,m]$ as a function on $[0,T]$ by ${\color{red}\pi}[\phi,m](t) = \int_0^\infty \phi(t,x)m(t,x)\dif x$.
\begin{assumption}[Nonlocal coupling, cont.]
	\label{as:coupling continuity}
	Let $\phi_1,\phi_2 \in L^\infty(0,\infty)$ and $m_1,m_2 \in L_{\geq 0}^1(0,\infty)$.
	Then for some constant $C > 0$ depending only on $\enVert{\phi_i}_\infty$ and $\enVert{m_i}_1$ for $i=1,2$, we have
	\begin{multline} \label{eq:[phi,m] Holder1}
	\abs{{\color{red}\pi}[\phi_1,m_1]-{\color{red}\pi}[\phi_2,m_2]} \leq C\left\{\int_0^\infty \del{\abs{\phi_1(x)}+\abs{\phi_2(x)}}\abs{m_1(x)-m_2(x)}\dif x \right.\\
	+ \left.\int_0^\infty \abs{\phi_1(x)-\phi_2(x)}\del{m_1(x)+m_2(x)}\dif x + \abs{\eta_1-\eta_2}\right\},
	\end{multline}
{\color{red}where $\eta_i:=\int_0^\infty m_i(x)\,dx$, $i=1,2$}.
\end{assumption}
\begin{assumption}[Hamiltonian]
	\label{as:H} $H \equiv H(t,\xi,\upsilon,\eta)\colon {\color{red}[0,T]} \times \bb{R} \times \bb{R} \times [0,\infty) \to [0,\infty)$  satisfies
	\begin{enumerate}
		\item $H$ is differentiable in $\xi$, and $D_\xi H$ is locally Lipschitz in both variables $\xi$ and $\upsilon$;
		\item $H$ and $-D_\xi H$ are both decreasing in $\xi$ (in particular, $H$ is convex in $\xi$);
		\item {\color{red}for every $t,\xi$, $H(t,\xi,\upsilon,\eta)$ and $D_\xi H(t,\xi,\upsilon,\eta)$ are bounded for bounded values of $\upsilon$ and $\eta$, 
		meaning that $\forall t\in [0,T]$ $\forall \xi\in \mathbb{R}$:  $\forall M>0$  $\exists N>0$ such that
		$|v| + \eta \leq N$ implies $H(t,\xi,\upsilon,\eta) + \big|D_\xi H(t,\xi,\upsilon,\eta)\big| \leq M$;}
		\item $D_{\xi\xi} H(t,\xi,\upsilon,\eta)$ is locally bounded;
		\item $H$ is locally 
		Lipschitz continuous in $\eta$;
		\item $H(T,\xi,\upsilon,\eta)$ depends only on $\xi$, and
		 $H(T,\xi,\upsilon,\eta) \equiv H_T(\xi)$ 
		is a decreasing convex function.
	\end{enumerate}
\end{assumption}
Examples of Hamiltonian satisfying Assumption \ref{as:H} are given below in Section \ref{sec:H examples}.

Our last assumption, in order to obtain smooth solutions up to the final time $T$, is known as a \emph{compatibility condition of order 1}, 
cf.~\cite[Section IV.5]{ladyzhenskaia1968linear}.
Briefly, if $u$ is a classical solution to the first equation in \eqref{eq:HJBs2} that is smooth up to the boundary, then $u_t(t,0) = 0$ for all $t < T$, so the remaining terms in the PDE must also sum to zero at $x = 0$.
Taking $t \to T$ we deduce from the PDE a constraint on the terminal condition $u_T$, which is stated as follows.	
\begin{assumption}
	The function $u_T$ satisfies
	\begin{equation} \label{eq:comp cond}
	\frac{1}{2}\sigma^2 u_T''(0) - ru_T(0) + H_T\del{u_T'(0)} + \s{I}[u_T](0) = 0.
	\end{equation}
\end{assumption}
\subsubsection{Assumptions on data for uniqueness} \label{subsec:ass-unique}
To obtain also uniqueness of the solution of the system of equations \eqref{eq:HJBs2}, we need to impose additional assumptions on the data, which we present in this subsubsection.

In mean field game theory, the standard approach to proving uniqueness is referred to as the Lasry-Lions monotonicity argument, following \cite[Theorem 2.4]{lasry07}.
In the present context, the main idea is as follows.
Let $(u_i,m_i), i=1,2$ be two solutions.
Let $w = u_1 - u_2$ and $\mu = m_1 - m_2$.
Then $w,\mu$ satisfy
\begin{equation} \label{eq:u differences}
w_t+\frac{1}{2}\sigma^2 w_{xx}-rw+{\s{I}}[w]= H(t,(u_1)_x,{\color{red}\pi}[(u_1)_x,m_1],\eta_1)- H(t,(u_2)_x,{\color{red}\pi}[(u_2)_x,m_2],\eta_2)
\end{equation}
and
\begin{equation}\label{eq:m differences}
\mu_t-\frac{1}{2}\sigma^2 \mu_{xx}+\del{D_\xi H(t,u_{2x},{\color{red}\pi}[u_{2x},m_2],\eta_2)m_2 - D_\xi H(t,u_{1x},{\color{red}\pi}[u_{1x},m_1],\eta_1)m_1}_x-{\s{I}}^*[\mu]=0
\end{equation}
with zero initial/boundary conditions.
Using $\mu$ as a test function in \eqref{eq:u differences} and using $w$ as a test function in \eqref{eq:m differences}, subtracting, and integrating by parts, we obtain
{\small
	\begin{equation} \label{eq:energy differences}
	\begin{split}
	0 &= \! \int_0^T \!\!\int_0^\infty \! \! e^{-rt}\del{H(t,u_{2x},{\color{red}\pi}[u_{2x},m_2],\eta_2)-H(t,u_{1x},{\color{red}\pi}[u_{1x},m_1],\eta_1)  -D_\xi H(t,u_{1x},{\color{red}\pi}[u_{1x},m_1],\eta_1) (u_{2x}-u_{1x})}m_1 \dif x \dif t\\
	& \ \ + \int_0^T \! \! \int_0^\infty \! \! e^{-rt}\del{H(t,u_{1x},{\color{red}\pi}[u_{1x},m_1],\eta_1)\!-\!H(t,u_{2x},{\color{red}\pi}[u_{2x},m_2],\eta_2)  -D_\xi H(t,u_{2x},{\color{red}\pi}[u_{2x},m_2],\eta_2) (u_{1x}\!-\!u_{2x})}m_2 \dif x \dif t.
	\end{split}
	\end{equation}}
If $H$ did not depend on $m$, then the strict convexity of $H$ would imply the positivity of each integrand on the right-hand side.
Then uniqueness would be almost immediate (note, however, that we do not have an additional monotone coupling term as in the case of \cite{lasry07}).
Here the path is not so straightforward, due to the additional dependence on nonlocal term.

Our uniqueness result will rely on two assumptions, which we may informally refer to as \emph{uniform convexity} and \emph{smallness}.
Uniform convexity of the Hamiltonian in $\xi$ refers to the existence and strict positivity of $D_{\xi\xi} H(t,u_x,{\color{red}\pi}[u_x,m],\eta)$ for all solutions, which allows us to deduce from \eqref{eq:energy differences} an estimate for $\int \abs{u_{1x}-u_{2x}}^2 (m_1+m_2)$.
Smallness refers to the continuous dependence on variables ${\color{red}\pi}[u_x,m],\eta$ and is stated precisely in Assumption \ref{as:small Lipschitz}.
%
%
%
\begin{assumption}
	\label{as:smoothness}
	For any solution $(u,m)$ of \eqref{eq:HJBs2}, $D_{\xi\xi} H(t,u_x,{\color{red}\pi}[u_x,m],\eta)$
	has a positive lower bound depending only on the data.
\end{assumption}
We refer to Subsection~\ref{sec:H examples} for examples of classes of Hamiltonians in the Bertrand and Cournot competitions that satisfy Assumption~\ref{as:smoothness}.
\begin{remark}
	Assumption \ref{as:smoothness}
	can be exploited to derive additional regularity of $m$ in the System \eqref{eq:HJBs2}.
	Indeed, in typical situations, when in addition to Assumption \ref{as:smoothness} we also have that
	$D_{\xi\xi} H(t,\xi,\upsilon,\eta)$ is a continuous function,  one can see that then  $D_{\xi\xi}H(t,u_x,{\color{red}\pi}[u_x,m],\eta)$ 
	is also \emph{H\"older} continuous.
	Then since \eqref{eq:FP expanded} holds, we can apply e.g.~Lemma \ref{lem:existence Holder} to prove that $m \in C^{1+\alpha/2,2+\alpha}$.
	However, this is merely incidental; our main concern is proving uniqueness.
\end{remark}
%
%
Moreover, our smallness assumption takes the following form:
\begin{assumption}
	\label{as:small Lipschitz}
	$H$ and its first derivative are locally Lipschitz in all variables, and there exists a fixed constant $\varepsilon > 0$ such that
	\begin{equation} \label{eq:epsilon Lipschitz}
	\begin{split}  
	&\abs{D_\xi H(t,\xi_1,\upsilon_1,\eta_1) - D_\xi H(t,\xi_2,\upsilon_2,\eta_2)} 
	\leq C(R)\del{\abs{\xi_1-\xi_2} + \varepsilon\abs{\upsilon_1-\upsilon_2} + \varepsilon\abs{\eta_1-\eta_2}},\\
	&\abs{D_\upsilon H(t,\xi_1,\upsilon_1,\eta_1) - D_\upsilon H(t,\xi_2,\upsilon_2,\eta_2)} 
	\leq C(R)\varepsilon\del{\abs{\xi_1-\xi_2} + \abs{\upsilon_1-\upsilon_2} + \abs{\eta_1-\eta_2}},\\
	&\abs{D_\eta H(t,\xi_1,\upsilon_1,\eta_1) - D_\eta H(t,\xi_2,\upsilon_2,\eta_2)} 
	\leq C(R)\varepsilon\del{\abs{\xi_1-\xi_2} + \abs{\upsilon_1-\upsilon_2} + \abs{\eta_1-\eta_2}},\\
	&\abs{D_\upsilon H(t,\xi_1,\upsilon_1,\eta_1)}, \abs{D_\eta H(t,\xi_1,\upsilon_1,\eta_1)} \leq C(R)\varepsilon
	\end{split}
	\end{equation}
	for all $\abs{\xi_1},\abs{\xi_2}, {\color{red}\abs{v_1}, \abs{v_2}}, 
	\abs{\eta_1},\abs{\eta_2} \leq R$.
\end{assumption}
\begin{remark}\label{rem:as:small-epsilon}
	To obtain uniqueness of the system of equation \eqref{eq:HJBs2}, we need to guarantee that Assumption~\ref{as:small Lipschitz} holds true with respect to a small enough $\varepsilon$, see Theorem~\ref{thm:main-result}.
\end{remark}
\begin{remark}\label{rem:ass-small-example}
	Assumption \ref{as:small Lipschitz} can be verified for all of the examples given in Section \ref{sec:H examples}, where $\varepsilon > 0$ is precisely the parameter appearing in Equation \eqref{eq:coeffs}.
\end{remark}
Now we are able to state our main result of this work.
\subsection{Main Result}\label{subsec:main-result}
In this section we provide in Theorem~\ref{thm:main-result} our main result of this paper which provides existence and uniqueness of a (classical) solution for the system of equations~\eqref{eq:HJBs2}. To that end, let us first define what we call a solution for \eqref{eq:HJBs2}.
\begin{definition}[Solution to the PDE system] \label{def:solution}		
	We say that $(u,m)$ is a (classical) solution to \eqref{eq:HJBs2} provided that
	\begin{enumerate}
		\item $u \in C^{1+\alpha/2,2+\alpha}([0,T] \times [0,\infty))$  for some $\alpha \in (0,1)$,
		\item $m \in C([0,T];L^1([0,\infty))) \cap W_p^{1,2}([0,T] \times [0,\infty))$ for some $p > 3$,
		\item $m\geq 0$ on $[0,T]\times [0,\infty)$, 
		\item $(u,m)$ satisfy system \eqref{eq:HJBs2} pointwise {\color{magenta} a.e.}.
	\end{enumerate}
\end{definition}
{\color{red}
\begin{remark}\label{rem:ConseqAssSolution}
Note that any solution $(u,m)$ to \eqref{eq:HJBs2} satisfies that
 $[0,T] \ni t \mapsto \eta(t):=\int_0^\infty m(t,x)\,dx \in [0,1]$ is decreasing with  $\eta(0)=1$ (see Lemma~\ref{lem:FP basics}), and that 
 $u_x \in W_p^{1,2}([0,T] \times [0,\infty))$ for some $p > 3$ (see Lemma~\ref{lem:uxW21p}).
\end{remark}
}
\begin{theorem}\label{thm:main-result}
	Consider the system of equations \eqref{eq:HJBs2} and let the data satisfy the assumptions presented in Subsubsection~\ref{subsec:ass-exist}. Then the following holds true.
	\begin{enumerate}
		\item \label{thm:main-result-ex} 
		There exists a {\color{magenta}(classical)} solution $(u,m)$ for \eqref{eq:HJBs2}.
		\item \label{thm:main-result-uni} There exists $\varepsilon_0>0$ (depending on the data) such that if in addition 
		the data also satisfies  the assumptions introduced in Subsubsection~\ref{subsec:ass-unique} with the $\varepsilon>0$  appearing in Assumption~\ref{as:small Lipschitz} satisfying $\varepsilon\leq \varepsilon_0$, then the {\color{magenta}(classical)} solution $(u,m)$ is unique.
	\end{enumerate}
\end{theorem}
%
%
%
%
\begin{remark}[A remark on the difficulty of this model] \label{rem:difficulty}
	In System \eqref{eq:HJBs2}, the Hamilton-Jacobi equation features a Hamiltonian with unspecified growth in the gradient variable; indeed, many of our examples 
	are superquadratic or even exponential; we refer to Subsection~\ref{sec:H examples}.
	It is known that generally superquadratic Hamiltonians can cause two phenomena prohibiting the existence of global-in-time solutions: ``gradient blow-up" and ``loss of boundary conditions".
	See \cite{porretta2017lossOfBC,porretta2018blow} and references therein.
	Consider the simple equation
	\begin{equation} \label{eq:heat exp}
	u_t = u_{xx} + e^{-u_x}
	\end{equation}
	with, say, Dirichlet boundary conditions at $x = 0$.
	This problem has been analyzed in \cite{zhang2010rate} (see also \cite{zhang2013gradient}), where a general criterion on the initial condition is given under which gradient blow-up eventually occurs.
	On the other hand, suppose the data were constructed favorably so that an a priori bound $u_x \geq -C$ could be established.
	Then \eqref{eq:heat exp} could be treated as a linear parabolic equation with a source that has an a priori bound.
	By a standard bootstrapping method, the global solutions would be classical.
	
	In what follows, we will take precisely this approach to the Hamilton-Jacobi equation; hence we can actually hope for classical solutions to the system (for nice enough initial data).
	However, there is still the problem of the \emph{coupling}.
	Instead of merely $e^{-u_x}$, here we have roughly speaking to deal with the following formula 
	\begin{equation*} \label{eq:HG formula}
	\begin{cases}
	H(t,u_x,m)
	=e^{(-u_x+
		{\color{red}\pi}[mu_x])}
	\\\\
	{\color{red}\pi}[mu_x]=\int_0^{\infty}u_x(t,x)m(t,x)\dif x
	\end{cases}
	\end{equation*}
	Thus even if $u_x$ is bounded from below, that does not tell us how to handle $\int u_x m$.
	For this we will use a duality trick deriving from mean field games.
	See Section \ref{sec:energy estimates}.
\end{remark}
\subsection{Examples of Hamiltonians} \label{sec:H examples}

The purpose of this section is to prove the existence of a large class of Hamiltonians satisfying Assumption \ref{as:H}.
We will do this by proposing certain abstract conditions on the demand (or inverse demand) function, which happen to be satisfied by a family of particular examples already proposed in the literature.
We consider separately the cases of Bertrand and Cournot competition. 
Moreover, we provide for Bertrand and Cournot competitions classes of examples for which the uniform convexity 
condition presented in Assumption~\ref{as:smoothness} is satisfied.

\subsubsection{Bertrand competition} \label{sec:Bertrand}
The following assumptions on the demand schedule \eqref{eq:demand schedule} are made in the case of Bertrand type competition.
Before proceeding to state the assumptions, we note that the demand schedule $\s{D}^\eta(p,\bar{p})$ is supposed to depend on $\eta$, the proportion of players still playing the game, and $\bar{p}$, the average price.
For technical reasons, we will actually be more interested in the \emph{aggregate} price $\pi = \eta\bar{p}$.
As part of our first assumption, we will rewrite $\s{D}^{\eta}(p,\bar{p})$ as $\s{D}(\eta;p,\pi)$.
A second technical point is that $\s{D}$ will also depend on time, in such a way that at the final time $t = T$, the demand does not depend on the market but only on $p$.
\begin{assumption}
	\label{as:demand}
	$\s{D}^{\eta}(p,\bar{p})$ can be written as $\s{D}(t,\eta;p,\pi)$ where $\pi = \eta \bar{p}$ and $\s{D}$ is twice continuously differentiable in all variables with
	\begin{equation} \label{eq:D dec in p}
	D_p \s{D}(t,\eta;p,\pi) < 0 \quad \ \forall \eta, p, \pi \geq 0.
	\end{equation}
	Moreover, $\s{D}(T,\eta;p,\pi) = \s{D}_T(p)$, i.e.\  $\s{D}(T,\eta;p,\pi)$ only depends on $p$ at time $T$.
\end{assumption}
\begin{assumption}
	\label{as:g1}
	The function
	\begin{equation} \label{eq:maximizer function}
	g(t,\eta;p,\pi) := \frac{\s{D}(t,\eta;p,\pi)}{D_p \s{D}(t,\eta;p,\pi)} + p
	\end{equation}
	is strictly increasing in $p$, and in fact there exists a constant $\delta_0 > 0$ such that
	\begin{equation} \label{eq:g increasing}
	D_p g(t,\eta;p,\pi) \geq \delta_0 \quad \ \forall p,\pi.
	\end{equation}
	One may also express this condition as
	\begin{equation} \label{eq:relative prudence}
	\frac{\s{D}D_{pp}\s{D}}{(D_p \s{D})^2} \leq 2 - \delta_0.
	\end{equation} 
\end{assumption}
The quantity $\frac{\s{D}D_{pp}\s{D}}{(D_p \s{D})^2}$ is known in economics as \emph{relative prudence}.
Condition \eqref{eq:relative prudence} is exactly that which appears in \cite[Proposition 2.5]{harris2010games} in order to guarantee existence of a unique Nash equilibrium.

{\color{red} The following assumption is {\color{blue} rather} technical, but we point out that it is satisfied by the examples considered in this section; we refer to Subsection~\ref{subsubsec:examples}.}
{\color{blue} We can interpret it as measuring the degree to which the demand is more sensitive to changes in an individual sale price than to the market price, which reflects the imperfect substitutability of goods.}
\begin{assumption}
	\label{as:g2}
	$g(t,\eta;p,\bar p)$ is strictly decreasing in $\bar p$, and there exists a constant $\beta \in (0,1)$ such that
	\begin{equation} \label{eq:pbar vs p}
	\abs{D_{\pi} g(t,\eta;p,\pi)} \leq \eta^{-1}\beta D_p g(t,\eta;p,\pi) \quad \forall \eta,p,\pi > 0.
	\end{equation}
	This could also be expressed as
	\begin{equation}
	\abs{D_p \s{D} D_{\pi}\s{D} - \s{D}D_{p\pi}\s{D}} \leq \eta^{-1}\beta\del{2(D_p \s{D})^2 - \s{D}D_{pp}\s{D}}.
	\end{equation}
	Moreover, $g(t,\eta;p,\bar p)$ is differentiable in $\eta$ with locally bounded derivative.
\end{assumption}

\begin{lemma} \label{lem:pstar}
{\color{red} Let the Assumptions~\ref{as:demand}, \ref{as:g1}, and~\ref{as:g2} hold, and}	define $F(t,\eta;\xi,\pi,p):= \s{D}(t,\eta;p,\pi)(p-\xi)$ for $\xi \in \bb{R},\pi \geq 0,$ and $p \geq 0$.
	Then there exists a function $p^* \equiv p^*(\eta;\xi,\pi) \geq 0$ defined on $[0,\infty) \times \bb{R} \times [0,\infty)$ such that
	\begin{equation}
	F(t,\eta;\xi,\pi,p^*) = \sup_{p \geq 0} F(t,\eta;\xi,\pi,p) =: H(t,\xi,\pi,\eta).
	\end{equation}
	The function $p^*$ is continuously differentiable in $(\xi,\pi)$ on the set where $\xi > g(t,\eta;0,\pi)$, and it is identically equal to zero on the set where $\xi \leq \frac{\s{D}(t,\eta;0,\pi)}{D_p \s{D}(t,\eta;0,\pi)}$.
	The derivatives $D_\xi p^*$ and $D_{\pi} p^*$ are defined on the set where $\xi \neq g(t,\eta;0,\pi)$ and satisfy the following bounds:
	\begin{equation}
	\label{eq:derivatives of pstar}
	0 \leq D_\xi p^* \leq \delta_0^{-1}, 
	\quad \quad \quad 
	\abs{D_{\pi} p^*} \leq \beta\eta^{-1}.
	\end{equation}
	Thus $p^*$ is Lipschitz continuous for bounded values of $\eta$.
	
	In addition, $H(t,\xi,\pi,\eta)$ is differentiable in $\xi$ on the set where $\xi \neq g(t,\eta;0,\pi)$, and $D_\xi H(t,\xi,\pi,\eta) = -\s{D}(t,\eta;p^*,\pi)$. 
	Thus, $H(t,\xi,\pi,\eta)$ is decreasing and convex with respect to $\xi$, and $D_\xi H(t,\xi,\pi)$ is locally Lipschitz with respect to both $\xi$ and $\pi$.
\end{lemma}

\begin{proof}
	We compute, using \eqref{eq:maximizer function},
	\begin{equation}
	D_p F(t,\xi,\pi,p) = -D_p \s{D}(t,\eta;p,\pi)\del{\xi - g(t,\eta;p,\pi)}.
	\end{equation}
	We then use \eqref{eq:g increasing} and the implicit function theorem.
	In particular, $p^*$ is defined by
	\begin{equation}
	g(t,\eta;p^*,\pi) = \xi
	\end{equation}
	on the set where $\xi > g(t,\eta;0,\pi)$; otherwise we set $p^* = 0$.
	By implicit differentiation,
	\begin{equation}
	D_\xi p^* = \frac{1}{D_p g(t,\eta;p^*,\pi)}, 
	\quad \quad 
	D_{\pi} p^* = - \frac{D_{\pi} g(t,\eta;p^*,\pi)}{D_p g(t,\eta;p^*,\pi)},
	\quad \quad
	D_{\eta} p^* = - \frac{D_{\eta} g(t,\eta;p^*,\pi)}{D_p g(t,\eta;p^*,\pi)}
	\end{equation}
	on the set where $\xi > g(t,\eta;0,\pi)$, while $D_\xi p^* = D_{\pi} p^* = D_\eta p^* = 0$ on the set where $\xi < g(t,\eta;0,\pi)$.
	The estimates \eqref{eq:derivatives of pstar} as well as the local boundedness of $D_\eta p^*$ follow by \eqref{eq:g increasing} and \eqref{eq:pbar vs p}.
	
	We differentiate $H(t,\xi,\pi) = F(t,\xi,\pi,p^*)$ by implicit differentiation to deduce the remaining claims.
\end{proof}

\begin{lemma}
	\label{lem:market clearing}
	{\color{red} Let the Assumptions~\ref{as:demand}, \ref{as:g1}, and~\ref{as:g2} hold, and let $p^* \equiv p^*(\eta;\xi,\pi)$ be the function defined in Lemma~\ref{lem:pstar}.}
	Given $\phi \in L^\infty(0,\infty)$ and $m \in L^1(0,\infty)$ such that $m \geq 0$, set $\eta = \int_0^\infty m(x)\dif x.$
	Then there exists a unique $\pi \geq 0$ such that
	\begin{equation} \label{eq:market clearing}
	\pi = \int_0^\infty p^*(\eta;\phi(x),\pi)m(x)\dif x.
	\end{equation}
	The solution satisfies the upper bound
	\begin{equation} \label{eq:pbar upper bound}
	\pi \leq \frac{1}{(1-\beta)\delta_0\eta}\int_0^\infty \abs{\phi(x)}m(x)\dif x + \frac{1}{1-\beta}p^*(\eta;0,0)
	\leq \frac{\enVert{\phi}_\infty}{(1-\beta)\delta_0} + \frac{1}{1-\beta}p^*(\eta;0,0).
	\end{equation}
\end{lemma}

\begin{proof}
	For $\pi \geq 0$ define
	\begin{equation}
	f(\pi) = \int_0^\infty p^*(\eta;\phi(x),\pi)m(x)\dif x.
	\end{equation}
	Our goal is to prove that $f$ has a unique fixed point.
	Observe that by \eqref{eq:derivatives of pstar} we have
	\begin{equation} \label{eq:contraction}
	\abs{f'(\pi)} \leq \int_0^\infty \abs{D_{\pi} p^*(\eta;\phi(x),\pi)}m(x)\dif x \leq \beta < 1
	\end{equation}
	and
	\begin{equation} \label{eq:f upper bound}
	f(\pi) \leq \frac{1}{\delta_0\eta}\int_0^\infty \abs{\phi(x)}m(x)\dif x + \beta \pi + p^*(\eta;0,0).
	\end{equation}
	By the contraction mapping theorem and \eqref{eq:contraction}, $f$ has a fixed point, for which \eqref{eq:pbar upper bound} must hold by \eqref{eq:f upper bound}.
\end{proof}

By Lemma \ref{lem:market clearing}, we can define ${\color{red}\pi}[\phi,m]:={\color{red}\int_0^\infty p^*(\eta;\phi(x),\pi)m(x)\dif x}$ {\color{red}as introduced in \eqref{eq:market clearing}} for any $(\phi,m) \in L^\infty(0,\infty) \times L^1(0,\infty)$, and this defines a function satisfying Assumption \ref{as:coupling}.
Moreover, $H(t,\xi,{\color{red}\pi}[\phi,m],\eta)$ satisfies Assumption \ref{as:H} by Lemma \ref{lem:pstar}.
We now need to show that ${\color{red}\pi}[\phi,m]$ satisfies Assumption \ref{as:coupling continuity}.

\begin{lemma} \label{lem:p bar Holder}
	{\color{red} Let the Assumptions~\ref{as:demand}, \ref{as:g1}, and~\ref{as:g2} hold, and for  any $(\phi,m) \in L^\infty(0,\infty) \times L^1(0,\infty)$ let
		$\pi[\phi,m]:=\int_0^\infty p^*(\eta;\phi(x),\pi)m(x)\dif x$ be defined as in \eqref{eq:market clearing}.}
		 Then the map $(\phi,m) \mapsto {\color{red}\pi}[\phi,m]$ satisfies Assumption \ref{as:coupling continuity}.
\end{lemma}

\begin{proof}
	Let $q$ be as in Assumption \ref{as:uT} and let $q' = q/(q-1)$ be its H\"older conjugate.
	Suppose $\phi_1,\phi_2 \in L^\infty(0,\infty) \cap L^q(0,\infty)$ and $m_1,m_2 \in L_{\geq 0}^1(0,\infty) \cap L^{q'}(0,\infty)$ (so that $\phi_1,\phi_2 \in L^{q'}$ and $m_1,m_2 \in L^q$ as well), and set $\eta_i = \int_0^\infty m_i(x)\dif x$ for $i=1,2$.
	Let $\pi_i = {\color{red}\pi}[\phi_i,m_i], i =1,2$.
	We have
	\begin{equation*}
	\begin{split}
	\abs{\pi_1-\pi_2}
	&= \abs{\int_0^\infty p^*(\eta_1;\phi_1(x),\pi_1)m_1(x)\dif x - \int_0^\infty p^*(\eta_2;\phi_2(x),\pi_2)m_2(x)\dif x}\\
	&\leq
	\int_0^\infty \abs{p^*(\eta_1;\phi_1(x),\pi_1)-p^*(\eta_1;\phi_2(x),\pi_1)}m_1(x)\dif x\\
	& \quad
	+ \int_0^\infty \abs{p^*(\eta_1;\phi_2(x),\pi_1)-p^*(\eta_1;\phi_2(x),\pi_2)}m_1(x)\dif x\\
	& \quad
	+ \int_0^\infty \abs{p^*(\eta_1;\phi_2(x),\pi_2)-p^*(\eta_2;\phi_2(x),\pi_2)}m_1(x)\dif x\\
	& \quad
	+ \abs{\int_0^\infty p^*(\eta_2;\phi_2(x),\pi_2)\del{m_1(x)-m_2(x)}\dif x}\\
	&\leq 
	\delta_0^{-1}\int_0^\infty \abs{\phi_1(x)-\phi_2(x)}m_1(x)\dif x
	+ \beta\abs{\pi_1-\pi_2}
	+ L\abs{\eta_1-\eta_2}\eta_1\\
	& \quad
	+ \delta_0^{-1}\int_0^\infty \abs{\phi_2(x)}\abs{m_1(x)-m_2(x)}\dif x
	+ p^*(\eta_2;0,\pi_2)\abs{\eta_1-\eta_2}
	\end{split}
	\end{equation*}
	where $L$ is a Lipschitz constant for $p^*$ in the $\eta$ variable determined by an upper bound on $\enVert{\phi_2}_\infty,\pi_2,\eta_2$ and $\eta_1$.
	Subtract $\beta\abs{\pi_1-\pi_2}$ from both sides and divide by $1-\beta$ to get
	\begin{equation*}
	\begin{split}
	\abs{\pi_1-\pi_2}
	&\leq 
	\frac{1}{\delta_0 (1-\beta)}\delta_0^{-1}\int_0^\infty \abs{\phi_1(x)-\phi_2(x)}m_1(x)\dif x
	\\
	& \quad
	+ \frac{1}{\delta_0 (1-\beta)}\int_0^\infty \abs{\phi_2(x)}\abs{m_1(x)-m_2(x)}\dif x
	+ \frac{p^*(\eta_2;0,\pi_2)+L\eta_1}{1-\beta}\abs{\eta_1-\eta_2}.
	\end{split}
	\end{equation*}
	Thus \eqref{eq:[phi,m] Holder1} holds.
\end{proof}

\subsubsection{Examples}\label{subsubsec:examples}
As a motivating example, we take a demand function of the following form 
\begin{equation} \label{eq:Dexp}
\s{D}^{\eta}(p,\bar{p})=e^{(a(t,\eta)-p+c(t,\eta)\bar{p})}
\end{equation}
with $a(t,\eta)+c(t,\eta)=1$. The coefficients depending on $\eta$ are defined by
\begin{equation} \label{eq:coeffs}
a(t,\eta)=\frac{1}{1+\varepsilon(t)\eta},
\quad \quad \quad 
c(\eta)=\frac{\varepsilon(t)\eta}{1+\varepsilon(t)\eta}
\end{equation}
for a smooth, non-negative function $\varepsilon(t)$ such that $\varepsilon(t) \leq \varepsilon(0)$ and $\varepsilon(T) = 0$.
It is then straightforward to verify Assumptions \ref{as:demand}, \ref{as:g1}, and \ref{as:g2}; in this case $\delta_0 = 1$ and $\beta \in (0,1)$ can be arbitrarily small.
Indeed, in this case we simply have
\begin{equation} \label{eq:gexp}
g(t,\eta;p,\pi) = p-1,
\quad \quad 
D_p g = 1,
\quad \quad 
D_{\pi} g = 0,
\quad \quad
D_\eta g = 0.
\end{equation}

Another example would be a power function of the form
\begin{equation} \label{eq:Dpowerlaw}
\s{D}^{\eta}(p,\bar{p})=\del{1+\frac{a(t,\eta)-p+c(t,\eta)\bar{p}}{\rho}}_+^\rho
\end{equation}
for some $\rho \geq 2$.
(One can also take $\rho > 1$, although this is not twice continuously differentiable; to fully satisfy the assumptions, it would suffice to take a smooth approximation with the same behavior as this function.)
Here we have
\begin{equation*}
\begin{split}
&g(t,\eta;p,\pi) = p - \del{1+ \frac{a+\frac{c}{\eta}\pi-p}{\rho}}_+,
\quad \quad
D_p g(t,\eta;p,\pi) = 1 + \frac{1}{\rho}\chi_{\rho+ a+\frac{c}{\eta}\pi > p}, \\
& 
D_{\pi} g(t,\eta;p,\pi) = -\frac{c}{\rho\eta}\chi_{\rho+ a+\frac{c}{\eta}\pi > p}, 
\quad \quad
D_\eta g(t,\eta;p,\pi) = \frac{\varepsilon + \varepsilon^2\pi}{\del{1+\varepsilon \eta}^2}\pi\chi_{\rho+ a+\frac{c}{\eta}\pi > p}.
\end{split}
\end{equation*}
In particular, we observe that Assumptions \ref{as:demand}, \ref{as:g1}, and \ref{as:g2} are satisfied with $\delta_0 = 1$ and $\beta = \max c(t,\eta) = \frac{\varepsilon(0)}{1+\varepsilon(0)} \in (0,1)$.

\begin{remark}
	Although the case of exponential growth as in \eqref{eq:Dexp} would seem to represent the most challenging example of the Hamiltonians given above, based on the growth rate with respect to the gradient variable, nevertheless in other respects this case is much simpler than for power laws like \eqref{eq:Dpowerlaw} because of \eqref{eq:gexp}.
\end{remark}
Now,  let us  provide a sufficient condition for the Hamiltonian in the Bertrand competition to satisfy the uniform convexity 
condition introduced in Assumption~\ref{as:smoothness}.
\begin{lemma}[Uniform convexity for Bertrand competition] \label{lem:smoothnessBertrand}
	Suppose $H$ is defined as in Section \ref{sec:Bertrand} with demand schedule $\s{D}^\eta(p,\bar p)$ given by \eqref{eq:Dexp}.
	Let $c_0 = -\min\{\min u_T',0\}$ and assume $c_0 < 1$.
	Then Assumption \ref{as:smoothness} is satisfied.
\end{lemma}

\begin{proof}
	Let $(u,m)$ be any solution of \eqref{eq:HJBs2} and recall from Lemma \ref{lem:ux lower bound} that $u_x \geq -c_0$.
	Recall that the function $p^*(\eta;\xi,\pi)$ defined in Lemma \ref{lem:pstar} is given by solving $g(t,\eta;p^*,\pi) = \xi$ on the set $\xi > g(t,\eta;0,\pi)$, on which set it is continuously differentiable.
	By \eqref{eq:gexp}, $g(t,\eta;p,\pi) = p-1$ and thus it follows that $g(t,\eta;0,\pi) = -1 < -c_0 \leq \min u_x$.
	We therefore conclude that $p^*$ is given precisely by $p^*(\eta;\xi,\pi) = \xi+1$.
	By implicit differentiation, we have
	\begin{equation}
	D_{\xi\xi} H(t,\xi,\pi,\eta) = -D_\xi \s{D}(t,\eta;\xi+1,\pi) = \exp\del{a(t,\eta)-p+\frac{c(t,\eta)}{\eta}\pi},
	\end{equation}
	which is continuous and even smooth in all variables.
	The claim follows.
\end{proof}
\begin{remark}
	Under the hypotheses of Lemma \ref{lem:smoothnessBertrand}, we can find an explicit formula for the Hamiltonian.
	Using the formula for $\bar p$ given in Lemma \ref{lem:market clearing}, we deduce after some calculation that in fact
	\begin{equation}
	H(t,u_x,{\color{red}\pi}[u_x,m],\eta)(t,x) = \exp\del{\frac{c(t,\eta)}{\eta}\int_0^\infty u_x(t,y)m(t,y)\dif y - u_x(t,x)}.
	\end{equation}
\end{remark}
\subsubsection{Cournot competition} \label{sec:Cournot}

The following assumptions are made in the case of Cournot competition, in which each firm uses quantity produced as a control.
Thus we consider price as a function of quantity, as in \eqref{eq:inverse demand schedule}.

\begin{assumption}
	\label{as:P}
	$P = P(t,q,Q) > 0$ is twice continuously differentiable on $[0,\infty) \times [0,\infty) \times [0,\infty)$.
\end{assumption}

\begin{assumption}
	\label{as:h}
	For all $q,Q \geq 0$ define
	\begin{equation} \label{eq:h def}
	h(t,q,Q) = qD_q P(t,q,Q) + P(t,q,Q).
	\end{equation}
	We assume that $h$ is strictly decreasing in both variables, i.e.\ for all $q > 0$, $Q > 0$ we assume that
	\begin{equation}
	q D_{qq} P(t,q,Q) + 2D_q P(t,q,Q) < 0, \quad \quad qD_{qQ}P(t,q,Q) + D_Q P(t,q,Q) < 0. 
	\end{equation}
	Moreover, we assume that both $\frac{1}{D_q h(t,q,Q)}$ and $\frac{D_Qh(t,q,Q)}{D_q h(t,q,Q)}$ remain bounded as $(q,Q) \to (0,0)$.
\end{assumption}

Now for the rest of this subsubsection, define $F(t,\xi,Q,q) := q(P(t,q,Q)-\xi)$.
\begin{lemma}[Unique optimal quantity] \label{lem:qstar}
	There exists a unique $q^* = q^*(t,\xi,Q) \geq 0$ such that $F(t,\xi,Q,q^*)$ $= \sup_{q \geq 0} F(t,\xi,Q,q) =: H(t,\xi,Q)$ for all $t \geq 0, \xi \in \bb{R}, Q \geq 0$.
	The function $q^*$ is continuously differentiable in $(\xi,Q)$ on the set where $\xi < P(t,0,Q)$ and identically zero on the set where $\xi \geq P(t,0,Q)$.
	The derivatives $D_\xi q^*$ and $D_Q q^*$ are well-defined and non-positive on the set $\xi \neq P(0,Q)$.
	In particular $q^*$ is decreasing and locally Lipschitz in both $\xi$ and $Q$, and in fact $q^*$ is $\beta$-Lipschitz with respect to the $Q$ variable.
	Moreover, $H(t,\xi,Q)$ is differentiable in $\xi$ on the set $\xi \neq P(t,0,Q)$ with $D_\xi H(t,\xi,Q) = -q^*$; in particular, $H$ is convex and decreasing in $\xi$, and $D_\xi H$ is locally Lipschitz with respect to both variables.
\end{lemma}

\begin{proof}
	We compute that $D_q F(t,\xi,Q,q) = h(t,q,Q) - \xi$, using \eqref{eq:h def}.
	The proof now follows from using $D_q h < 0, D_Q h < 0$, and the implicit function theorem, cf.~the proof of Lemma \ref{lem:pstar}.
	By implicit differentiation we can deduce
	\begin{equation}
	D_\xi q^* = \frac{1}{D_q h(t,q^*,Q)} < 0,
	\quad \quad \quad
	D_Q q^* = - \frac{D_Q h(t,q^*,Q)}{D_q h(t,q^*,Q)} < 0
	\end{equation}
	in the region where $\xi < h(t,0,Q)$, while $D_\xi q^* = D_Q q^* = 0$ in the region $\xi > h(t,0,Q)$.
	From Assumption \ref{as:h} we deduce that $q^*(t,\xi,Q)$ is locally Lipschitz on $\intco{0,\infty} \times \bb{R} \times [0,\infty)$.
	We omit the remaining details, which follow the same way as in Lemma \ref{lem:pstar}.
\end{proof}

\begin{lemma}[Unique aggregate quantity] \label{lem:market clearing Cournot}
	Fix $t \in \intco{0,\infty}$.
	Let $\phi \in L^\infty(0,\infty)$ and $m \in L^1(0,\infty), m \geq 0$.
	Then there exists a unique $Q \geq 0$ such that
	\begin{equation}\label{eq:market clearing Cournot}
	Q = \int_0^\infty q^*(t,\phi(x),Q)m(x) \dif x.
	\end{equation}
	Moreover, $Q$ satisfies the a priori estimate
	\begin{equation} \label{eq:Q upper bound}
	Q \leq q^*(t,\inf \phi,0)\enVert{m}_{L^1}.
	\end{equation}
\end{lemma}

\begin{proof}
	Let $f(Q) = Q - \int_0^\infty q^*(t,\phi(x),Q)m(x) \dif x$.
	We claim that $f(Q) = 0$ for a unique $Q \geq 0$.
	Note that $f(0) \leq 0$ because $q^* \geq 0$.
	Since by Lemma \ref{lem:qstar} we have $D_Q q^* \leq 0$, it follows that $f' \geq 1$.
	The claim follows, and we deduce \eqref{eq:market clearing Cournot}.
	Estimate \eqref{eq:Q upper bound} follows from the fact that $q^*$ is decreasing in both $\xi$ and $Q$.
\end{proof}

By Lemma \ref{lem:market clearing Cournot} we can define a map $(\phi,m) \mapsto {\color{red}\pi}[\phi,m]$ from $L^\infty(0,\infty) \times L_{\geq 0}^1(0,\infty)$ by setting ${\color{red}\pi}[\phi,m] = Q$,
where $Q$ is the unique solution of \eqref{eq:market clearing Cournot}.

\begin{lemma} \label{lem:Q Holder}
	The function $(\phi,m) \mapsto {\color{red}\pi}[\phi,m]$ satisfies Assumptions \ref{as:coupling} and \ref{as:coupling continuity}.
\end{lemma}

\begin{proof}
	Assumption \ref{as:coupling} is satisfied because of \eqref{eq:Q upper bound}.
	
	Now we verify Assumption \ref{as:coupling continuity}.
	Let $\phi_1,\phi_2 \in L^\infty(0,\infty) \cap L^q(0,\infty)$ and $m_1,m_2 \in L_{\geq 0}^1(0,\infty) \cap L^{q'}(0,\infty)$, so that $\phi_1,\phi_2 \in L^{q'}$ and $m_1,m_2 \in L^q$ as well, and set $\eta_i = \int_0^\infty m_i(x)\dif x, i=1,2$.
	Let $Q_1 = {\color{red}\pi}[\phi_1,m_1] = \int_0^\infty q^*(\phi_1(x),Q_1)m_1(x)\dif x$ 
	and 
	$Q_2 = {\color{red}\pi}[\phi_2,m_2] = \int_0^\infty q^*(\phi_2(x),Q_2)m_2(x)\dif x$.
	Without loss of generality assume $Q_1 \geq Q_2$.
	Since $q^*$ is decreasing in the $Q$ variable, it follows that
	\begin{equation*}
	\begin{split}
	Q_1-Q_2 &\leq \int_0^\infty q^*(\phi_1(x),Q_1)m_1(x)\dif x - \int_0^\infty q^*(\phi_2(x),Q_1)m_2(x)\dif x\\
	&= 
	\int_0^\infty \del{q^*(\phi_1(x),Q_1)-q^*(0,Q_1)}\del{m_1(x)-m_2(x)}\dif x
	+ q^*(0,Q_1)\del{\eta_1-\eta_2}\\
	& \quad 
	+ \int_0^\infty \del{q^*(\phi_1(x),Q_1)-q^*(\phi_2(x),Q_1)}m_2(x)\dif x,
	\end{split}
	\end{equation*}
	and since $q^*$ is locally Lipschitz, there exist a constant $C$ depending on $\enVert{\phi_i}_\infty, i=1,2$ such that
	\begin{equation}
	Q_1 - Q_2 \leq C\int_0^\infty \abs{\phi_1(x)}\abs{m_1(x)-m_2(x)}\dif x
	+ C\abs{\eta_1-\eta_2} + C\int_0^\infty \abs{\phi_1(x)-\phi_2(x)}m_2(x)\dif x.
	\end{equation}
	Thus \eqref{eq:[phi,m] Holder1} holds.	
\end{proof}
Finally, we note that $H(t,\xi,{\color{red}\pi}[\phi,m],\eta):= H(t,\xi,Q)$ satisfies Assumption \ref{as:H}; note that the H\"older continuity with respect to $\eta$ is trivial, since there is no explicit dependence on $\eta$.

\subsubsection{Examples}

Consider a smooth, non-negative function $\varepsilon(t)$ such that $\varepsilon(t) \leq \varepsilon(0)$ and $\varepsilon(T) = 0$.
Let $\rho \in [0,1]$.
For $\rho > 0$ set
\begin{equation} \label{eq:inverse demand power law}
P(t,q,Q) = \frac{1}{\rho}\del{1-(q+\varepsilon(t) Q)^\rho}
\end{equation}
and for $\rho = 0$
\begin{equation} \label{eq:inverse demand logarithmic law}
P(t,q,Q) = -\log(q+\varepsilon(t) Q).
\end{equation}
This is exactly the set of examples proposed by Chan and Sircar in \cite{chan2017fracking}, except that $\varepsilon(t)$ is assumed to be time-dependent.
It is straightforward to check that
\begin{equation}
D_q h(t,q,Q) = -\frac{(1+\rho)q + 2\varepsilon(t) Q}{(q+\varepsilon(t) Q)^{2-\rho}},
\quad
D_Q h(t,q,Q) = -\frac{\rho q + \varepsilon(t) Q}{(q+\varepsilon(t) Q)^{2-\rho}},
\end{equation}
and hence Assumptions \ref{as:P} and \ref{as:h} are satisfied.

We finish this subsection by providing a sufficient condition for the Hamiltonian in the Cournot competition to satisfy the uniform convexity 
condition introduced in Assumption~\ref{as:smoothness}.
\begin{lemma}[Uniform convexity for Cournot competition] \label{lem:smoothnessCournot}
	Let $H$ be defined as in Section \ref{sec:Cournot} with inverse demand schedule $P(t,q,Q)$ given either by \eqref{eq:inverse demand logarithmic law} for $\rho = 0$ or by \eqref{eq:inverse demand power law} when $\rho > 0$.
	If $\rho \geq 0$ and $\varepsilon > 0$ are sufficiently small, then Assumption \ref{as:smoothness} is satisfied.
\end{lemma}

\begin{proof}
	Let $(u,m)$ be any solution of \eqref{eq:HJBs2}.
	We have  by Lemma~\ref{lem:qstar} that
	\begin{equation*}
	H(t,u_x,{\color{red}\pi}[u_x,m],\eta) = q^*(t,u_x,Q)\big(P(q^*(t,u_x,Q),Q)-u_x\big) \quad \mbox{ and } \quad -D_\xi H(t,u_x,{\color{red}\pi}[u_x,m],\eta) = q^*(t,u_x,Q),
\end{equation*}
 where $Q$ is given by applying Lemma  \ref{lem:market clearing Cournot} with $\phi = u_x$, and $q^*(\xi,Q)$ is decreasing in both variables.
	Let $c_0 = -\min\{\min u_T',0\}$ and note that from Lemma \ref{lem:ux lower bound} we have $u_x \geq -c_0$.
	It follows that $q^*(t,u_x,Q) \leq q^*(t,-c_0,0)$.
	By direct computation using formula \eqref{eq:inverse demand logarithmic law} or \eqref{eq:inverse demand power law}, we get
	\begin{equation}
	q^*(t,-c_0,0) =
	\begin{cases}
	\del{\frac{1+\rho c_0}{1+\rho}}^{1/\rho} & \text{if}~\rho > 0,\\
	e^{c_0-1} & \text{if}~\rho = 0,
	\end{cases}
	\quad \quad \
	P(q^*(t,-c_0,0),0)+c_0 =
	\begin{cases}
	\frac{1+\rho c_0}{1+\rho} & \text{if}~\rho > 0,\\
	1 & \text{if}~\rho = 0,
	\end{cases}
	\end{equation}
	and note that in all cases $q^*(-c_0,0) \leq e^{c_0}$.
	Thus, using the fact that $P(t,q,Q) \leq P(t,q,0)$, we have
	\begin{align*}
	H(t,u_x,{\color{red}\pi}[u_x,m],\eta) &= q^*(t,u_x,Q)(P(q^*(t,u_x,Q),Q)-u_x)\\
	&\leq q^*(t,u_x,Q)(P(q^*(t,u_x,Q),0)+c_0)\\
	&\leq q^*(t,-c_0,0)(P(q^*(t,-c_0,0),0)+c_0)\\
	&\leq e^{c_0}(1+c_0),
	\end{align*}
	where in the second to last line we have used the definition of $q^*(t,\xi,Q)$ as the maximizing argument of $q(P(t,q,Q)-\xi)$.
	It follows that the $L^\infty$ norms of $H(t,u_x,{\color{red}\pi}[u_x,m],\eta)$ and $D_\xi H(t,u_x,{\color{red}\pi}[u_x,m],\eta)$ depend only on $c_0$, even as $\rho$ and $\varepsilon(t)$ vary!
	
	Therefore, Lemmas \ref{lem:uW21p-local} and \ref{lem:uxW21p} are applicable to show that $\max u_x \leq \xi_0$, where $\xi_0$ is a constant depending only on $\sigma,r,\enVert{u_T}_{C^2},\enVert{u_T}_{W^{2,1}_q}$ and $c_0$.
	If we can ensure, by fixing $\rho,\varepsilon$ small enough, that $\xi_0 < P(t,0,Q(t))$ for all $t$, then it follows from Lemma \ref{lem:qstar} that $D_{\xi\xi} H(t,u_x,{\color{red}\pi}[u_x,m],\eta) = -D_\xi q^*(t,u_x,Q)$ is a continuous, positive function.
	It suffices to impose the following condition on $\rho$ and $\varepsilon$:
	\begin{equation}
	\label{eq:rho epsilon small}
	\xi_0 < P(0,0,e^{c_0}) =
	\begin{cases}
	\frac{1}{\rho}\del{1-\varepsilon(0)^\rho e^{\rho c_0}} &\text{if}~\rho > 0,\\
	-\log \varepsilon(0) - c_0 &\text{if}~\rho = 0.
	\end{cases}
	\end{equation}
	Using the fact that $Q(t) \leq q^*(t,-c_0,0) \leq e^{c_0}$ by Lemma \ref{lem:market clearing Cournot}, \eqref{eq:rho epsilon small} implies the desired constraint ensuring that $D_{\xi\xi} H(t,u_x,{\color{red}\pi}[u_x,m],\eta)$ exists and is continuous.
	We then obtain the following lower bound:
	\begin{equation}
	D_{\xi\xi} H(t,u_x,{\color{red}\pi}[u_x,m],\eta) \geq \min\{-D_\xi q^*(t,\xi,Q) : t \geq 0, -c_0 \leq \xi \leq \xi_0, 0 \leq Q \leq q^*(-c_0,0)\} > 0.
	\end{equation}
\end{proof}
%
%
%
%
\subsection{Linear demand schedules and variational structure}

In this subsection we make some remarks on the special case when the demand schedule is linear, i.e.
\begin{equation} \label{eq:linear demand}
\s{D}^\eta(p,\bar p) = a(t,\eta) + c(t,\eta)\bar p - p
\end{equation}
with $a,c$ given as in \eqref{eq:coeffs}.
If we consider $m(x)$ to be the current distribution of states, then we have that $\eta = \int_0^\infty m(x)\dif x$ and $\bar p = \int_0^\infty p(x)m(x)\dif x$.
If we change variables, setting $q = \s{D}^\eta(p,\bar p)$ and $Q = \int_0^\infty q(x)m(x)\dif x$, then a simple calculation yields
\begin{equation} \label{eq:linear inverse demand}
p = P(t,q,Q) = 1 - (q+\varepsilon Q),
\end{equation}	
a linear inverse demand schedule.
Thus, Cournot and Bertrand competition turn out to be merely inverses of one another, i.e.~they are equivalent (see \cite{chan2015bertrand}).
In fact, when there is no constraint on the quantity or price, i.e.~both negative prices and production are allowed, the first-order optimality condition can be used to derive an explicit form:
\begin{equation} \label{eq:quadratic H explicit}
\begin{split}
H(t,u_x,{\color{red}\pi}[u_x,m],\eta) 
= q^*(t,u_x,Q)^2 
= \frac{1}{4}\del{1-\varepsilon Q - u_x}^2 
&= \frac{1}{4}\del{a(\eta) + c(\eta)\bar p - u_x}^2 
\\
&= \frac{1}{4}\del{\frac{2}{2+\varepsilon \eta} + \frac{\varepsilon}{2+\varepsilon \eta}\int_0^\infty u_x m \dif y - u_x}^2,
\end{split}
\end{equation}
see \cite{chan2015bertrand,graber2018existence}.

In this case, uniqueness of solutions is \emph{unconditional}, as pointed out in \cite{graber2018variational}; see also \cite{graber2018commodities}.
We rely on a simple expression for the Hamiltonian and its derivative:
\begin{equation} \label{eq:quadratic H}
H(t,u_x,{\color{red}\pi}[u_x,m],\eta) = q^*(t,u_x,Q)^2, \ D_\xi H(t,u_x,{\color{red}\pi}[u_x,m],\eta) = -q^*(t,u_x,Q).
\end{equation}
Then for two solutions $(u_1,m_1)$ and $(u_2,m_2)$, \eqref{eq:energy differences} becomes
\begin{equation} \label{eq:energy differences quadratic}
\begin{split}
0 
&= \int_0^T \int_0^\infty e^{-rt}\del{q^*(t,u_{2x},Q_2)^2-q^*(t,u_{1x},Q_1)^2 +q^*(t,u_{1x},Q_1) (u_{2x}-u_{1x})}m_1(t,x) \dif x \dif t\\
& \quad 
+ \int_0^T \int_0^\infty e^{-rt}\del{q^*(t,u_{1x},Q_1)^2-q^*(t,u_{2x},Q_2)^2  +q^*(t,u_{2x},Q_2) (u_{1x}-u_{2x})}m_2(t,x) \dif x \dif t,
\end{split}
\end{equation}
and the proof is completed essentially by completing the square; see \cite[Theorem 2]{graber2018variational}.
Note that the argument holds even when constraints are imposed; indeed, if we replace the demand schedule \eqref{eq:linear demand} by its positive part, or equivalently consider the inverse demand \eqref{eq:linear inverse demand} only over $q \geq 0$, then we still obtain formula \eqref{eq:quadratic H}, although the formula is less explicit; then the argument holds as in step 3 of the proof of \cite[Theorem 1.5]{graber2018commodities}.

As a final remark, we note that this is related to the \emph{variational structure} of the problem--see \cite{graber2018variational,bonnans2019schauder} for further discussion.
Indeed, in this case we can write System \eqref{eq:HJBs2} as the first-order optimality condition for the problem of maximizing
\begin{equation} \label{control of FP}
\begin{split}
J(m,q) 
&= \int_0^T \int_0^L e^{-rt} \left(q^2(t,x)  - q(t,x)\right)m(t,x)\dif x \dif t
\\
& \quad
+ \frac{\varepsilon}{2}\int_0^T e^{-rt}\left(\int_0^L  q(t,y)m(t,y)\dif y\right)^2\dif t - \int_0^L e^{-rT} u_T(x)m(T,x)\dif x
\end{split}
\end{equation}
for $(m,q)$ in the class $\s{K}$, defined as the set of all pairs such that  $m \in L_{\geq 0}^1([0,T] \times [0,L])$, $q \in L^2([0,T] \times [0,L])$, and $m$ is a weak solution to the Fokker-Planck equation
\begin{equation}
\label{fokker planck condition}
m_t -  \frac{\sigma^2}{2}m_{xx} - \s{I}^*[m] - (qm)_x = 0, \qquad \qquad m(0) = m_0,
\end{equation}
equipped with Dirichlet boundary conditions.
We refer to \cite[Proposition 3]{graber2018variational} for the proof, and we note that it works if we include the constraint $q \geq 0$ in the definition of $\s{K}$.
%
%
%
%
\section{Auxiliary lemmas} \label{sec:auxiliary}
First, we want to establish some preliminary facts which are apparently well-known but are not easily found in the literature.
Our first two lemmas essentially say that $\s{I}$ is relatively compact with respect to the Laplacian.
\begin{lemma} \label{lem:compactness estimate}
	Let $u \in W^{1,2}_p([0,T] \times \bb{R})$ for some $p \in (1,\infty)$. Then
	for all $\varepsilon > 0$ there exists $C_\varepsilon > 0$ such that
	\begin{equation}
	\enVert{\s{I}[u]}_p \leq \varepsilon\enVert{u_{xx}}_p + C_\varepsilon\big(\enVert{u}_p+ \enVert{u_x}_p\big).
	\end{equation}
\end{lemma}

\begin{proof}
	See, e.g., in \cite[Lemma~2.2.1]{garroni2002second}.
\end{proof}

\begin{lemma} \label{lem:compactness estimate Holder}
	Let $u \in C^{1+\alpha/2,2+\alpha}([0,T] \times \bb{R})$ for some $\alpha \in (0,1)$.
	Then for all $\varepsilon > 0$ there exists $C_\varepsilon > 0$ such that
	\begin{equation} \label{eq:Holder compactness}
	\enVert{\s{I}[u]}_{C^{\alpha/2,\alpha}} \leq \varepsilon\enVert{u_{xx}}_{C^{\alpha/2,\alpha}} + C_\varepsilon\big(\enVert{u}_{C^{\alpha/2,\alpha}}+\enVert{u_x}_{C^{\alpha/2,\alpha}}\big).
	\end{equation}
\end{lemma}

\begin{proof}
	Decompose $\s{I}[u]$ into two terms 
	by writing for any $0<r<1$ that
	\begin{equation} \label{Irho-1}
	\begin{split}
	\s{I}[u](t,x)
	&=\int_{\bb{R} \setminus B_r(0)}[u(t,x+z)- u(t,x)-u_x(t,x)\,z\mathbf{1}_{\{|z|\leq 1\}}]\,F(dz)\\
	& \quad 
	+
	\int_{B_r(0)}[u(t,x+z)- u(t,x)-u_x(t,x)\,z]\,F(dz)
	\end{split}
	\end{equation}
	For the first term, set $C_r:=\int_{\bb{R} \setminus B_r(0)} 2\,F(dz)<~\infty$  to obtain that
	\begin{equation}\label{eq:big stuff}
	\abs{\int_{\bb{R} \setminus B_r(0)}[u(t,x+z)- u(t,x)-u_x(t,x)\,z\mathbf{1}_{\{|z|\leq 1\}}]\, F(dz)}
	\leq  C_r\big(\enVert{u}_\infty +\enVert{u_x}_\infty\big)<\infty.
	\end{equation}
	For the second term, we have for every  $0<r<1$ that
	\begin{equation} \label{eq:small stuff}
	\begin{split}
	\abs{\int_{B_r(0)}[u(t,x+z)- u(t,x)-u_x(t,x)\,z\mathbf{1}_{\{|z|\leq 1\}}]\,F(dz)}
	&= \abs{\int_{-r}^r\int_0^1 \int_{0}^1 u_{xx}(t,x+\tau \sigma z)\sigma z^2\dif \tau \dif \sigma \,F(dz)}\\
	&\leq \enVert{u_{xx}}_\infty  \int_{-r}^r\int_0^1 \int_{0}^1 \sigma\dif \tau \dif \sigma \,z^2\,F(dz)\\
	&=: \enVert{u_{xx}}_\infty c_r,
	\end{split}
	\end{equation}
	where $c_r$ is a finite constant 
	which satisfies that $c_r \to 0$ as $r\to 0$.
	Therefore, we obtain from \eqref{eq:big stuff} and \eqref{eq:small stuff} that
	\begin{equation} \label{eq:bound on Irho}
	\enVert{\s{I}[u]}_\infty \leq  C_r\big(\enVert{u}_\infty +\enVert{u_x}_\infty\big) 
	+ c_r\enVert{u_{xx}}_\infty.
	\end{equation}
	Now, observe that for any  $(\tau,y)\in[0,T]\times \mathbb R\setminus(0,0)$ we have 
	for $w(t,x) := u(t+\tau,x+y)-u(t,x)$ that $\s{I}[u](t+\tau,x+y) - \s{I}[u](t,x) = \s{I}[w](t,x)$.
	Therefore, we can apply \eqref{eq:bound on Irho} with respect to $w$ to get
	\begin{equation*}
	\begin{split} 
	\sup_{t,x}\, \abs{\s{I}[u](t+\tau,x+y) - \s{I}[u](t,x)} 
	&\leq
	C_r\big(\enVert{w}_\infty +\enVert{w_x}_\infty\big) 
	+ c_r\enVert{w_{xx}}_\infty 
	\\
	&\leq \del{C_r\big(\enVert{u}_{C^{\alpha/2,\alpha}} +\enVert{u_x}_{C^{\alpha/2,\alpha}}\big) 
		+ c_r\enVert{u_{xx}}_{C^{\alpha/2,\alpha}}}\del{\abs{\tau}^{\alpha/2} + \abs{y}^\alpha},
	\end{split}
	\end{equation*}
	which implies that
	\begin{equation} \label{eq:Holder bound on Irho}
	\intcc{\s{I}[u]}_{C^{\alpha/2,\alpha}} \leq 	C_r\big(\enVert{u}_{C^{\alpha/2,\alpha}} +\enVert{u_x}_{C^{\alpha/2,\alpha}}\big) 
	+ c_r\enVert{u_{xx}}_{C^{\alpha/2,\alpha}}.
	\end{equation}
	Adding \eqref{eq:bound on Irho} to \eqref{eq:Holder bound on Irho} and letting $r$ be small enough, we obtain \eqref{eq:Holder compactness}.
\end{proof}

Our next lemma is {\color{red} a particular case of the well-known} 
maximum principle for a nonlocal equation involving~$\s{I}$; {\color{red} see, e.g., also \cite{barles2008second}.}
\begin{lemma}
	[Maximum principle] \label{lem:max principle}
	Let $w \in C^{1,2}([0,T] \times [0,\infty))$ be a classical solution of
	\begin{equation}
	-w_t - \frac{1}{2}\sigma^2 w_{xx} - \s{I}[w] \leq 0
	\end{equation}
	satisfying  $w(t,x) \leq 0 \ \forall x \leq 0 \ \forall t \in [0,T]$, and  $w(T,x) \leq 0 \ \forall x \in \bb{R}$.
	Then we have $\sup w \leq 0$.
\end{lemma}

\begin{proof}
	We imitate the proof of \cite[Proposition 2.1]{souplet2006global}. 
	Let $\phi(x) = (1+x^2)^{\frac{s}{2}}$, where $s \in (0,1)$ is defined in Assumption~\ref{as:levy}, and note that 
	$$\phi'(x) = sx(1+x^2)^{\frac{s}{2} - 1}, \qquad  \phi''(x) = s(1+x^2)^{\frac{s}{2} - 1} + s(s-2)x^2(1+x^2)^{\frac{s}{2} - 2}.$$
	Then $\abs{\phi(x)} \leq 2\abs{x}^s$ for $\abs{x} \geq 1$ while $\abs{\phi''(x)} \leq s(2-s)$ for all $x$. 
	Therefore, we have that
	\begin{equation*}
	\begin{split}
	\s{I}[\phi](x) 
	&=
	\int_{|z|\leq 1} \intcc{\phi(x+z)-\phi(x)-\phi'(x)z}\,F(dz)
	+ \int_{|z|>1} \intcc{\phi(x+z)-\phi(x)}\,F(dz)\\
	&\leq
	\int_{|z|\leq 1} \int_0^1 \int_0^1 \phi''(x+\tau \sigma z) z^2 \sigma \, d\tau\,d\sigma \,F(dz)
	+ 2\int_{|z|>1} \intcc{|x+z|^s + |x|^s}\,F(dz)\\
	&\leq s(2-s)\int_{|z|\leq 1}  z^2\,F(dz) 
	+ 2\int_{|z|>1} \abs{z}^s F(dz) + 4F\big(\{|z|>1\})|x|^s<\infty.
	\end{split}
	\end{equation*}
	Now let $\psi:\mathbb{R} \to \mathbb{R}_+$ be a smooth function for which there exists a constant $0<C_\psi<\infty$ such that 
	\begin{equation*}
	\begin{split}
	|y|^s\leq 1 + \psi(y) \quad \forall y \in \mathbb{R}, 
	\qquad \qquad
	\frac{1}{2}\sigma^2 T\sup_{y \in \mathbb{R}}|\psi''(y)|\leq  C_\psi.
	\end{split}
	\end{equation*}
	Then, we have for all $x \in \mathbb{R}$ that
	\begin{equation}
	\label{eq:I-psi-estim}
	\begin{split}
	\s{I}[\phi](x)
	&\leq 
	s(2-s)\int_{|z|\leq 1}  z^2\,F(dz) 
	+ 2\int_{|z|>1} \abs{z}^s F(dz) + 4F\big(\{|z|>1\}) (1+\psi(x))\\
	&:= C_1 + C_2\, \psi(x).
	\end{split}
	\end{equation}
	Now, define the smooth function $M:\mathbb{R} \to \mathbb{R}_+$ by
	\begin{equation}
	\label{eq:def-fct-M}
	M(x):= \frac{1}{2}\sigma^2s(2-s) + C_1 + C_2 \psi(x) + C_2C_\psi +1
	\end{equation}
	and  define $v(t,x) = w(t,x) - \varepsilon(M(x)(T-t) + \phi(x))$ for some fixed $\varepsilon>0$. Then, using that  $|\phi''|\leq s(2-s)$, \eqref{eq:I-psi-estim}, and \eqref{eq:def-fct-M} ensures 
	for any $(t,x)$ that
	\begin{equation}
	\label{eq:v inequality}
	\begin{split}
	-v_t(t,x) - \frac{1}{2}\sigma^2 v_{xx}(t,x) - \s{I}[v](t,x) &\leq \varepsilon\del{\frac{1}{2}\sigma^2\phi''(x) + \s{I}[\phi] -M(x) +\frac{1}{2}\sigma^2 (T-t)M''(x)}\\
	&\leq -\varepsilon.
	\end{split}
	\end{equation}
	Next, since $w$ is bounded and $\lim_{x \to \infty} \phi(x) = +\infty$, it follows that $v(t,x) \to -\infty$ as $x \to \infty$, uniformly in $t \in [0,T]$.
	Thus there exists a point $(t_0,x_0) \in [0,T] \times [0,\infty)$ such that $v(t_0,x_0) = \sup v$.
	Assume by contradiction that $\sup v > 0$.
	Then by the boundary conditions on $w$, we see that $t_0 < T$ and $x_0 > 0$. This implies 
	that $v_t(t_0,x_0) \leq 0$, $v_x(t_0,x_0) = 0$, $v_{xx}(t_0,x_0) \leq 0$, as well as $\s{I}[v](t_0,x_0) \leq 0$, since
	\begin{equation}
	\s{I}[v](t_0,x_0) = \int_{-\infty}^\infty \del{v(t_0,x_0+z)-v(t_0,x_0)-v_x(t_0,x_0)\,z\mathbf{1}_{\{|z|\leq 1\}}}\,F(dz) \leq 0,
	\end{equation}
	using the fact that $v(t_0,x_0) \geq v(t_0,x_0+ z)$.
	We thus obtain a contradiction with \eqref{eq:v inequality} and hence conclude that  $\sup v \leq 0$. Letting $\varepsilon \to 0$ then ensures that also $\sup w \leq 0$.
\end{proof}

Before extending the maximum principle for a nonlocal equation to solutions in $W^{1,2}_p$,
we provide some asymptotic properties needed to prove the convergence of certain integrals.

\begin{lemma} \label{lem:goes to zero}
	Let $u \in C^{\alpha/2,\alpha}([0,T] \times [0,\infty)) \cap L^q([0,T] \times [0,\infty))$ for some $\alpha \in (0,1), q > 1$.
	Then $\sup_{t\in [0,T]} \abs{u(x,t)} \to 0$ as $x \to \infty$.
\end{lemma}

\begin{proof}
	Let $M > 0$, $r \in (0,1)$. Observe that
	\begin{equation*}
	\abs{u(t,x)} \leq \frac{1}{2r^2}\int_{[0,T] \cap B_r(t)} \int_{x-r}^{x+r} \abs{u(\tau,\xi)}\dif \xi \dif \tau + \intcc{u}_{C^{\alpha/2,\alpha}}\del{r^{\alpha/2} + r^\alpha}
	\end{equation*}
	By applying H\"older's inequality, we see for any $x\geq M+1$ that 
	\begin{equation*}
	\begin{split}
	\sup_{t \in [0,T]} \abs{u(t,x)} 
	&\leq 
	\frac{(2rT)^{(q-1)/q}}{2r^2}\del{\int_0^T \int_M^\infty \abs{u(\tau,\xi)}^q\dif \xi \dif \tau}^{1/q} + \intcc{u}_{C^{\alpha/2,\alpha}}\del{r^{\alpha/2} + r^\alpha}\\
	& \quad \to
	\intcc{u}_{C^{\alpha/2,\alpha}}\del{r^{\alpha/2} + r^\alpha}
	\end{split}
	\end{equation*}
	as $M \to \infty$.
	Letting $r \to 0$ we obtain the desired result.
\end{proof}

\begin{lemma}
	[Maximum principle for weak solutions] \label{lem:max principle weak}
	Let $p \in (1,\infty)$, $C > 0$, and
	suppose that $w \in W^{1,2}_p([0,T] \times [0,\infty))$ is a solution {\color{red}in the $L^p$-sense} of
	\begin{equation}\label{eq:max-prin-weak}
	-w_t - \frac{1}{2}\sigma^2 w_{xx} - \s{I}[w] \leq C\abs{w_x}, \ w(t,x) \leq 0 \ \forall x \leq 0 \ \forall t \in [0,T], \ w(T,x) \leq 0 \ \forall x \in \bb{R}.
	\end{equation}
	Then $\sup w \leq 0$.
\end{lemma}

\begin{proof}
	We use a multiplier method to show that $w\leq0$. 
	Fix $t\in [0,T]$,
	set $q := \max\{2,p\}$, and define
	$w_+ := \max\{w,0\}$.
	For any $R>0$,  we multiply \eqref{eq:max-prin-weak} by $w_+^{q-1}$, then
	integrate over $[t,T] \times [0,R]$ and apply the Peter--Paul inequality to get that
	\begin{equation}
	\begin{split}
	&\frac{1}{q}\int_0^R w_+^q(t,x)\dif x 
	-\int_t^T \frac{\sigma^2}{2}w_{x}(s,R)w_+^{q-1}(s,R)\,ds
	+ \frac{\sigma^2}{2}(q-1)\int_t^T \int_0^R (w_+)_x^2(s,x) w_+^{q-2}(s,x) \dif x \dif s\\
	&\quad - \int_t^T \int_0^R \s{I}[w](s,x) w_+^{q-1}(s,x) \dif x \dif s \\
	&\leq C \int_t^T \int_0^R \abs{w_{x}}(s,x)w_+^{q-1}(s,x) \dif x \dif s\\
	&\leq \frac{\sigma^2}{4}\int_t^T \int_0^R (w_+)_x^2(s,x) w_+^{q-2}(s,x)  \dif x \dif s
	+ \frac{C^2}{\sigma^2}\int_t^T \int_0^R (w_+(t,x))^q \dif x \dif s.
	\end{split}
	\end{equation}
	By Sobolev-type embedding theorems,  see e.g.~\cite[Lemma II.3.4]{ladyzhenskaia1968linear} $w,w_{x} \in L^p \cap C^{\alpha/2,\alpha}$ for some $\alpha \in (0,1)$.
	By Lemma \ref{lem:goes to zero} it follows that $w(t,x) \to 0$ as $x \to \infty$ uniformly in $t$, and thus $w_+(t,R) \to 0$ as $R \to \infty$, uniformly in $t$.
	It also follows that $w_+^{q-1} \in L^{k}$ for any $k \geq p/(q-1)$, in particular for $k = p/(p-1)$; since $\s{I}[w] \in L^p$ by Lemma \ref{lem:compactness estimate}, 
	$\int_0^T \int_0^\infty \s{I}[w]w^{q-1}_+(s,x)\,dx\,ds$ converges absolutely when $T$ goes to infinity.
	By letting $R \to \infty$, we therefore obtain that 
	\begin{equation} \label{eq:w_+ estimate}
	\begin{split}
	&\frac{1}{q}\int_0^\infty w_+^q(t,x)\dif x+ \frac{\sigma^2}{4}\int_t^T \int_0^\infty (w_+)_x^2(s,x) w_+^{q-2}(s,x) \dif x \dif s - \int_t^T \int_0^\infty \s{I}[w](s,x) w_+^{q-1}(s,x) \dif x \dif s \\
	&\leq \frac{C^2}{\sigma^2}\int_t^T \int_0^\infty (w_+(s,x))^q \dif x \dif s.
	\end{split}
	\end{equation}
	Finally, we want to show that $\int_t^T \int_0^\infty \s{I}[w](s,x) w_+^{q-1}(s,x)\leq 0$.
	To that end, first note that
	\begin{equation}
	\label{eq:max-split}
	\begin{split}
	&\int_t^T \int_0^\infty \s{I}[w](s,x)\, w_+^{q-1}(s,x) \dif x \dif s\\
	&= \lim_{r \to 0} \int_t^T \int_0^\infty \int_{|z|\geq r} \del{w(s,x+z)-w(s,x)-w_x(s,x)\,z\mathbf{1}_{\{|z|\leq 1\}}} w_+^{q-1}(s,x) F(\dif z) \dif x \dif s.
	\end{split}
	\end{equation}
	Therefore, it suffices to show for any $0<r<1$ that
	\begin{equation*}
	\int_t^T \int_0^\infty \int_{|z|\geq r} \del{w(s,x+z)-w(s,x)-w_x(s,x)\,z\mathbf{1}_{\{|z|\leq 1\}}} w_+^{q-1}(s,x) F(\dif z) \dif x \dif s\leq 0.
	\end{equation*}
	Fix any $0<r<1$ and observe that
	\begin{equation*}
	\begin{split}
	&	\int_t^T \int_0^\infty \int_{|z|\geq r} \del{w(s,x+z)-w(s,x)-w_x(s,x)\,z\mathbf{1}_{\{|z|\leq 1\}}} w_+^{q-1}(s,x) F(\dif z) \dif x \dif s\\
	&=
	\int_t^T \int_0^\infty \int_{|z|\geq r} \del{w(s,x+z)-w(s,x)} w_+^{q-1}(s,x) F(\dif z) \dif x \dif s\\
	& \quad 
	+\int_t^T \int_0^\infty \int_{r\leq |z| \leq 1}\del{-w_x(s,x)z} w_+^{q-1}(s,x) F(\dif z) \dif x \dif s. 
	\end{split}
	\end{equation*}
	For the first summand, notice that by Fubini's theorem
	\begin{equation}
	\label{eq:max-first-1}
	\begin{split}
	&\int_t^T \int_0^\infty \int_{|z|\geq r} \del{w(s,x+z)-w(s,x)} w_+^{q-1}(s,x) F(\dif z) \dif x \dif s\\
	&\leq \int_t^T \int_{|z|\geq r} \int_0^\infty  |w(s,x+z) w_+^{q-1}(s,x)|  \dif x\, F(\dif z)\, \dif s 
	-\int_t^T  \int_{|z|\geq r} \|w_+(s,\cdot)\|_{L^q(0,\infty)}^q\,F(\dif z)\,\dif s.
	\end{split}
	\end{equation}
	Moreover, by H\"older's inequality, we have that
	\begin{equation}
	\label{eq:max-first-2}
	\begin{split}
	& \int_t^T \int_{|z|\geq r} \int_0^\infty  |w(s,x+z) w_+^{q-1}(s,x)|  \dif x\, F(\dif z)\, \dif s \\
	&\leq \int_t^T \int_{|z|\geq r} \bigg(\int_0^\infty  |w(s,x+z)|^q  \dif x\bigg)^{\nicefrac{1}{q}} \bigg(\int |w_+(s,x)|^{q}\,\dif x\bigg)^{\frac{q-1}{q}}\, F(\dif z)\, \dif s\\
	& \leq \int_t^T  \int_{|z|\geq r} \|w_+(s,\cdot)\|_{L^q(0,\infty)}^q\,F(\dif z)\,\dif s.
	\end{split}
	\end{equation}
	Therefore, we conclude from \eqref{eq:max-first-1} and \eqref{eq:max-first-2} that the first summand satisfies
	\begin{equation}
	\label{eq:max-first}
	\int_t^T \int_0^\infty \int_{|z|\geq r} \del{w(s,x+z)-w(s,x)} w_+^{q-1}(s,x) F(\dif z) \dif x \dif s\leq 0.
	\end{equation}
	For the second summand, as $C_r:=\int_{r\leq |z| \leq 1} z \, F(dz)<\infty$ and using the boundary conditions, we have that
	\begin{equation}
	\label{eq:max-second}
	\begin{split}
	\int_t^T \int_0^\infty \int_{r\leq |z| \leq 1}\del{-w_x(s,x)z} w_+^{q-1}(s,x) F(\dif z) \dif x \dif s
	& = C_r
	\int_t^T \int_0^\infty  -w_x(s,x) w_+^{q-1}(s,x)  \dif x \dif s\\
	&= - \tfrac{C_r}{q}\int_t^T \int_0^\infty  \big(w_+^{q}(s,x)\big)_x  \dif x \dif s\\
	&=0.
	\end{split}
	\end{equation}
	Therefore, we conclude from \eqref{eq:max-split} together with \eqref{eq:max-first} and \eqref{eq:max-second} that indeed
	\begin{equation}
	\int_t^T \int_0^\infty \s{I}[w](s,x)\, w_+^{q-1}(s,x) \dif x \dif s\leq 0.
	\end{equation}
	This and \eqref{eq:w_+ estimate} ensures that 
	\begin{equation} 
	\begin{split}
	\frac{1}{q}\int_0^\infty w_+^q(t,x)\dif x
	&\leq \frac{C^2}{\sigma^2}\int_t^T \int_0^\infty (w_+(s,x))^q \dif x \dif s,
	\end{split}
	\end{equation}
	which by Gronwall's Lemma implies that $\int_0^\infty w_+^q(t,x)\dif x$. As $t\in [0,T]$ was arbitrary and by the continuity $w$, we can thus conclude that
	$w \leq 0$.
\end{proof}

We now state two basic existence results for parabolic equations with a nonlocal term:
\begin{equation}\label{eq:parabolic}
\begin{array}{lc}
u_t+\frac{1}{2}\sigma^2 u_{xx} + {\s{I}}[u] + a(t,x)u_x + b(t,x)u=f,&  0<t<T,\ 0<x<\infty\\
u(T,x)=u_T(x), & 0\leq x< \infty\\
u(t,x)=0, & 0\leq t\leq T,  -\infty< x\leq 0.
\end{array}
\end{equation}
We have stated the problem backwards in time, starting at time $t = T$, only because we will most often apply these results to the Hamilton-Jacobi equation in \eqref{eq:HJBs2}, which is backwards in time; the result is exactly the same if we reverse time and start at time $t = 0$.
Both lemmas essentially follow from the fact that the nonlocal operator is a relatively compact perturbation of the Laplacian, according to \ref{lem:compactness estimate}.	
\begin{lemma} \label{lem:existence Lp}
	Let $p \in (1,\infty)$ and let $f \in L^p([0,T] \times [0,\infty))$, $u_T \in W^{2-2/p}_p([0,\infty))$. Moreover, let 
	$a(t,x) \in L^r_{loc}([0,T]\times [0,\infty))$, $b(t,x)\in L^s_{loc}([0,T]\times [0,\infty))$, where $r > \max\{p,3\}, s > \max\{p,3/2\}$.
	Then there exists a unique weak solution $u \in W^{1,2}_p([0,T]\times [0,\infty))$ satisfying the boundary value problem \eqref{eq:parabolic}.		
	Moreover, the following estimate holds:
	\begin{equation}
	\label{eq:parabolic estimate}
	\enVert{u}_{W^{1,2}_p} \leq C\del{\enVert{f}_p + \enVert{u_T}_{W^{2-2/p}_p}}
	\end{equation}
	where $C = C(\enVert{a}_{L^r_{loc}},\enVert{b}_{L^s_{loc}},T)$.
\end{lemma}

\begin{proof}		
	We follow an argument similar to that of \cite[Proposition 3.11]{cirant2019fractionalmfg}.
	For any small parameter $\tau > 0$ and any
	$m>0$ 
	we set $Q_\tau = (T-\tau,T) \times (0,\infty)$ and
	\begin{equation*}
	X_{k,\tau} := \{u \in W^{1,2}_p(Q_\tau) :  u(T,\cdot) = u_T(\cdot), \ \enVert{u}_{W^{1,2}_p(Q_\tau)} \leq k\}.
	\end{equation*}
	For $z \in W_p^{1,2}(Q_\tau)$  define $Jz$ to be the unique solution $w \in W^{1,2}_p(Q_\tau)$ of
	\begin{equation}\label{eq:parabolic fixed point}
	\begin{array}{lc}
	w_t+\frac{1}{2}\sigma^2 w_{xx} + a(t,x)w_x + b(t,x)w+{\s{I}}[z]=f,&  0<t<T,\ 0<x<\infty\\
	w(T,x)=u_T(x), & 0\leq x< \infty\\
	w(t,0)=0, & 0\leq t\leq T.
	\end{array}
	\end{equation}
	By \cite[Theorem IV.9.1]{ladyzhenskaia1968linear}, $Jz$ is well-defined and satisfies the estimate
	\begin{equation}
	\enVert{Jz}_{W^{1,2}_p(Q_\tau)} \leq C\del{\enVert{f}_{L^p(Q_\tau)} + \enVert{u_T}_{W^{2-2/p}_p} + \enVert{{\s{I}}[z]}_{L^p(Q_\tau)}}.
	\end{equation}
	By Lemma \ref{lem:compactness estimate} we have
	\begin{equation*}
	\enVert{{\s{I}}[z]}_{L^p(Q_\tau)}\leq \delta\enVert{z}_{W^{1,2}_p(Q_\tau)} +  C_\delta \big( \enVert{z}_{L^p(Q_\tau)}+ \enVert{z_x}_{L^p(Q_\tau)}\big).
	\end{equation*}
	On the other hand, by \cite[Lemma II.3.3]{ladyzhenskaia1968linear} we have
	\begin{equation*}
	\enVert{z_x}_p \leq \gamma\enVert{z}_{W^{1,2}_p(Q_\tau)} 
	+ 
	C_\gamma  \enVert{z}_{L^p(Q_\tau)}
	\end{equation*}
	for arbitrary $\gamma > 0$.
	It follows that
	\begin{equation}
	\label{eq:lin1}
	\begin{split}
	\enVert{Jz}_{W^{1,2}_p(Q_\tau)}
	&\leq 
	C\Big(\!\enVert{f}_{L^p(Q_\tau)}
	+ 
	\enVert{u_T}_{W^{2-2/p}_p(Q_\tau)}
	+ 
	(\delta+C_\delta
	\gamma)\enVert{z}_{W^{1,2}_p(Q_\tau)}
	+ C_\delta(1+ C_\gamma) \enVert{z}_{L^p(Q_\tau)}		 	  
	\!\Big).
	\end{split}
	\end{equation}
	Moreover, by integrating
	\begin{equation*}
	z(t,x) = u_T(x) - \int_t^T \partial_t z(s,x)\dif s
	\end{equation*}
	for any $t \in (T-\tau,T)$ and using Jensen's inequality, we have
	\begin{equation*}
	\begin{split}
	\enVert{z}_{L^p(Q_\tau)}^p
	&\leq 
	2^p  \int_{T-\tau}^T\int_0^\infty |u_T(x)|^p \, dx\,dt
	+ 
	2^p \int_{T-\tau}^T\int_0^\infty \bigg|\int_t^T \partial_t z(s,x)\dif s\bigg|^p \, dx\,dt\\
	&\leq 
	2^p \tau \enVert{u_T}^p_{L^p} 
	+ 
	2^p \int_{T-\tau}^T\int_0^\infty (T-t)^{p-1}\int_t^T |\partial_t z(s,x)|^p\dif s \, dx\,dt\\
	&\leq 
	2^p \tau \enVert{u_T}^p_{L^p} 
	+ 
	2^p \tau^{p-1}\int_{T-\tau}^T\int_{T-\tau}^T\int_0^\infty  |\partial_t z(s,x)|^p \, dx\,dt \, \dif s\\
	&\leq
	2^p \tau \enVert{u_T}^p_{L^p} 
	+ 
	2^p \tau^{p} \enVert{z}^p_{W^{1,2}_p(Q_\tau)}
	\end{split}
	\end{equation*}
	and so
	\begin{equation*}
	\begin{split}
	\enVert{z}_{L^p(Q_\tau)} 
	\leq 2\tau^{1/p}\enVert{u_T}_{L^p} + 2\tau \enVert{z}_{W^{1,2}_p(Q_\tau)}.
	\end{split}
	\end{equation*}
	Thus, we see from \eqref{eq:lin1} that
	\begin{equation} \label{eq:contraction1}
	\begin{split}
	&\enVert{Jz}_{W^{1,2}_p(Q_\tau)}\\
	&\leq 
	C\del{\enVert{f}_{L^p(Q_\tau)} + \enVert{u_T}_{W^{2-2/p}_p} + 2C_\delta(1+ C_\gamma)\tau^{1/p}\enVert{u_T}_{L^p} + \big(\delta
		+C_\delta
		\gamma
		+ 2C_\delta(1+ C_\gamma) \tau\big)\enVert{z}_{W^{1,2}_p(Q_\tau)}}.
	\end{split}
	\end{equation}
	Now if $\tilde z$ is another element of $W^{1,2}_p(Q_\tau)$, then $Jz - J\tilde z$ is the solution of \eqref{eq:parabolic fixed point} with $f$ and $u_T$ replaced by zero, hence
	\begin{equation} \label{eq:contraction2}
	\enVert{Jz - J\tilde z}_{W^{1,2}_p(Q_\tau)} \leq C\big(\delta
	+C_\delta
	\gamma
	+ 2C_\delta(1+ C_\gamma) \tau\big)\enVert{z-\tilde z}_{W^{1,2}_p(Q_\tau)}.
	\end{equation}
	By setting $k > C\del{\enVert{f}_{L^p(Q_\tau)} + \enVert{u_T}_{W^{2-2/p}_p}}$, then taking $\delta$ and $\gamma$ small enough, we see that $\tau$ can be chosen small enough so that \eqref{eq:contraction1} and \eqref{eq:contraction2} imply $J$ is a well-defined contraction on $X_{k,\tau}$.
	We apply the contraction mapping theorem to get a fixed point $u = Ju$ in $X_{k,\tau}$.
	Note that there exists a constant $C_{k,\tau}$ such that
	\begin{equation} \label{eq:fixed point estimate}
	\enVert{u}_{W^{1,2}_p(Q_\tau)} \leq C_{k,\tau}\del{\enVert{f}_{L^p(Q_\tau)} + \enVert{u_T}_{W^{2-2/p}_p}}.
	\end{equation}
	
	To complete the proof, we divide $(0,T]$ into subintervals of length no greater than $\tau$, i.e.~$(0,T] = \cup_{i=0}^N (t_i,t_{i+1}]$ where $t_0 = 0$, $t_{N+1} = T$, and $0 < t_{i+1} - t_i < \tau$ for all $i$.
	Set $Q_i = (t_i,t_{i+1}) \times (0,\infty)$ and
	\begin{equation*}
	X_{k,i} := \{u \in W^{1,2}_p(Q_i) : u(T) = u_T, \ \enVert{u}_{W^{1,2}_p(Q_i)} \leq k\}.
	\end{equation*}
	By the argument above, we can find a fixed point $u_N=Ju_N$ in $X_{k,N}$ satisfying
	\begin{equation*}
	\enVert{u_N}_{W^{1,2}_p(Q_N)} \leq C_{k,\tau}\del{\enVert{f}_{L^p(Q_N)} + \enVert{u_T}_{W^{2-2/p}_p}}.
	\end{equation*}
	Now replace $u_T$ with $u_N(t_N,\cdot)$; by the trace theorem \cite[Lemma II.3.4]{ladyzhenskaia1968linear} we have
	\begin{equation*}
	\enVert{u_N(t_N,\cdot)}_{W^{2-2/p}_p} \leq C\enVert{u_{N}}_{W^{1,2}_p(Q_N)}.
	\end{equation*}
	Thus we obtain $u_{N-1} \in X_{k,N-1}$, which is a fixed point of the same operator defined up to time $t_N$ rather than $T$, satisfying
	\begin{equation*}
	\enVert{u_{N-1}}_{W^{1,2}_p(Q_{N-1})} \leq C_{k,\tau}\del{\enVert{f}_{L^p(Q_N)} + \enVert{u_N}_{W^{1,2}_p(Q_N)}},
	\end{equation*}
	where the value of $C_{k,\tau}$ may have increased.
	Iterate this argument finitely many times to get a sequence $u_1,\ldots,u_N$, through which we can define $u$ in a piecewise way, i.e.~$u(t,x) = u_i(t,x)$ for $t \in (t_i,t_{i+1}]$.
	Then $u$ satisfies the desired properties.
	Its uniqueness follows from the uniqueness of the fixed point on any small interval.
\end{proof}

\begin{lemma} \label{lem:existence Holder}
	Let $\alpha \in (0,1)$ and let $f \in C^{\alpha/2,\alpha}([0,T] \times [0,\infty))$, $u_T \in C^{2+\alpha}([0,\infty))$.
	Assume $a,b \in C^{\alpha/2,\alpha}([0,T] \times [0,\infty))$.
	Moreover, assume that the compatibility condition of order 1 is satisfied, namely
	\begin{equation}
	\label{eq:compatibility general}
	\frac{1}{2}\sigma^2 u_T''(0) + \s{I}[u_T](0) + a(T,0)u_T'(0) + b(T,0)u_T(0) = f(T,0).
	\end{equation}
	Then there exists a unique classical solution $u \in C^{1+\alpha/2,2+\alpha}([0,T] \times [0,\infty))$ to the boundary value problem \eqref{eq:parabolic}.
	Moreover, the following estimate holds:
	\begin{equation}
	\label{eq:parabolic estimate Holder}
	\enVert{u}_{C^{1+\alpha/2,2+\alpha}} \leq C\del{\enVert{f}_{C^{\alpha/2,\alpha}} + \enVert{u_T}_{C^{2+\alpha}}}
	\end{equation}
	where $C = C\del{\enVert{a}_{C^{\alpha/2,\alpha}},\enVert{b}_{C^{\alpha/2,\alpha}}}$.
\end{lemma}

\begin{proof}
	The proof is analogous to that of Lemma \ref{lem:existence Lp}, this time appealing to \cite[Theorem IV.5.2]{ladyzhenskaia1968linear} for the classical result and then applying Lemma \ref{lem:compactness estimate Holder} to get a fixed point.
	We omit the details.
\end{proof}

\section{A priori estimates}\label{sec:apriori}

The goal of this section is to prove a priori estimates on solutions to System \eqref{eq:HJBs2}.
In fact, we will consider a slight modification:
\begin{equation}\label{eq:HJBs2 lambda}
\begin{array}{lc}
u_t+\frac{1}{2}\sigma^2 u_{xx}-ru+\lambda H(t,u_x,{\color{red}\pi}[u_x,m],\eta)+{\s{I}}[u]=0,&  0<t<T,\ 0<x<\infty\\
m_t-\frac{1}{2}\sigma^2 m_{xx}-\lambda(D_\xi H(t,u_x,{\color{red}\pi}[u_x,m],\eta)m)_x-\s{I}^*[m]=0,& 0<t<T,\ 0<x<\infty\\
m(0,x)=\lambda m_0(x),\ u(T,x)=\lambda u_T(x), & 0\leq x< \infty\\
u(t,x)=m(t,x)=0, & 0\leq t\leq T,  -\infty< x\leq 0
\end{array}
\end{equation}
where $\lambda \in [0,1]$ is fixed.

We assume throughout that $(u,m)$ is a {\color{magenta}(}classical{\color{magenta})} solution to the system in the sense of Definition \ref{def:solution} and that the assumptions on the data hold true.

Let us give some motivation for this definition and, simultaneously, describe how to obtain a priori estimates in the specified function spaces.
If we can prove that the term $H(t,\xi,{\color{red}\pi}[u_x,m],\eta)$ is a priori bounded in the space $C^{\alpha/2,\alpha}$, then by Lemma \ref{lem:existence Holder} we obtain an a priori estimate on $u$ in the classical space $C^{1+\alpha/2,2+\alpha}([0,T] \times [0,\infty))$; the corresponding estimate on $m$ follows.
To do this, it is convenient to obtain bounds on $u_x$, which can be viewed as a solution to another parabolic equation.
Indeed, differentiating (formally) the Hamilton-Jacobi equation with respect to $x$, we get
\begin{equation} \label{eq:grad HJ}
u_{xt} + \frac{\sigma^2}{2}u_{xxx} - ru_x - \lambda D_\xi H(t,\xi,{\color{red}\pi}[u_x,m],\eta)u_{xx} + \s{I}[u_x] = 0
\end{equation}
in the sense of distributions.
If we assume $D_\xi H(t,\xi,{\color{red}\pi}[u_x,m],\eta)$ is bounded, then by classical theory, i.e.~Lemma \ref{lem:existence Lp}, we will be able to show that $u_x$ is a strong solution of \eqref{eq:grad HJ} with an estimate in $W_q^{1,2}$ for some $q > 1$, which immediately implies a H\"older estimate on $u_x$.
We can then also regard the Fokker-Planck equation as a parabolic equation with coefficients bounded in $L^q$, which implies $m$ has an estimate in $W^{1,2}_q$ as well.
All of this regularity comes in very useful when justifying the crucial estimate of the nonlocal term ${\color{red}\pi}[u_x, m] = \int u_x(t,x)m(t,x)\dif x$.
For this we want to use the fact that Equation \eqref{eq:grad HJ} is in \emph{duality} with the Fokker-Planck equation, so that a natural test function in \eqref{eq:grad HJ} would be $m$.
Once this estimate is obtained, the remaining estimates follow by bootstrapping.

\subsection{Standard estimates}
To begin with, we will use the maximum principle to get our first a priori bounds.

\begin{lemma}[Bounds on $u$] \label{lem:u bounded}
	For all $(t,x) \in [0,T] \times [0,\infty)$, we have 
	\begin{equation} \label{eq:u bounded}
	0 \leq u(t,x) \leq Te^{rT}\enVert{H(t,u_x,{\color{red}\pi}[u_x,m],\eta)}_\infty + \enVert{u_T}_\infty
	\end{equation}
\end{lemma}
\begin{proof}
	First set $v(t,x) = -e^{r(T-t)}u(t,x)$. Then as by Assumption~\ref{as:uT}  $u_T\leq 0$ and as by assumption $H(t,u_x,{\color{red}\pi}[u_x,m],\eta)\geq0$, we see that $v$ satisfies the assumption of  Lemma \ref{lem:max principle} and hence $v\leq 0$. This, in turn implies $u\geq 0$. Second, define $w(t,x) = e^{r(T-t)}u(t,x) - e^{rT}(T-t)\enVert{H(t,u_x,{\color{red}\pi}[u_x,m],\eta)}_\infty - \enVert{u_T}_\infty$.
	Then  $w$ also satisfies the hypotheses of Lemma \ref{lem:max principle} and hence $w \leq 0$, which ensures the second inequality.
\end{proof}

\begin{lemma} \label{lem:ux lower bound}
	[Lower bound on $u_x$]
	For all $(t,x) \in [0,T] \times [0,\infty)$ we have $u_x(t,x) \geq -c_0$ where $c_0 = -\min\{\min u_T',0\}$.
\end{lemma}

\begin{proof}
	Since $u(t,0) = 0 = \min u$, it follows that $u_x(t,0+) \geq 0$ for all $t \in [0,T]$.
	Define $w = -e^{r(T-t)}u_x + \min\{\min u_T',0\}$ and observe that
	\begin{equation}
	-w_t - \frac{\sigma^2}{2}w_{xx}- \s{I}[w] \leq - \lambda D_\xi H(t,u_x,{\color{red}\pi}[u_x,m],\eta)w_{x},
	\ w(t,x) \leq 0 \ \forall x \leq 0, \ w(T,x) \leq 0 \ \forall x \in \bb{R}.
	\end{equation}
	From Lemma \ref{lem:max principle weak} we deduce $w \leq 0$.
	The claim follows.
\end{proof}

\begin{lemma}[Bounds on $m$]
	\label{lem:FP basics}
	For all $(t,x) \in [0,T] \times [0,\infty)$ we have $m(t,x) \geq 0, m_x(t,0) \geq 0$, and $[0,T] \ni t \mapsto \eta(t):=\int_0^\infty m(t,x)\,dx \in [0,1]$ is decreasing with  $\eta(0)=1$, 
\end{lemma}
\begin{proof}
	We use a multiplier argument similar to Lemma~\ref{lem:max principle weak} to deduce $m \geq 0$, see also \cite{porretta2015weak}.
	Since $m(t,0) = 0 = \min m$, it follows that $m_x(t,0) \geq 0$.
	Now, 
	since $m \in C([0,T];L^1([0,\infty))) \cap W^{1,2}_p$ for $p > 3$ we deduce from a similar calculation as in Lemma \ref{lem:max principle weak} 
	that $\int_0^\infty \s{I}^*[m] \dif x\leq 0$. Hence  we
	integrate the Fokker-Planck equation  to get that
	\begin{equation}
	\label{eq:int m}
	\od{}{t} \eta (t)= \od{}{t} \int_0^\infty m(t,x)\dif x \leq  -\frac{1}{2}\sigma^2 m_x(t,0).
	\end{equation}
	We conclude that 
	$t \mapsto \eta(t):=\int_0^\infty m(t,x)\,dx \in [0,1]$ is decreasing with  $\eta(0)=\int m_0(x)\dif x=1$.
\end{proof}

\subsection{Energy estimates} \label{sec:energy estimates}

A crucial point to be observed is that, while we have a lower bound on $u$, the norm $\enVert{u}_\infty$ has not yet been estimated, because we do not have a bound on $\enVert{H(t,u_x,{\color{red}\pi}[u_x,m],\eta)}_\infty$.
We can rectify this by finding an estimate on $\int_0^\infty u_x(t,x)m(t,x)\dif x$.
We want to use the fact that Equation \eqref{eq:grad HJ} is in duality with the Fokker-Planck equation.
Indeed, if $(u_x,m)$ were smooth and the boundary $u_x(t,0) = 0$ were satisfied, then formally (leaving aside whether integration by parts is valid here) we get
\begin{equation}
\od{}{t}\left(e^{-rt}\int_0^\infty u_x(t,x)m(t,x)\dif x\right) = 0
\ \Rightarrow \
e^{-rt}\int_0^\infty u_x(t,x)m(t,x)\dif x = e^{-rT}\int_0^\infty u_T'(x)m(T,x)\dif x.
\end{equation}
Using the fact that $u_x \geq -c_0$, this would actually imply $\int_0^\infty  \abs{u_x(t,x)}m(t,x)\dif x$ is bounded.

The boundary condition $u_x(t,0) = 0$ is not satisfied, but we do have an estimate $u_x(t,0) \geq 0$ as well as $m_x(t,0) \geq 0$.
We exploit this in the following
\begin{lemma} \label{lem:[u_xm]}
	For all $t \in [0,T]$,
	\begin{equation} \label{eq:[u_xm]}
	e^{-rt}\int_0^\infty u_x(t,x)m(t,x)\dif x = e^{-rT}\int_0^\infty u_T'(x)m(T,x)\dif x
	- \frac{1}{2}\sigma^2\int_t^T e^{-rs}u_x(\tau,0)m_x(\tau,0)\dif \tau.
	\end{equation}
	In particular,  we obtain that
	\begin{equation}
	\label{eq:[u_xm]bound}
	\abs{\int_0^\infty u_x(t,x)m(t,x)\dif x} 
	\leq \enVert{u_T'}_\infty + 2c_0.
	\end{equation}		
\end{lemma}

\begin{proof}
	Since $u_x,u_{xx} \in L^p \cap C^{\alpha/2,\alpha}$ we can 
	use 
	$m$ as a test function of the PDE  \eqref{eq:grad HJ} satisfied by  $u_x$
	to get		
	\begin{equation} \label{eq:[u_xm]2}
	\begin{split}
	&e^{-rt}\int_0^\infty u_x(t,x)m(t,x)\dif x \\
	&= 
	e^{-rT}\int_0^\infty u_T'(x)m(T,x)\dif x
	- \frac{\sigma^2}{2}\int_t^T \int_0^\infty e^{-rs}\left(u_xm_{xx} + u_{xx}m_x\right)\dif x \dif \tau
	\\
	&= e^{-rT}\int_0^\infty u_T'(x)m(T,x)\dif x
	- \frac{\sigma^2}{2}\int_t^T \int_0^\infty e^{-rs}\left(u_xm_{x}\right)_x\dif x \dif \tau
	\end{split}
	\end{equation}
	where we have used Equation \eqref{eq:grad HJ} satisfied by $u_x$.
	Since $u_xm_x$ is both H\"older continuous and in $L^q$, it follows from Lemma \ref{lem:goes to zero} that
	$$
	\int_t^T e^{-r\tau}\abs{u_xm_x}(\tau,R) \dif \tau \to 0
	$$
	as $R \to \infty$, and so
	$$
	\int_t^T \int_0^\infty e^{-rs}\left(u_xm_{x}\right)_x\dif x \dif \tau = \lim_{R \to \infty} \int_t^T \int_0^{R} e^{-rs}\left(u_xm_{x}\right)_x\dif x \dif \tau = -  \int_t^T e^{-rs}u_x(\tau,0)m_x(\tau,0)\dif \tau.
	$$
	Equation \eqref{eq:[u_xm]} hence follows from \eqref{lem:max principle weak}.
	Now we use the fact that $u_x(\tau,0)m_x(\tau,0) \geq 0$ to deduce
	\begin{equation} \label{eq:[u_xm]3}
	\int_0^\infty u_x(t,x)m(t,x)\dif x \leq e^{r(t-T)}\int_0^\infty u_T'(x)m(T,x)\dif x \leq e^{r(t-T)}\enVert{u_T'}_\infty.
	\end{equation}
	Moreover, using the fact that $u_x + c_0 \geq 0$, we deduce that
	\begin{equation} \label{eq:[u_xm]4}
	\int_0^\infty \abs{u_x(t,x)}m(t,x)\dif x \leq \int_0^\infty (u_x(t,x)+c_0)m(t,x)\dif x + c_0
	\leq e^{r(t-T)}\enVert{u_T'}_\infty + 2c_0,
	\end{equation}
	and thus we arrive at \eqref{eq:[u_xm]bound}.
\end{proof}

\begin{remark}
	The idea used to prove Lemma \ref{lem:[u_xm]} was also used in \cite{graber2018existence}, where cut-off functions were used to deal with two incompatible boundary conditions on a bounded interval.
	The idea of proving an $L^\infty$ bound on $\int u_x m$ is also found in \cite{bertucci2019some}; however in that case the boundary conditions are not an obstacle.
\end{remark}

\subsection{Further a priori bounds}

Combining Lemma \ref{lem:ux lower bound} and Lemma \ref{lem:[u_xm]} with Assumptions \ref{as:coupling} and \ref{as:H}, we can find a constant $C\big(\enVert{u_T'}_\infty,c_0\big)>0$ only depending on the given boundary condition (see Assumption~\ref{as:uT}) such that
\begin{equation} \label{eq:bound on H}
\max\big\{\enVert{H(t,u_x,{\color{red}\pi}[u_x,m],\eta)}_\infty,\, \enVert{D_\xi H(t,u_x,{\color{red}\pi}[u_x,m],\eta)}_\infty\big\} \leq C\big(\enVert{u_T'}_\infty,c_0\big).
\end{equation}
In other words, $u$ satisfies a parabolic equation with an a priori bound on the right-hand side.
Now we are in a position to use a ``bootstrapping" argument to prove full regularity of the solution $(u,m)$.

We start by pointing out that $u$ is \emph{locally} in $W^{1,2}_p$ for arbitrary $p > 1$.
We then show that, away from the boundary $x = 0$, $u_x$ has an estimate in $W^{1,2}_q$.
We combine these estimates to prove a \emph{global} H\"older estimate on $u_x$.
\begin{lemma} \label{lem:uW21p-local}
	Let $p \in (1,\infty)$ be arbitrary.
	For every $M > 0$ there exists a constant $ C_{p,M}>0$ such that
	\begin{equation} \label{eq:uW21p-local}
	\enVert{u}_{W^{1,2}_p([0,T] \times [0,M])} \leq C_{p,M}\del{\enVert{H(t,u_x,{\color{red}\pi}[u_x,m],\eta)}_{\infty} + \enVert{u_T}_{\s{C}^2}}.
	\end{equation}
	As a corollary, if $p > 3$ then  (by \cite[Lemma II.3.3]{ladyzhenskaia1968linear}) we have for $\alpha := 1 - \tfrac{3}{p}$  that
	\begin{equation} \label{eq:uxHolder-local}
	\enVert{u_x}_{C^{\alpha/2,\alpha}([0,T] \times [0,M])} \leq C_{p,M}\del{\enVert{H(t,u_x,{\color{red}\pi}[u_x,m],\eta)}_{\infty} + \enVert{u_T}_{\s{C}^2}}.
	\end{equation}
\end{lemma}
\begin{proof}
	Let $\phi(x) = (1+x)^{-2}$ be defined on $(0,\infty)$.
	Note that $\phi \in W^2_p(0,\infty)$ for all $p \geq 1$ with 
	\begin{equation}
	\label{eq:phi estimates}
	\enVert{\phi}_p = (2p-1)^{-1/p},
	\qquad
	\enVert{\phi'}_p = 2(3p-1)^{-1/p},
	\qquad
	\enVert{\phi''}_p = 6(4p-1)^{-1/p}.
	\end{equation}
	Moreover, observe that 
	\begin{equation} \label{eq:phi relations}
	\tfrac{\phi'(x)}{\phi(x)} = -2(1+x)^{-1},  
	\quad \quad 
	\tfrac{\phi''(x)}{\phi(x)} = 6(1+x)^{-2},
	\end{equation}
	which are both bounded functions (by $2$ and $6$, respectively) that are also in $L^p$ for any $p > 1$.
	
	Now, define $w := \phi u$.
	Then we see that $w$ satisfies
	\begin{equation}
	w_t +  \tfrac{1}{2}\sigma^2w_{xx} - \sigma^2\tfrac{\phi'}{\phi}w_x + \Big(\tfrac{1}{2}\sigma^2\tfrac{\phi''}{\phi} - \sigma^2\del{\tfrac{\phi'}{\phi}}^2 -r\Big)w  
	+ \s{I}[w] = -\phi \lambda H(t,u_x,{\color{red}\pi}[u_x,m],\eta) + \s{I}[\phi u] - \phi\s{I}[u].
	\end{equation}
	%
	Moreover, we apply Lemma \ref{lem:existence Lp} and the estimates \eqref{eq:phi estimates} on $\phi$ to get the estimate
	\begin{equation} \label{eq:phi u estimate}
	\begin{split}
	\enVert{\phi u}_{W^{1,2}_p} 
	&\leq 
	\widetilde C_p\del{\enVert{\phi H(t,u_x,{\color{red}\pi}[u_x,m],\eta)}_{p} + \enVert{\phi u_T}_{W^{2-2/p}_p} + \enVert{\s{I}[\phi u] - \phi\s{I}[u]}_{p}}\\
	&\leq 
	\widetilde C_p\del{\enVert{H(t,u_x,{\color{red}\pi}[u_x,m],\eta)}_{\infty}\enVert{\phi}_{p} + \enVert{u_T}_{C^2}\enVert{\phi}_{W^{2}_p} + \enVert{\s{I}[\phi u] - \phi\s{I}[u]}_{p}}\\
	&\leq 
	\widetilde C_p\del{\enVert{H(t,u_x,{\color{red}\pi}[u_x,m],\eta)}_{\infty} + \enVert{u_T}_{C^2} + \enVert{\s{I}[\phi u] - \phi\s{I}[u]}_{p}},
	\end{split}
	\end{equation}
	where the value of $\widetilde C_p$ has possibly increased in each line.
	
	Next, we need to estimate $\enVert{\s{I}[\phi u] - \phi\s{I}[u]}_{L^p}$. 
	To that hand, we start by decomposing the big and small jumps in the integro-part, namely for  any $0<r<1$, we write for any $t,x$
	\begin{equation} \label{eq:Irho decomposition}
	\s{I}[\phi u](t,x) - \phi(x)\s{I}[u](t,x) = I_0(t,x) + I_1(t,x),
	\end{equation}
	where
	\begin{equation}
	\label{eq:I0}
	\begin{split}
	I_0(t,x)
	&=
	\int_{\bb{R} \setminus B_r(0)}\big[(\phi u)(t,x+z)- (\phi u)(t,x)-(\phi u)_x(t,x)\,z\mathbf{1}_{\{|z|\leq 1\}}\big]\,F(dz)\\
	& \quad 
	-
	\phi(x)\int_{\bb{R} \setminus B_r(0)}\big[u(t,x+z)- u(t,x)-u_x(t,x)\,z\mathbf{1}_{\{|z|\leq 1\}}\big]\,F(dz)
	\end{split}
	\end{equation}
	and
	\begin{equation}
	\label{eq:I1}
	\begin{split}
	I_1(t,x)
	&=
	\int_{B_r(0)}\big[(\phi u)(t,x+z)- (\phi u)(t,x)-(\phi u)_x(t,x)\,z\big]\,F(dz)\\
	& \quad
	-
	\phi(x)\int_{B_r(0)}\big[u(t,x+z)- u(t,x)-u_x(t,x)\,z\big]\,F(dz).
	\end{split}
	\end{equation}
	Let us first estimate 	\eqref{eq:I0}.  To that hand, note that by definition
	\begin{equation*}
	\begin{split}
	I_0(t,x)
	&= \int_{\bb{R} \setminus B_r(0)}
	\big[
	(\phi u)(t,x+z) -\phi(x)u(t,x+z)
	-\phi_x(x)u(t,x)z\mathbf{1}_{\{|z|\leq 1\}}
	\big]\,F(dz).
	\end{split}
	\end{equation*}
	Therefore, using H\"older's inequality and as $\widetilde{C}_r:=\int_{\bb{R} \setminus B_r(0)} 1 F(dz)<\infty$
	we see that
	\begin{equation*}
	\begin{split}
	\big|I_0(t,x)|
	&\leq
	\int_{\bb{R} \setminus B_r(0)} 	\big|(\phi u)(t,x+z)\big|\, F(dz)
	+
	\int_{\bb{R} \setminus B_r(0)} \big[|\phi(x)|\,|u(t,x+z)|
	+
	|\phi_x(x)|\,|u(t,x)|\big]\, F(dz)\\
	&\leq
	\int_{\bb{R} \setminus B_r(0)} 	\big|(\phi u)(t,x+z)\big|\, F(dz)
	+
	\Vert u\Vert_\infty \int_{\bb{R} \setminus B_r(0)} 
	\big[|\phi(x)| +|\phi_x(x)|\big]
	\, F(dz)\\
	&\leq \widetilde{C_r}^{\nicefrac{1}{q}} \bigg(\int_{\bb{R} \setminus B_r(0)} 	\big|(\phi u)(t,x+z)\big|^p\, F(dz)\bigg)^{\nicefrac{1}{p}}
	+  	
	\Vert u\Vert_\infty\widetilde{C}_r \big[|\phi(x)| +|\phi_x(x)|\big],
	\end{split}
	\end{equation*}
	where here $q:=\nicefrac{p}{p-1}$ denotes the H\"older conjugate of $p$.
	As a consequence, we obtain by Fubini's theorem that
	\begin{equation*}
	\begin{split}
	\Vert I_0 \Vert_p^p
	& \leq 
	2^{2(p-1)}\int_0^T \int_{0}^{\infty}
	\bigg[  
	\widetilde{C_r}^{\nicefrac{p}{q}} \int_{\bb{R} \setminus B_r(0)} 	\big|(\phi u)(t,x+z)\big|^p\, F(dz) + \Vert u\Vert_\infty^p\widetilde{C}_r^p \big[|\phi(x)|^p +|\phi_x(x)|^p\big] \bigg]\, dx \,dt\\
	&=2^{2(p-1)}  \widetilde{C_r}^p \Vert\phi u \Vert_p^p
	+
	2^{2(p-1)} \widetilde{C}_r^p T \Vert u\Vert_\infty^p \big[\Vert \phi \Vert_p^p + \Vert \phi_x \Vert_p^p\big]\\
	&\leq
	2^{2(p-1)}  \widetilde{C_r}^p \Vert u \Vert_\infty^p \Vert\phi  \Vert_p^p
	+
	2^{2(p-1)} \widetilde{C}_r^p T \Vert u\Vert_\infty^p \big[\Vert \phi \Vert_p^p + \Vert \phi_x \Vert_p^p\big]\\
	\end{split}
	\end{equation*}
	This, together with the estimates on $\phi$ provided in \eqref{eq:phi estimates} stating that both $\Vert \phi \Vert_p$ and $\Vert \phi_x \Vert_p$ are finite,  as well as Lemma~\ref{lem:u bounded} ensures that
	\begin{equation}
	\label{eq:I0 estimate}
	\begin{split}
	\Vert I_0 \Vert_p
	\leq  C_r\Vert u \Vert_\infty
	&\leq C_r \Big(Te^{rT}\enVert{H(t,u_x,{\color{red}\pi}[u_x,m],\eta)}_\infty + \enVert{u_T}_\infty\Big)\\
	&\leq C_r \Big(\enVert{H(t,u_x,{\color{red}\pi}[u_x,m],\eta)}_\infty + \enVert{u_T}_\infty\Big)
	\end{split}
	\end{equation}
	for some constant $C_r>0$ whose value has possibly increased in the second line.

	Next, we estimate \eqref{eq:I1}. To that end, note that by Taylor expansion we have that
	\begin{equation*}
	\begin{split}
	&I_1(t,x) \\
	&= 
	\int_{-r}^r\int_0^1 \int_{0}^1 (\phi u)_{xx}(t,x+\tau \sigma z)\sigma z^2\dif \tau \dif \sigma\, F(dz)
	-
	\phi(x)\int_{-r}^r\int_0^1 \int_{0}^1 u_{xx}(t,x+\tau \sigma z)\sigma z^2 \dif \tau \dif \sigma \,F(dz).
	\end{split}
	\end{equation*}
	This together with the formula
	\begin{equation*}
	(\phi u)_{xx} = \phi u_{xx} + \tfrac{2\phi'}{\phi}(\phi u)_x + \del{\tfrac{\phi''}{\phi} - \tfrac{2(\phi')^2}{\phi^2}}\phi u
	\end{equation*}
	implies that
	\begin{equation*} 
	I_1(t,x)
	= \int_{-r}^r\int_0^1 \int_{0}^1 \del{\tfrac{2\phi'}{\phi}(\phi u)_x + \del{\tfrac{\phi''}{\phi} - \tfrac{2(\phi')^2}{\phi^2}}\phi u}\big(t,x+\tau \sigma z\big)\,\sigma z^2 \dif \tau \dif \sigma \,F(dz).
	\end{equation*}
	Therefore, using \eqref{eq:phi relations} and H\"older's inequality applied to the measure $\nu(dz):= z^2\dif \tau \dif \sigma\, F(dz)$, which is finite
	on $[-1,1]\times[0,1]\times[0,1]$, we obtain the estimates
	\begin{equation*}
	\begin{split}
	\abs{I_1(t,x)}
	&\leq 
	4\int_{-r}^r\int_0^1 \int_{0}^1 \abs{(\phi u)_x(t,x+\tau \sigma z)} z^2 \dif \tau \dif \sigma \,F(dz)\\
	& + \quad
	14\int_{-r}^r\int_0^1 \int_{0}^1 \abs{(\phi u)(t,x+\tau \sigma z)} z^2 \dif \tau \dif \sigma \,F(dz)\\
	&\leq 
	4\widetilde c_r^{\nicefrac{1}{q}}\bigg(\int_{-r}^r\int_0^1 \int_{0}^1 \abs{(\phi u)_x(t,x+\tau \sigma z)}^p z^2 \dif \tau \dif \sigma \,F(dz) \bigg)^{\nicefrac{1}{p}}
	\\
	& \quad 
	+ 14 \widetilde c_r^{\nicefrac{1}{q}}\bigg(\int_{-r}^r\int_0^1 \int_{0}^1 \abs{(\phi u)(t,x+\tau \sigma z)}^p z^2 \dif \tau \dif \sigma \,F(dz) \bigg)^{\nicefrac{1}{p}}
	\end{split}
	\end{equation*}
	where here $q:=\nicefrac{p}{p-1}$ denotes the H\"older conjugate of $p$ and $\widetilde c_r := \int_{-r}^r\int_0^1 \int_{0}^1 1\, z^2\dif \tau \dif \sigma\, F(dz)<\infty$, which satisfies that $\widetilde c_r \to 0$ as $r \to 0$.
	Therefore, we obtain by Fubini's theorem that
	\begin{equation}
	\begin{split}
	\Vert I_1 \Vert_p^p
	& \leq
	2^p 4^p \widetilde c_r^{\nicefrac{p}{q}}\int_0^T \int_0^\infty\int_{-r}^r\int_0^1 \int_{0}^1 \abs{(\phi u)_x(t,x+\tau \sigma z)}^p z^2 \dif \tau \dif \sigma \,F(dz)\, dx \, dt\\
	& \quad 
	+ 2^p 14^p \widetilde c_r^{\nicefrac{p}{q}}\int_0^T \int_0^\infty\int_{-r}^r\int_0^1 \int_{0}^1 \abs{(\phi u)(t,x+\tau \sigma z)}^p z^2 \dif \tau \dif \sigma \,F(dz)\, dx \, dt\\
	& \leq
	2^p 4^p \widetilde c_r^{{p}} \Vert (\phi u)_x \Vert_p^p
	+ 2^p 14^p \widetilde c_r^{{p}} \Vert \phi u \Vert_p^p,
	\end{split}
	\end{equation}
	from which we deduce
	\begin{equation} \label{eq:I1 estimate}
	\enVert{I_1}_p 
	\leq  
	c_r \big[ \Vert (\phi u)_x \Vert_p + \Vert \phi u \Vert_p\big]
	\leq 
	c_r \Vert \phi u \Vert_{W^{1,2}_p},
	\end{equation}
	where $c_r>0$ satisfies that $ c_r \to 0$ as $r \to 0$.
	
	We now plug \eqref{eq:I0 estimate} and \eqref{eq:I1 estimate} into \eqref{eq:Irho decomposition} to deduce
	\begin{equation} \label{eq:Irho difference stimate}
	\enVert{\s{I}[\phi u] - \phi\s{I}[u]}_{p}
	\leq 
	c_r \Vert \phi u \Vert_{W^{1,2}_p}
	+
	C_r\del{\enVert{H(t,u_x,{\color{red}\pi}[u_x,m],\eta)}_\infty + \enVert{u_T}_\infty}.
	\end{equation}
	Therefore, using \eqref{eq:Irho difference stimate} in \eqref{eq:phi u estimate}  shows that
	\begin{equation*}
	\begin{split}
	&\Vert \phi u \Vert_{W^{1,2}_p}\\
	&\leq 
	\widetilde C_p\del{\enVert{H(t,u_x,{\color{red}\pi}[u_x,m],\eta)}_{\infty} + \enVert{u_T}_{C^2} 
		+
		c_r \Vert \phi u \Vert_{W^{1,2}_p}
		+
		C_r\del{\enVert{H(t,u_x,{\color{red}\pi}[u_x,m],\eta)}_\infty + \enVert{u_T}_\infty}
	},
	\end{split}
	\end{equation*}
	which in turn ensures for small enough $r$, as $ c_r \to 0$ when $r \to 0$, that
	\begin{equation*}
	\begin{split}
	\Vert \phi u \Vert_{W^{1,2}_p}
	&\leq \tfrac{\widetilde C_p}{1-\widetilde C_pc_r}
	\del{(1+C_r)\big[\enVert{H(t,u_x,{\color{red}\pi}[u_x,m],\eta)}_{\infty} + \enVert{u_T}_{C^2} \big]
	}.
	\end{split}
	\end{equation*}
	As a consequence,	as $ c_r \to 0$ when $r \to 0$,	we can take $r$ small and $C_{p,T}$ large enough so that
	\begin{equation} 
	\enVert{\phi u}_{W^{1,2}_p} 
	\leq C_{p,T}\del{\enVert{H(t,u_x,{\color{red}\pi}[u_x,m],\eta)}_{L^\infty} + \enVert{u_T}_{C^2}}.
	\end{equation}
	
	Equation \eqref{eq:uW21p-local} now follows using the fact that $\phi$ is a smooth function.
\end{proof}
As another intermediary step, we derive some estimates in $W^{1,2}_p$ on $u_x$ and $m$ (rather than on $u$ itself).
\begin{lemma} \label{lem:uxW21p}
	Let $p > 3, q > 1$.
	Fix $\delta > 0$.
	Then we have estimates
	\begin{equation} \label{eq:uxW21p}
	\enVert{u_x}_{W^{1,2}_p\del{[0,T] \times (\delta,\infty)}} \leq C\del{\enVert{u_T}_{\s{C}^2},c_0,p,T,\delta}
	\end{equation}
	and 
	\begin{equation} \label{eq:mW21p}
	\enVert{m}_{W_q^{1,2}} \leq C\del{\enVert{m_0}_{W_q^{2-2/q}},\enVert{u_T}_{\s{C}^2},\enVert{u_T'}_{W^{2}_p},c_0,p,q,T}.
	\end{equation}
	
	As a consequence, $u_x, m$, and $m_x$ are H\"older continuous with
	\begin{equation}\label{eq:uxHolder}
	\enVert{u_x}_{C^{\alpha/2,\alpha}([0,T] \times [0,\infty))} \leq C\del{\enVert{u_T'}_{W^{2-2/p}_p}, \enVert{u_T}_{C^2},c_0,p,T,\delta}
	\end{equation}
	and
	\begin{equation} \label{eq:mHolder}
	\max\Big\{\enVert{m}_{C^{\beta/2,\beta}([0,T] \times [0,\infty))}, \enVert{m_x}_{C^{\beta/2,\beta}([0,T] \times [0,\infty))}
	\Big\}
	\leq C\del{\enVert{m_0}_{W_p^{2-2/q}},\enVert{u_T}_{\s{C}^2},\enVert{u_T'}_{W^{2}_p},c_0,p,T}
	\end{equation}
	for $\alpha = 1 - \frac{3}{p}, \beta = 1 - \frac{3}{q}$.
\end{lemma}
\begin{proof}
	Without loss of generality we will assume $\delta \leq 1$.
	Let $\psi(x)$ be a smooth function with support in $(0,\infty)$ such that $\abs{\psi'} \leq 2/\delta$, $\abs{\psi''} \leq 4/\delta^2$, and $\psi(x) = 1$ for all $x \geq \delta$.
	Set $w = \psi u_x$.
	Then, by using the notation ${\color{red}\{\psi,\s{I}\}}(u_x) := \psi \s{I}[u_x] - \s{I}[\psi u_x]$, we see that $w$ satisfies
	\begin{equation} \label{eq:HJgrad psi}
	\begin{split}
	&w_t -rw + \dfrac{\sigma^2}{2}w_{xx} + \s{I}[w]
	- {\color{red}\lambda}D_\xi H\del{t,u_x,{\color{red}\pi}[u_x,m],\eta}w_x
	\\
	&= \dfrac{\sigma^2}{2}\del{2\psi' u_{xx} + \psi'' u_x}
	+ \psi' {\color{red}\lambda}D_\xi H\del{t,u_x,{\color{red}\pi}[u_x,m],\eta}u_{x} + {\color{red}\{\psi,\s{I}\}}(u_x)
	\end{split}
	\end{equation}
	with $w(x,t) = 0$ for $x \leq 0$.
	Therefore, by setting
	\begin{equation}
	g(x,t) := \dfrac{\sigma^2}{2}\del{2\psi' u_{xx} + \psi'' u_x}
	+ \psi' {\color{red}\lambda}D_\xi H\del{t,u_x,{\color{red}\pi}[u_x,m],\eta}u_{x} + {\color{red}\{\psi,\s{I}\}}(u_x)
	\end{equation}
	we get
	\begin{equation*}
	\begin{split}
	\enVert{g}_{p} &\leq \frac{2\sigma^2}{\delta^2}\del{\enVert{u_{xx}}_{L^p([0,1] \times [0,T])} + \enVert{u_x}_{L^p([0,1] \times [0,T])}}\\
	& \quad 
	+ \frac{2}{\delta}\enVert{D_\xi H\del{\cdot,u_x,{\color{red}\pi}[u_x,m],\eta}}_\infty \enVert{u_{xx}}_{L^p([0,1] \times [0,T])}
	+ \enVert{{\color{red}\{\psi,\s{I}\}}(u_x)}_p
	\end{split}
	\end{equation*}
	Using \eqref{eq:uW21p-local} and \eqref{eq:Irho difference stimate} from Lemma \ref{lem:uW21p-local}, we see that for any $\varepsilon > 0$ we have
	\begin{equation} \label{eq:gLp}
	\begin{split}
	\enVert{g}_{p} 
	&\leq 
	\del{\frac{2\sigma^2}{\delta^2}+C_{p}\frac{2}{\delta}\enVert{D_\xi H\del{\cdot,u_x,{\color{red}\pi}[u_x,m],\eta}}_\infty}\del{\enVert{H(t,u_x,{\color{red}\pi}[u_x,m],\eta)}_{L^\infty} + \enVert{u_T}_{\s{C}^2}}\\
	& \quad 
	+ \varepsilon\enVert{w}_{W^{1,2}_p} +C_\varepsilon\del{T\enVert{H(t,u_x,{\color{red}\pi}[u_x,m],\eta)}_\infty + \enVert{u_T}_\infty}.
	\end{split}
	\end{equation}
	
	Applying Lemma~\ref{lem:existence Lp} to \eqref{eq:HJgrad psi} 
	we obtain the estimate
	\begin{equation}
	\enVert{w}_{W^{1,2}_p} \leq C\del{\enVert{g}_p + \enVert{\phi u_T}_{W^{2-2/p}_p}},
	\end{equation}
	where $C$ depends on $\enVert{D_\xi H\del{\cdot,u_x,{\color{red}\pi}[u_x,m],\eta}}_\infty$ and $r$.
	Combining \eqref{eq:gLp} with \eqref{eq:bound on H} and taking $\varepsilon$ small enough, we obtain
	\begin{equation}
	\enVert{w}_{W^{1,2}_p} \leq C\del{\enVert{u_T}_{\s{C}^2},c_0,p,T,\delta},
	\end{equation}
	from which we deduce \eqref{eq:uxW21p}.
	
	Equation \eqref{eq:uxHolder} now follows from \cite[Lemma II.3.3]{ladyzhenskaia1968linear} combined with \eqref{eq:uxHolder-local}.
	In particular, we have an a priori bound on $\enVert{u_x}_\infty$.
	
	Similarly, write the Fokker-Planck equation in the form
	\begin{equation} \label{eq:FP expanded}
	m_t
	-\frac{1}{2}\sigma^2 m_{xx}
	-\lambda D_\xi H(t,u_x,{\color{red}\pi}[u_x,m],\eta)m_x 
	- \lambda D_{\xi\xi} H(t,u_x,{\color{red}\pi}[u_x,m],\eta)u_{xx} m 
	-\s{I}^*[m]
	=0.
	\end{equation}
	By Assumption~\ref{as:H} and 
	\eqref{eq:uxW21p} and \eqref{eq:bound on H}, we have that $D_\xi H(t,u_x,{\color{red}\pi}[u_x,m],\eta)$ and
	$D_{\xi\xi}H(t,u_x,{\color{red}\pi}[u_x,m],\eta)u_{xx}$ are in $L^q_{loc}$ with
	\begin{equation}\label{eq:l4.7-1}
	\enVert{D_\xi H(t,u_x,{\color{red}\pi}[u_x,m],\eta)}_{L^q_{loc}} \leq \enVert{D_\xi H(t,u_x,{\color{red}\pi}[u_x,m],\eta)}_{L^\infty} \leq C\del{\enVert{u_T'}_\infty, c_0}
	\end{equation} 
	and 
	\begin{equation}\label{eq:l4.7-2}
	\begin{split}
	\enVert{D_{\xi\xi}H(t,u_x,{\color{red}\pi}[u_x,m],\eta)u_{xx}}_{L^q_{loc}} 
	&\leq C\enVert{D_{\xi\xi}H(t,u_x,{\color{red}\pi}[u_x,m],\eta)}_{L^\infty}\enVert{u_{xx}}_{L^p} \\
	&\leq C\del{\enVert{u_x}_\infty,\enVert{u_x}_{W^{1,2}_p\del{[0,T] \times (1,\infty)}},\enVert{u}_{W^{1,2}_p\del{[0,T] \times [0,1]}}}.
	\end{split}
	\end{equation}
	Thus Lemma \ref{lem:existence Lp} implies
	\begin{equation}
	\enVert{m}_{W_q^{1,2}} \leq C\del{\enVert{m_0}_{W_q^{2-2/q}},\enVert{D_\xi H(t,u_x,{\color{red}\pi}[u_x,m],\eta)}_{L^q_{loc}},\enVert{D_{\xi\xi}H(t,u_x,{\color{red}\pi}[u_x,m],\eta)u_{xx}}_{L^q_{loc}},q,T},
	\end{equation}
	which by using \eqref{eq:l4.7-1} and \eqref{eq:l4.7-2} can be rewritten as \eqref{eq:mW21p}.		
	Equation \eqref{eq:mHolder} now follows from \cite[Lemma II.3.3]{ladyzhenskaia1968linear} as before.
\end{proof}

The last step is to prove that the nonlocal coupling term ${\color{red}\pi}[u_x, m]$ is H\"older continuous in time.
We will need a preliminary lemma.
\begin{lemma}
	\label{lem:eta continuity}
	The mass function $\eta(t) = \int_0^\infty m(t,x)\dif x$ is H\"older continuous with
	\begin{equation} \label{eq:eta holder}
	\abs{\eta(t_1)-\eta(t_2)} \leq C\abs{t_1-t_2}^{\alpha} \quad \forall t_1,t_2 \in [0,T],
	\end{equation}
	where $C,\alpha$ depend only on the data.
\end{lemma}

\begin{proof}
	Fix $t_1,t_2 \in [0,T]$ with $t_1 \neq t_2$. For any $\varepsilon > 0$ let $\zeta\equiv\zeta^{(\varepsilon)}
	\in C^2([0,\infty)$ be such that $\zeta \equiv 1$ on $[2\varepsilon,\infty)$, $\zeta \equiv 0$ on $[0,\varepsilon]$, $0 \leq \zeta \leq 1$,  $\abs{\zeta'} \leq 2/\varepsilon$, and $\abs{\zeta''} \leq 4/\varepsilon^2$.
	We use $\zeta$ as a test function in Equation \eqref{eq:HJBs2 lambda}(ii) to get
	\begin{equation} \label{eq:zeta dif}
	\begin{split}
	&\int_0^\infty \zeta(x)m(t_1,x)\dif x - \int_0^\infty \zeta(x)m(t_2,x)\dif x\\
	&= 
	\int_{t_2}^{t_1}\int_0^\infty \del{-\zeta' m_x + \s{I}[\zeta]m - \zeta'D_\xi H\del{u_x,{\color{red}\pi}[u_x,m]}m}\dif x \dif t.
	\end{split}
	\end{equation}
	Let $q > 1$ and let $q' := q/(q-1)$ be its H\"older conjugate.
	Applying H\"older's inequality to \eqref{eq:zeta dif} we have
	\begin{equation} \label{eq:zeta dif1}
	\begin{split}
	&\abs{\int_0^\infty \zeta(x)m(t_1,x)\dif x - \int_0^\infty \zeta(x)m(t_2,x)\dif x}\\
	&\leq 
	\del{\int_{t_2}^{t_1}\int_\varepsilon^{2\varepsilon} \abs{\zeta'}^{q'} \dif x \dif t}^{1/q'}\del{\enVert{m_x}_{L^q} + \enVert{D_\xi H\del{u_x,{\color{red}\pi}[u_x,m]}m}_{L^q}}
	+ \abs{\int_{t_2}^{t_1}\int_0^\infty \s{I}[\zeta]m\dif x \dif t}.
	\end{split}
	\end{equation}		
	To deal with the nonlocal term, we apply Fubini's theorem and integrate by  to obtain that
	\begin{equation*}
	\begin{split}
	\int_{t_2}^{t_1}\int_0^\infty \s{I}[\zeta]m\dif x \dif t
	&= \int_{t_2}^{t_1}\int_0^\infty \int_{-1}^1 \int_0^1 \int_0^1  \zeta''(x+\tau' \tau z) \tau z^2 m(t,x)\dif \tau' \dif \tau F(\dif z) \dif x \dif t\\
	& \quad 
	+ \int_{t_2}^{t_1}\int_0^\infty \int_{\bb{R}\setminus (-1,1)} \del{\zeta(x+z)-\zeta(x)}m(t,x)F(\dif z) \dif x \dif t\\
	&= -\int_{t_2}^{t_1}\int_0^\infty \int_{-1}^1 \int_0^1 \int_0^1 \zeta'(x+\tau' \tau z) \tau z^2 m_x(t,x)\dif \tau' \dif \tau F(\dif z) \dif x \dif t\\
	& \quad + \int_{t_2}^{t_1}\int_0^\infty \int_{\bb{R}\setminus (-1,1)} \del{\zeta(x+z)-\zeta(x)}m(t,x)F(\dif z) \dif x \dif t.
	\end{split}
	\end{equation*}
	Hence, using the estimate that $0\leq\int m(t,x)\dif x \leq 1$ for all $t$, we have for  $C_1 := \int_{-1}^1 z^2 F(\dif z)$ and $C_2 := 2\int_{\bb{R}\setminus (-1,1)} F(\dif z)$ that
	\begin{equation} \label{eq:zeta dif2}
	\abs{\int_{t_2}^{t_1}\int_0^\infty \s{I}[\zeta]m\dif x \dif t}
	\leq C_1 \del{\int_{t_2}^{t_1}\int_0^\infty \abs{\zeta'}^{q'} \dif x \dif t}^{1/q'}\enVert{m_x}_{L^q} + C_2\abs{t_1-t_2}.
	\end{equation}
	Plugging \eqref{eq:zeta dif2} into \eqref{eq:zeta dif1} we get
	\begin{equation} \label{eq:zeta dif3}
	\begin{split}
	&\abs{\int_0^\infty \zeta(x)m(t_1,x)\dif x - \int_0^\infty \zeta(x)m(t_2,x)\dif x}\\
	&\leq \del{\int_{t_2}^{t_1}\int_\varepsilon^{2\varepsilon} \abs{\zeta'}^{q'} \dif x \dif t}^{1/q'}\del{(C_1+1)\enVert{m_x}_{L^q} + \enVert{D_\xi H\del{u_x,{\color{red}\pi}[u_x,m]}m}_{L^q}}
	+ C_2\abs{t_1-t_2}.
	\end{split}
	\end{equation}	
	Moreover, recall that	$\enVert{m_x}_{L^p}$ and $\enVert{D_\xi H\del{u_x,{\color{red}\pi}[u_x,m]}m}_{L^q}$ are both estimated by Lemma~\ref{lem:uxW21p}.
	Now observe that by the choice of $\zeta$, we have that
	\begin{equation}
	\del{\int_{t_2}^{t_1}\int_\varepsilon^{2\varepsilon} \abs{\zeta'}^{q'} \dif x \dif t}^{1/q'}
	\leq 2\varepsilon^{-1/q}\abs{t_1-t_2}^{1/q'}.
	\end{equation}
	Therefore by \eqref{eq:zeta dif3} there exists a constant $C$ depending only on the data such that
	\begin{equation} \label{eq:zeta dif4}
	\abs{\int_0^\infty \zeta(x)m(t_1,x)\dif x - \int_0^\infty \zeta(x)m(t_2,x)\dif x}
	\leq C\del{\varepsilon^{-1/q}\abs{t_1-t_2}^{1/q'} + \abs{t_1-t_2}}.
	\end{equation}
	On the other hand, as $0\leq \zeta \leq 1$ with $\zeta=1$ on $[2\varepsilon,\infty)$  and using the $L^\infty$ bound on $m$ from \eqref{eq:mHolder}, we have
	\begin{equation} \label{eq:zeta dif5}
	\abs{\int_0^\infty \del{1-\zeta(x)}m(t_1,x)\dif x - \int_0^\infty \del{1-\zeta(x)}m(t_2,x)\dif x}
	\leq C\varepsilon.
	\end{equation}
	Now set $\varepsilon = \abs{t_1-t_2}^\alpha$ where $\alpha = \frac{q-1}{q+1}$.
	Adding together \eqref{eq:zeta dif4} and \eqref{eq:zeta dif5} we then deduce
	\begin{equation}
	\abs{\eta(t_1) - \eta(t_2)}
	\leq C\del{\abs{t_1-t_2}^{\alpha} + \abs{t_1-t_2}}.
	\end{equation}
	Estimate \eqref{eq:eta holder} follows.
\end{proof}
\begin{lemma}
	\label{lem:[uxm]Holder}
	Let $p > 3$ and define $q = p' = \frac{p}{p-1} \in (1,3/2)$ to be its H\"older conjugate.
	Then there exists 
	$\alpha \in (0,1)$ and a constant $C$ depending only on the data such that
	\begin{equation}\label{eq:[uxm]Holder}
	\enVert{{\color{red}\pi}[u_x,m]}_{C^\alpha([0,T])} \leq C.
	\end{equation}
\end{lemma}

\begin{proof}
	First we establish estimates on $u_x$ and $m$ in $C^\alpha([0,T];L^p(0,\infty))$ and $C^\alpha([0,T];L^q(0,\infty))$, respectively.
	Set $\gamma = 1 - \frac{3}{p}$.
	Given any $t_1,t_2 \in [0,T]$ and any $\phi \in L^{q'}(0,\infty)$, we use Lemma \ref{lem:uxW21p} and  H\"older's inequality to estimate
	\begin{equation*}
	\begin{split}
	\int_0^\infty &\del{u_x(t_1,x)-u_x(t_2,x)}\phi(x)\dif x\\
	&= \int_0^1 \del{u_x(t_1,x)-u_x(t_2,x)}\phi(x)\dif x + \int_1^\infty \int_{t_2}^{t_1}u_{xt}(t,x)\phi(x)\dif x\\
	&\leq \enVert{u_x}_{C^{\gamma/2,\gamma}}\enVert{\phi}_{L^{q'}(0,1)}\abs{t_1-t_2}^{\gamma/2} + \enVert{u_{xt}}_{L^p([0,T] \times [1,\infty))} \enVert{\phi}_{L^{q}(1,\infty)}\abs{t_1-t_2}^{1/q}\\
	&\leq C\enVert{\phi}_{L^{q}(0,\infty)}\del{\abs{t_1-t_2}^{1/q}+\abs{t_1-t_2}^{\gamma/2}},
	\end{split}
	\end{equation*}
	where $C = C\del{\enVert{u_T'}_{W^{2-2/p}_p}, \enVert{u_T}_{C^2},c_0,p,T}$.
	Therefore, since $\phi \in L^q(0,\infty)$ is arbitrary it follows that
	\begin{equation} \label{eq:ux time holder}
	\enVert{u_x(t_1,\cdot)-u_x(t_2,\cdot)}_{L^p(0,\infty)} \leq C \abs{t_1-t_2}^{\alpha}, \ \alpha \leq \min\cbr{\gamma/2,1/q}.
	\end{equation}
	By a similar argument, we directly get that
	\begin{equation} \label{eq:m time holder}
	\enVert{m(t_1,\cdot)-m(t_2,\cdot)}_{L^q(0,\infty)} \leq \enVert{m_t}_{L^q([0,T] \times [0,\infty)} \abs{t_1-t_2}^{1/p} \leq C\abs{t_1-t_2}^{1/p}
	\end{equation}
	again using Lemma \ref{lem:uxW21p}.
	As $u_T' \in L^p(0,\infty)$ it also follows that $t \mapsto u_x(t,\cdot)$ is bounded from $[0,T]$ to $L^p(0,\infty)$.
	Similarly, since $m_0 \in L^q(0,\infty)$, it follows that $m(t,\cdot)$ is bounded from $[0,T]$ to $L^q(0,\infty)$.
	
	Now by Assumption \ref{as:coupling continuity}, Lemma \ref{lem:eta continuity}, and Equations \eqref{eq:ux time holder} and \eqref{eq:m time holder}, we have
	\begin{equation*}
	\begin{split}
	\abs{{\color{red}\pi}[u_x,m](t_1)-{\color{red}\pi}[u_x,m](t_2)} 
	&\leq C\sup_{t \in [0,T]}\enVert{u_x(t,\cdot)}_{L^{p}(0,\infty)}\enVert{m(t_1,\cdot)-m(t_2,\cdot)}_{L^{q}(0,\infty)}\\
	&\quad + C\sup_{t \in [0,T]}\enVert{m(t,\cdot)}_{L^{q}(0,\infty)}\enVert{u_x(t_1,\cdot)-u_x(t_2,\cdot)}_{L^{p}(0,\infty)} + \abs{\eta(t_1)-\eta(t_2)}\\
	&\leq C\abs{t_1-t_2}^\alpha,
	\end{split}
	\end{equation*}
	where $\alpha$ is minimum of $\gamma/2,1/p,$ and the exponent appearing in \eqref{eq:eta holder}.	
\end{proof}

Finally, we deduce smoothness of both $u$ and $m$.
\begin{lemma}
	\label{lem:smoothness}
	Let $\alpha \in (0,1)$ be sufficiently small, and let $p > 1$ be arbitrarily large.
	Then any solution $(u,m)$ of \eqref{eq:HJBs2 lambda} satisfies
	\begin{equation} \label{eq:smoothness}
	\enVert{u}_{C^{1+\alpha/2,2+\alpha}} + \enVert{m}_{W^{1,2}_p} + \enVert{m}_{C^{\alpha}([0,T];L^1)} \leq C\del{m_0,u_T,c_0,\alpha,\varepsilon,\sigma,T}.
	\end{equation}
\end{lemma}

\begin{proof}
	By Lemmas \ref{lem:uxW21p}, \ref{lem:[uxm]Holder}, and \ref{lem:eta continuity} combined with Assumption \ref{as:H}, it follows that $H(t,u_x,{\color{red}\pi}[u_x,m],\eta)$ is H\"older continuous with norm estimated by a constant depending on the data.
	Now we apply Lemma \ref{lem:existence Holder} to Equation \eqref{eq:HJBs2 lambda}(i) to get
	\begin{equation} \label{eq:u regular}
	\enVert{u}_{C^{1+\alpha/2,2+\alpha}} \leq C\del{\enVert{u_T}_{C^{2+\alpha}} + \enVert{H(t,u_x,{\color{red}\pi}[u_x,m],\eta)}_{C^\alpha([0,T] \times [0,\infty))}}.
	\end{equation}
	Then it follows that all the coefficients in Equation \eqref{eq:HJBs2 lambda}(ii), equivalently \eqref{eq:FP expanded}, are bounded by a constant depending on the data.
	We can use Lemma \ref{lem:existence Lp} applied to the Fokker-Planck equation in the form of \eqref{eq:FP expanded} to get
	\begin{equation} \label{eq:m regular}
	\enVert{m}_{W^{1,2}_p} \leq C\del{\enVert{m_0}_{W^{2}_p},\enVert{D_\xi H(t,u_x,{\color{red}\pi}[u_x,m],\eta)}_\infty,
		\enVert{D_{\xi\xi} H(t,u_x,{\color{red}\pi}[u_x,m],\eta)u_{xx}}_\infty
	}.
	\end{equation}
	Now using \cite[Lemma II.3.4]{ladyzhenskaia1968linear} we deduce that $\enVert{m_x(\cdot,0)}_{W^{\frac{1}{2}-\frac{1}{2p}}(0,T)}$ is estimated by a constant depending on the data.
	Plugging this estimate into \eqref{eq:int m}, we see that $\enVert{m}_{C^{\alpha}([0,T];L^1)}$ is estimated as well.
	Thus, by the locally boundedness assumption on $D_\xi H$ and $D_{\xi\xi} H$ imposed on the Hamiltonian in Assumption \ref{as:H}, together with the bounds on each component in the Hamiltonian derived in Lemma~\ref{lem:uxW21p}, \ref{lem:eta continuity}, and \ref{lem:[uxm]Holder}, we see that  \eqref{eq:smoothness} is established.
\end{proof}

\section{Proof of existence for the system of equations\label{sec:existenceP} \eqref{eq:HJBs2}}

To prove the existence of solutions of \eqref{eq:HJBs2}, we will use the Leray-Schauder fixed point theorem.
\begin{proof}[Proof of Theorem~\ref{thm:main-result}(1)]
	Let $p > 1$ be arbitrarily large and $\alpha \in (0,1)$ sufficiently small as in Lemma \ref{lem:smoothness}. 
	Consider the Banach space $X$ defined as the set of all pairs $(u,m)$ with $u,u_x \in C^{\alpha/2,\alpha}([0,T] \times [0,\infty))$ and $m \in C^\alpha([0,T];L^1(0,\infty)) \cap C^{\alpha/2,\alpha}([0,T] \times [0,\infty))$, with norm
	\begin{equation}
	\enVert{(u,m)}_X = \enVert{u}_{C^{\alpha/2,\alpha}} + \enVert{u_x}_{C^{\alpha/2,\alpha}} + \enVert{m}_{C^{\alpha/2,\alpha}} + \enVert{m}_{C^{\alpha}([0,T];L^1)}.
	\end{equation}
	Let $\s{T}: X \times [0,1] \to X$ be given by $\s{T}(u,m,\lambda) = (w,\mu)$, where $(w,\mu)$ is the solution of
	\begin{equation}\label{eq:fixed pt problem}
	\begin{array}{lc}
	w_t+\frac{1}{2}\sigma^2 w_{xx}-rw+\lambda H(t,u_x,{\color{red}\pi}[u_x,m],\eta)+{\s{I}}[w]=0,&  0<t<T,\ 0<x<\infty\\
	\mu_t-\frac{1}{2}\sigma^2 \mu_{xx}-\lambda(D_\xi H(t,w_x,{\color{red}\pi}[w_x,m],\eta)\mu)_x-{\s{I}}^*[\mu]=0,& 0<t<T,\ 0<x<\infty\\
	\eta(t) = \int_0^\infty m(t,x)\dif x, & 0 < t < T,\\
	\mu(0,x)=\lambda m_0(x),\ w(T,x)=\lambda u_T(x), & 0\leq x< \infty\\
	w(t,x)=\mu(t,x)=0, & 0\leq t\leq T,  -\infty< x\leq 0.
	\end{array}
	\end{equation}
	We check that $\s{T}$ is well-defined.
	Indeed, given $(u,m,\lambda) \in X \times [0,1]$, it follows from Assumption \ref{as:H} that $H(t,u_x,{\color{red}\pi}[u_x,m],\eta) \in C^{\alpha/2,\alpha}$ with
	\begin{equation}
	\enVert{H(t,u_x,{\color{red}\pi}[u_x,m],\eta)}_{C^{\alpha/2,\alpha}} \leq C\del{\enVert{(u,m)}_X},
	\end{equation}
	where we used that by Assumption~\ref{as:H}, the Hamiltonian $H$ is locally Lipschitz continuous together with the continuity assumption on the nonlocal coupling term in Assumption~\ref{as:coupling continuity}.
	Therefore, we deduce from Lemma \ref{lem:existence Holder} that there exists a unique solution $w$ to the first equation in \eqref{eq:fixed pt problem}, satisfying an estimate
	\begin{equation} \label{eq:w estimate}
	\enVert{w}_{C^{1+\alpha/2,2+\alpha}} \leq C\del{\enVert{(u,m)}_X, \enVert{u_T}_{C^{2+\alpha}}}.
	\end{equation}
	Now, we can expand the equation for $\mu$ to get
	\begin{equation}
	\mu_t-\frac{1}{2}\sigma^2 \mu_{xx}-\lambda D_{\xi} H(t,w_x,{\color{red}\pi}[w_x,m],\eta)\mu_x - \lambda D_{\xi\xi} H(t,w_x,{\color{red}\pi}[w_x,m],\eta)w_{xx}\mu 
	-{\s{I}}[\mu]=0.
	\end{equation}
	Observe that by Assumption \ref{as:H} and \eqref{eq:w estimate} we have that
	\begin{equation*}
	\begin{split}
	\max\Big\{\enVert{D_\xi H(t,w_x,{\color{red}\pi}[w_x,m],\eta)}_{\infty},\, \enVert{D_{\xi\xi} H(t,w_x,{\color{red}\pi}[w_x,m],\eta)w_{xx}}_{\infty}\Big\}
	&\leq C\del{\enVert{w}_{C^{1+\alpha/2,2+\alpha}}} \\
	&\leq C\del{\enVert{(u,m)}_X, \enVert{u_T}_{C^{2+\alpha}}}.
	\end{split}
	\end{equation*}
	Thus by Lemma \ref{lem:existence Holder} there exists a unique solution $\mu$, satisfying the estimate
	\begin{equation} \label{eq:mu estimate1}
	\enVert{\mu}_{W^{1,2}_p} \leq C\del{\enVert{(u,m)}_X, \enVert{u_T}_{C^{2+\alpha}},\enVert{m_0}_{W^2_p}}.
	\end{equation}
	Using Lemma \ref{lem:FP basics}, Sobolev embeddings (see e.g. ~\cite[Lemma II.2.1]{ladyzhenskaia1968linear}) and the fact that $p$ is sufficiently large, we also deduce that
	\begin{equation} \label{eq:mu estimate2}
	\enVert{\mu}_{C^\alpha([0,T];L^1) \cap C^{\alpha/2,\alpha}} \leq C\del{\enVert{(u,m)}_X, \enVert{u_T}_{C^{2+\alpha}},\enVert{m_0}_{W^2_p}}.
	\end{equation}
	Hence $\s{T}$ is well-defined.
	Moreover, we have  $\s{T}(u,m,0) = (0,0)$ for all $(u,m) \in X$, i.e.~$w=\mu = 0$ by uniqueness for parabolic equations when $\lambda = 0$ in \eqref{eq:fixed pt problem}.
	
	Next, we show that $\s{T}$ is continuous.
	Let $(u_n,m_n,\lambda_n) \to (u,m,\lambda)$ in $X \times [0,1]$.
	Let $(w_n,\mu_n) = \s{T}(u_n,m_n,\lambda_n)$ and $(w,\mu) = \s{T}(u,m,\lambda)$.
	Then $\tilde w = w_n - w$ satisfies
	{\small \begin{equation}\label{eq:T differences}
		\begin{array}{lc}
		\tilde w_t+\frac{1}{2}\sigma^2 \tilde w_{xx}-r\tilde w+{\s{I}}[\tilde w]=-\lambda_n H(t,(u_n)_x,{\color{red}\pi}[(u_n)_x,m_n],\eta_n)+\lambda H(t,u_x,{\color{red}\pi}[u_x,m],\eta),&  0<t<T,\ 0<x<\infty\\
		w(T,x)=(\lambda_n-\lambda) u_T(x), & 0\leq x< \infty\\
		w(t,x)=0, & 0\leq t\leq T,  -\infty< x\leq 0
		\end{array}
		\end{equation}}
	where, as usual, $\eta(t) = \int_0^\infty m(t,x)\dif x$ and $\eta_n(t) = \int_0^\infty m_n(t,x)\dif x$.
	Using Assumption \ref{as:H} we deduce that
	\begin{equation}
	\enVert{H(t,(u_n)_x,m_n)-H(t,u_x,m)}_{C^{\alpha/2,\alpha}} \leq C\del{\enVert{(u_n,m_n)}_X, \enVert{(u,m)}_X}\enVert{(u_n-u,m_n-m)}_X
	\end{equation}
	so that, by Lemma \ref{lem:existence Holder} we obtain that
	\begin{equation} \label{eq:w1-w2}
	\begin{split}
	&\enVert{w_n-w}_{C^{1+\alpha/2,2+\alpha}}\\ 
	&\leq C\del{\enVert{(u_n,m_n)}_X, \enVert{(u,m)}_X,\enVert{u_T}_{C^{2+\alpha}}}\del{\enVert{(u_n-u,m_n-m)}_X + \abs{\lambda_n-\lambda}}.
	\end{split}
	\end{equation}
	Hence we conclude that $w_n \to w$ in $C^{1+\alpha/2,2+\alpha}$.
	Moreover,	using \eqref{eq:w1-w2} one similarly deduces that 
	\begin{equation} \label{eq:DxiHn}
	\enVert{D_\xi H(t,(w_n)_x,m_n,\eta_n)-D_\xi H(t,w_x,m,\eta)}_{C^{\alpha/2,\alpha}} \leq C\del{\enVert{(u_n,m_n)}_X, \enVert{(u,m)}_X}\enVert{(u_n-u,m_n-m)}_X.
	\end{equation}
	We note that \eqref{eq:mu estimate1} holds for each $\mu_n$, and as by assumption $(u_n,m_n)$ converges to $(u,m)$ in $X$, we obtain  a uniform estimate on  the sequence $\mu_n$ in $W_p^{1,2}$.
	On any subsequence, we can hence pass to a further subsequence that converges uniformly to a limit, which we will call $\bar \mu$.
	By \eqref{eq:DxiHn} $\bar \mu$ must satisfy
	\begin{equation}
	\bar \mu_t-\frac{1}{2}\sigma^2 \bar \mu_{xx}-\lambda(D_\xi H(t,w_x,{\color{red}\pi}[w_x,m],\eta)\bar \mu)_x-{\s{I}}[\bar \mu]=0
	\end{equation}
	in a weak sense, with $\bar \mu(0,x) = \lambda m_0(x)$.
	By uniqueness for this equation, it follows that $\bar \mu = \mu$.
	We deduce that $\mu_n \rightharpoonup \mu$ weakly in $W^{1,2}_p$ for sufficiently large $p$ and thus also $\mu_n \to \mu$ strongly in $C^\alpha([0,T];L^1) \cap C^{\alpha/2,\alpha}$.
	It follows that $\s{T}$ is continuous.
	To see that it is compact, note that if $(u_n,\mu_n)$ is bounded in $X$ and $\lambda_n \in [0,1]$, then $(w_n,\mu_n) = \s{T}(u_n,\mu_n,\lambda_n)$ is
	relatively compact 
	in $X$ by \eqref{eq:w estimate},\eqref{eq:mu estimate1}, and \eqref{eq:mu estimate2} and the Arzel\`a-Ascoli Theorem.
	
	Finally, we now invoke the a priori estimate of Lemma \ref{lem:smoothness}: any $(u,m)$ satisfying $(u,m) = \s{T}(u,m,\lambda)$ for some $\lambda \in [0,1]$ must satisfy $\enVert{(u,m)}_X \leq C$, where $C$ is the constant appearing in \eqref{eq:smoothness}.
	Therefore, by the Leray-Schauder fixed point theorem (see e.g.~\cite[Theorem 11.6]{gilbarg2015elliptic}), it follows that there exists some $(u,m) \in X$ such that $\s{T}(u,m,1) = (u,m)$. Moreover,
	by Lemma~\ref{lem:FP basics}, we know that $m\geq 0$.
	Hence we conclude that $(u,m)$ is a solution of \eqref{eq:HJBs2}.
\end{proof}
\section{Proof of uniqueness for the system of equations \eqref{eq:HJBs2}}\label{sec:uniqueness}
In this section we prove the uniqueness for the system of equations \eqref{eq:HJBs2}. 
{\color{red} Our uniqueness result is part of a wide class of  results in mean-field games which state that when the coupling among the players is rather small,
	then uniqueness can be proved.}
{\color{blue} Although our smallness condition may appear restrictive, we point out that it is precisely under this regime  that the numerical computations performed by Chan and Sircar \cite{chan2015bertrand,chan2017fracking} are supposed to hold.}
%
\begin{proof}[Proof of Theorem~\ref{thm:main-result}~(2)]
	Assume {\color{blue} $0 \leq \varepsilon \leq \varepsilon_0 \leq 1$}, where we specify $\varepsilon_0$ later in the proof.
	Let $(u_i,m_i), i=1,2$ be two solutions, and let $w = u_1 - u_2$ and $\mu = m_1 - m_2$, which satisfy \eqref{eq:u differences} and \eqref{eq:m differences}, respectively.
	Define $V = D_\xi H(t,u_{2x},{\color{red}\pi}[u_{2x},m_2],\eta_2) - D_\xi H(t,u_{1x},{\color{red}\pi}[u_{1x},m_1],\eta_1)$.
	We can rewrite \eqref{eq:m differences} as
	\begin{equation}\label{eq:m differences1}
	\mu_t-\frac{1}{2}\sigma^2 \mu_{xx}+\del{Vm_2 - D_\xi H(t,u_{1x},{\color{red}\pi}[u_{1x},m_1],\eta_1)\mu}_x-{\s{I}}^*[\mu]=0.
	\end{equation}
	Multiply \eqref{eq:m differences1} by $\mu$ and integrate by parts
	we have
	\begin{equation}\label{eq:uniq-start}
	\begin{split}
	&\frac{1}{2}\od{}{t}\int_0^\infty \mu^2(t,x)\dif x + \frac{\sigma^2}{2}\int_0^\infty \mu_x^2(t,x)\dif x\\
	&=
	\int_0^\infty \del{Vm_2 - D_\xi H(t,u_{1x},{\color{red}\pi}[u_{1x},m_1],\eta_1)\mu}\mu_x \dif x
	+ \int_{0}^\infty I[\mu](t,x)\mu(t,x)\,dx
	\end{split}
	\end{equation}
	Moreover, to see that $\int_{0}^\infty I[\mu](t,x)\mu(t,x)\,dx\leq 0$ we note that for any $0<r< 1$ we have that
	\begin{equation}
	\label{eq:uniq-int-spalt}
	\begin{split}
	\int_{0}^\infty I[\mu](t,x)\mu(t,x)\,dx
	&=
	\int_{0}^\infty \int_{-r}^r [\mu(t,x+z)-\mu(t,x)-\mu_x(t,x) z]\,\mu(t,x)\,F(dz)\,dx\\
	& \quad
	+ 	\int_{0}^\infty \int_{|z|>r} [\mu(t,x+z)-\mu(t,x)]\,\mu(t,x)\,F(dz)\,dx\\
	&	\quad
	+ \int_{0}^\infty \int_{r\leq |z|\leq 1} [-\mu_x(t,x)z]\,\mu(t,x)\,F(dz)\,dx.
	\end{split}
	\end{equation}
	For the first summand, notice that for any $0<r<1$ we have by twice applying a Taylor expansion, integration by parts, and Fubini's theorem that
	\begin{equation*}
	\begin{split}
	&\int_{0}^\infty \int_{-r}^r [\mu(t,x+z)-\mu(t,x)-\mu_x(t,x) z]\,\mu(t,x)\,F(dz)\,dx\\
	&=
	\int_{0}^\infty  \int_{-r}^r \int_{0}^1 \int_{0}^1 [\mu_{xx}(t,x+t_1 t_2 z)]\,\mu(t,x)\,t_1 z^2\,dt_2\,dt_1\,F(dz)\,dx\\
	&=
	-	\int_{0}^\infty  \int_{-r}^r \int_{0}^1 \int_{0}^1 \mu_{x}(t,x+t_1 t_2 z)\,\mu_x(t,x)\,t_1 z^2\,dt_2\,dt_1\,F(dz)\,dx\\
	&\leq
	\int_{-r}^r \int_{0}^1 \int_{0}^1 \int_{0}^\infty  \big|\mu_{x}(t,x+t_1 t_2 z)\,\mu_x(t,x)\big|\,dx\,t_1 \,z^2\,dt_2\,dt_1\,F(dz)\\
	\end{split}
	\end{equation*}
	Moreover, by the Cauchy-Schwarz inequality and as $\mu_x(t,y)=0$ for all $y\leq 0$, we see that
	\begin{equation}
	\label{eq:uniq-first}
	\begin{split}
	&
	\int_{-r}^r \int_{0}^1 \int_{0}^1 \int_{0}^\infty  \big|\mu_{x}(t,x+t_1 t_2 z)\,\mu_x(t,x)\big|\,dx\,t_1 \,
	z^2\,dt_2\,dt_1\,F(dz)\\
	& \leq 
	\int_{-r}^r \int_{0}^1 \int_{0}^1 \bigg(\int_{0}^\infty  \big|\mu_{x}(t,x+t_1 t_2 z)\big|^2\,\,dx\bigg)^{\nicefrac{1}{2}} \bigg(\int_{0}^\infty  \big|\mu_{x}(t,x) \big|^2\,\,dx\bigg)^{\nicefrac{1}{2}}\,t_1
	z^2\,dt_2\,dt_1\,F(dz)\\
	& \leq 
	\int_{-r}^r \int_{0}^1 \int_{0}^1 \|\mu_{x}(t,\cdot)\|_{L^2(0,\infty)}^2 \,t_1
	z^2\,dt_2\,dt_1\,F(dz)\\ 
	& = 
	c_r \|\mu_{x}(t,\cdot)\|_{L^2(0,\infty)}^2,
	\end{split}
	\end{equation} 
	where $c_r \to 0$ as $r \to 0$. For the second term, we can apply the same argument as in \eqref{eq:uniq-first} but with $\mu(t,x+z)$ instead of $\mu_{x}(t,x+t_1 t_2 z)$ to see that
	\begin{equation}
	\label{eq:uniq-second}
	\begin{split}
	&\int_{0}^\infty \int_{|z|>r} [\mu(t,x+z)-\mu(t,x)]\,\mu(t,x)\,F(dz)\,dx\\
	&\leq  \int_{|z|>r} \int_{0}^\infty  |\mu(t,x+z)\mu(t,x)|\,dx\,F(dz) - \int_{|z|>r} \|\mu(t,\cdot)\|_{L^2(0,\infty)}^2\,F(dz)\\
	&\leq 0.
	\end{split}
	\end{equation} 
	For the third term, as $C_r:=\int_{r\leq |z|\leq 1} z \, F(dz)<\infty$, we can use integration by parts and the boundary conditions to see that 
	\begin{equation}
	\label{eq:uniq-third}
	\begin{split}
	\int_{0}^\infty \int_{r\leq |z|\leq 1} [-\mu_x(t,x)z]\,\mu(t,x)\,F(dz)\,dx
	&= C_r
	\int_{0}^\infty  -\mu_x(t,x)\,\mu(t,x)\,dx\\
	&= C_r \int_{0}^\infty  \mu(t,x)\,\mu_x(t,x)\,dx =0.
	\end{split}
	\end{equation}
	Therefore, we conclude from \eqref{eq:uniq-int-spalt}--
	\eqref{eq:uniq-third} that indeed 
	\begin{equation*}
	\int_{0}^\infty I[\mu](t,x)\mu(t,x)\,dx\leq 0.
	\end{equation*}
	We hence deduce from \eqref{eq:uniq-start} that
	\begin{equation}\label{eq:uniq-start-1}
	\begin{split}
	\frac{1}{2}\od{}{t}\int_0^\infty \mu^2(t,x)\dif x + \frac{\sigma^2}{2}\int_0^\infty \mu_x^2(t,x)\dif x
	&\leq
	\int_0^\infty \del{Vm_2 - D_\xi H(t,u_{1x},{\color{red}\pi}[u_{1x},m_1],\eta_1)\mu}\mu_x(t,x) \dif x.
	\end{split}
	\end{equation}
	Next, we use Assumption~\ref{as:H}, the a priori bounds from Lemma \ref{lem:smoothness}, and Assumption~\ref{as:small Lipschitz} to deduce that
	\begin{equation}\label{eq:D-xi-mu estimate}
	\abs{D_\xi H(t,u_{1x},{\color{red}\pi}[u_{1x},m_1],\eta_1)\mu} \leq C
	\end{equation}
	and
	\begin{equation} \label{eq:V estimate}
	\abs{V} \leq C\del{\abs{u_{1x}-u_{2x}} + \varepsilon\abs{{\color{red}\pi}[u_{1x},m_1]-{\color{red}\pi}[u_{2x},m_2]} + \varepsilon\abs{\eta_1-\eta_2}}
	\end{equation}
	for some constant $C$ depending only on the data. Therefore, 
	using Peter--Paul 
	inequality, we get that
	\begin{equation*}
	\begin{split}
	&\od{}{t}\int_0^\infty \mu^2(t,x)\dif x + \frac{\sigma^2}{2}\int_0^\infty \mu_x^2(t,x)\dif x\\
	&\leq C\int_0^\infty \mu^2(t,x) \dif x
	+ C\int_0^\infty \abs{u_{1x}(t,x)-u_{2x}(t,x)}^2 m_2^2(t,x) \dif x 
	\\
	& \quad 
	+ C\varepsilon\abs{{\color{red}\pi}[u_{1x},m_1]-{\color{red}\pi}[u_{2x},m_2]}^2\int_0^\infty m_2^2(t,x) \dif x 
	+ C\varepsilon\abs{\eta_1(t)-\eta_2(t)}^2\int_0^\infty m_2^2(t,x) \dif x.
	\end{split}
	\end{equation*}
	Moreover, by Lemma \ref{lem:smoothness} we have an a priori bound on $m_1,m_2$ in all $L^p$-spaces.
	Thus, using that $m_2\in L^\infty([0,T];L^2(\mathbb{R}))$, we see that
	\begin{equation} \label{eq:mu2 estimate}
	\begin{split}
	&\od{}{t}\int_0^\infty \mu^2(t,x)\dif x + \frac{\sigma^2}{2}\int_0^\infty \mu_x^2(t,x)\dif x\\
	&\leq C\int_0^\infty \mu^2(t,x) \dif x
	+ C\int_0^\infty \abs{u_{1x}(t,x)-u_{2x}(t,x)}^2 m_2(t,x) \dif x + C\varepsilon\abs{{\color{red}\pi}[u_{1x},m_1]-{\color{red}\pi}[u_{2x},m_2]}^2\\
	& \quad
	+ C\varepsilon\abs{\eta_1(t)-\eta_2(t)}^2.
	\end{split}
	\end{equation}
	Furthermore, by Assumption \ref{as:coupling continuity} and the a priori bounds on $u_{1x},u_{2x},m_1,m_2$, we have
	\begin{equation} \label{eq:[ux,m]differences}
	\begin{split}
	&\big|{\color{red}\pi}[u_{1x},m_1]-{\color{red}\pi}[u_{2x},m_2]\big|^2\\
	&\leq C\int_0^\infty \mu^2(t,x)\dif x
	+ C\int_0^\infty \big|u_{1x}(t,x) - u_{2x}(t,x)\big|^2 \del{m_1(t,x)+m_2(t,x)}\dif x
	+ C\big|\eta_1(t)-\eta_2(t)\big|^2.
	\end{split}
	\end{equation}
	Thus by first multiplying \eqref{eq:mu2 estimate} with an integrating factor of the form $e^{-Ct}$ and then
	integrating \eqref{eq:mu2 estimate}
	we get
	\begin{equation} \label{eq:mu2 estimate2}
	\begin{split}
	&\int_0^\infty \mu^2(t,x)\dif x + \int_0^t\int_0^\infty \mu_x^2(\tau,x)\dif x \dif \tau\\
	&\leq C\int_0^t\int_0^\infty \abs{u_{1x}(\tau,x) - u_{2x}(\tau,x)}^2 \del{m_1(t,x)+m_2(t,x)}\dif x \dif \tau + C\varepsilon\int_0^t\abs{\eta_1(\tau)-\eta_2(\tau)}^2 \dif \tau.
	\end{split}
	\end{equation}
	Next, to estimate $\eta_1-\eta_2 = \int_0^\infty \mu \dif x$, we use an argument similar in spirit to that of Lemma \ref{lem:eta continuity}.
	That is, let $\zeta(x)$ be a smooth function such that $0 \leq \zeta \leq 1$, $\zeta'$ and $\zeta''$ have compact support in $(0,M)$ for some $M > 0$, and $\zeta \equiv 1$ for $x \geq M$.
	We multiply \eqref{eq:m differences1} by $\zeta$ and integrate by parts to get that
	\begin{equation} \label{eq:eta differences1}
	\begin{split}
	\int_0^\infty \! \zeta(x)\mu(t,x)\dif x 
	&= \int_0^t\!\int_0^\infty\!\!\del{-\tfrac{\sigma^2}{2}\mu_x(\tau,x) + Vm_2(\tau,x) - D_\xi H(t,u_{1x},{\color{red}\pi}[u_{1x},m_1],\eta_1)\mu(\tau,x)}\zeta'(x) \dif x \dif \tau \\
	& \quad
	+ \int_0^t\int_0^\infty \s{I}[\zeta](\tau,x)\mu(\tau,x) \dif x\dif \tau.
	\end{split}
	\end{equation}
	Now, to deal with the nonlocal term, we apply a Taylor expansion and integration by parts to see that
	\begin{equation*}
	\begin{split}
	\int_0^t\int_0^\infty\!\! \s{I}[\zeta](\tau,x)\mu(\tau,x)\dif x \dif \tau
	&=
	\int_0^t\int_0^\infty\!\! \int_{-1}^1 \int_0^1 \int_0^1  \zeta''(x+\tau'' \tau' z) \tau' z^2 \mu(\tau,x)\dif \tau'' \dif \tau' F(\dif z) \dif x\dif \tau \\
	& \quad 
	+ \int_0^t\int_0^\infty \int_{\bb{R}\setminus (-1,1)} \del{\zeta(x+z)-\zeta(x)}\mu(\tau,x)\,F(\dif z) \dif x\dif \tau \\
	&= -\int_0^t\int_0^\infty \!\!\int_{-1}^1 \int_0^1 \int_0^1 \zeta'(x+\tau'' \tau' z) \tau' z^2 \mu_x(\tau,x)\dif \tau'' \dif \tau' \,F(\dif z) \dif x\dif \tau \\
	& \quad 
	+ \int_0^t\int_0^\infty \int_{\bb{R}\setminus (-1,1)} \del{\zeta(x+z)-\zeta(x)}\mu(\tau,x)\,F(\dif z) \dif x \dif \tau
	\end{split}
	\end{equation*}
	To estimate the first term, observe that by Fubini's theorem, we have for $C_1:= \int_{-1}^1 z^2 \,F(\dif z)$ that
	\begin{equation}
	\abs{\int_0^t\!\int_0^\infty \!\!\int_{-1}^1 \int_0^1 \!\int_0^1 \!  \zeta'(x+\tau'' \tau' z) \tau' z^2 \mu_x(\tau,x)\dif \tau'' \dif \tau' \,F(\dif z) \dif x \dif \tau}
	\leq C_1\enVert{\zeta'}_2\!\del{\int_0^t \! \int_0^\infty\!\! \mu_x^2(\tau,x) \dif x \dif \tau\!}^{1/2}\!.
	\end{equation}
	Moreover, for the second term, 
	we use that
	$0 \leq \zeta \leq 1$ and 
	$\zeta \equiv 1$ for $x \geq M$ 
	to obtain that
	\begin{equation*}
	\begin{split}
	&\int_0^t\int_0^\infty \int_{\{|z|>1\}} \del{\zeta(x+z)-\zeta(x)}\mu(\tau,x)\,F(\dif z) \dif x\dif \tau \\
	&= \int_0^t\int_0^M \int_{\{|z|>1\}} \del{\zeta(x+z)-\zeta(x)}\mu(\tau,x)\,F(\dif z) \dif x \dif \tau\\
	& \quad
	+ \int_0^t\int_{-\infty}^{-1} \int_M^{M-z}\del{\zeta(x+z)-1}\mu(\tau,x)\dif x \,F(\dif z)\dif \tau \\
	&\leq \del{2M^{1/2}C_2+ T^{1/2}C_3}\del{\int_0^t\int_0^\infty \mu^2(\tau,x) \dif x\dif \tau}^{1/2},	
	\end{split}
	\end{equation*}
	where $C_2 := \int_{\{|z|>1\}} F(\dif z)$ and $C_3:= \int_{-\infty}^{-1} \abs{z}^{1/2} F(\dif z)$ is finite due to Assumption~\ref{as:levy}.
	Therefore, we obtain that
	\begin{equation} \label{eq:nonlocal estimate uniqueness}
	\begin{split}
	&\int_0^t\int_0^\infty \s{I}[\zeta](\tau,x)\mu(\tau,x)\dif x \dif \tau\\
	&\leq 
	C_1\enVert{\zeta'}_2\del{\int_0^t \int_0^\infty \mu_x^2(\tau,x) \dif x \dif \tau}^{1/2}\!+ \,\del{2M^{1/2}C_2+ T^{1/2}C_3}\del{\int_0^t\int_0^\infty \mu^2(\tau,x) \dif x\dif \tau}^{1/2}
	\end{split}
	\end{equation}
	Moreover, using \eqref{eq:nonlocal estimate uniqueness} in \eqref{eq:eta differences1}, \eqref{eq:D-xi-mu estimate}, and that we have a priori bounds in all $L^p$-spaces for $m_2$ by Lemma~\ref{lem:smoothness}
	we deduce that
	\begin{equation} \label{eq:eta differences2}
	\begin{split}
	\int_0^\infty \zeta(x)\mu(t,x)\dif x
	\leq C\del{\enVert{\zeta'}_2+1}\Bigg[\del{\int_0^t\!\int_0^\infty \!\!\mu^2(\tau,x) \dif x\dif \tau}^{1/2} \!\!\!&+ \del{\int_0^t\int_0^\infty \mu_x^2(\tau,x) \dif x\dif \tau}^{1/2}\\ &+\del{\int_0^t\int_0^\infty V^2(\tau,x) \dif x\dif \tau}^{1/2}\Bigg].
	\end{split}
	\end{equation}
	On the other hand, as $1-\zeta(x) = 0$ for $x \geq M$, we have
	\begin{equation}\label{eq:eta differences3}
	\int_0^\infty \del{1-\zeta(x)}\mu(t,x)\dif x \leq M^{1/2}\del{\int_0^\infty\mu^2(t,x)\dif x}^{1/2}.
	\end{equation}
	Therefore, by adding together \eqref{eq:eta differences2} and \eqref{eq:eta differences3} we get
	\begin{equation}
	\begin{split}
	\label{eq:eta differences4}
	\eta_1(t) - \eta_2(t)
	= \int_0^\infty \mu(t,x)\dif x
	&\leq M^{1/2}\del{\int_0^\infty\mu^2(t,x)\dif x}^{1/2}\\
	&\quad +C\del{\enVert{\zeta'}_2+1}\Bigg[\del{\int_0^t\!\int_0^\infty\! \mu^2(\tau,x) \dif x\dif \tau}^{1/2}
	\\
	& \quad
	+ \del{\int_0^t\!\int_0^\infty\! \mu_x^2(\tau,x) \dif x\dif \tau}^{1/2} 
	+\del{\int_0^t\!\int_0^\infty\! V^2(\tau,x) \dif x\dif \tau}^{1/2}\Bigg].
	\end{split}
	\end{equation}
	Note that by the same argument, we get the same estimate for $\eta_2(t) - \eta_1(t)$ on the left-hand side of\eqref{eq:eta differences4}.
	Moreover, observe that for any fixed  $M > 1$ we can  choose $\zeta$ in such a way that it additionally satisfies $\abs{\zeta'} \leq 1$, which in turn implies that $\enVert{\zeta'}_2 \leq M^{1/2}$.
	%
	Therefore, we square both sides of \eqref{eq:eta differences4}, using \eqref{eq:V estimate} and \eqref{eq:[ux,m]differences}, to  obtain that
	\begin{equation*}
	\begin{split}
	&\abs{\eta_1(t) - \eta_2(t)}^2\\
	&\leq 
	C\del{\int_0^\infty\mu^2(t,x)\dif x + \int_0^t\int_0^\infty \mu^2(\tau,x) \dif x\dif \tau + \int_0^t\int_0^\infty \mu_x^2(\tau,x) \dif x\dif \tau +\int_0^t\int_0^\infty V^2(\tau,x) \dif x\dif \tau}\\
	&\leq C\del{\int_0^\infty\mu^2(t,x)\dif x + \int_0^t\int_0^\infty \mu^2(\tau,x) \dif x\dif \tau + \int_0^t\int_0^\infty \mu_x^2(\tau,x) \dif x\dif \tau}\\
	& \quad	 +C\del{\int_0^t\int_0^\infty \abs{u_{1x}(\tau,x)-u_{2x}(\tau,x)}^2(m_1(\tau,x)+m_2(\tau,x)) \dif x\dif \tau
		+ \int_0^t\abs{\eta_1(\tau) - \eta_2(\tau)}^2 \dif \tau}.
	\end{split}
	\end{equation*}
	Therefore Gronwall's Lemma ensures that for each time $t$, we have
	\begin{equation}
	\begin{split}
	\label{eq:eta differences5a}
	\abs{\eta_1(t) - \eta_2(t)}^2
	&\leq C\del{\int_0^\infty\mu^2(t,x)\dif x + \int_0^t\int_0^\infty \mu^2(\tau,x) \dif x\dif \tau + \int_0^t\int_0^\infty \mu_x^2(\tau,x) \dif x\dif \tau}\\
	& \quad
	+C\int_0^t\int_0^\infty \Big[\abs{u_{1x}(\tau,x)-u_{2x}(\tau,x)}^2\big(m_1(\tau,x)+m_2(\tau,x)\big)\Big] \dif x\dif \tau.
	\end{split}
	\end{equation}
	This ensures that for each time $t$, we have that
	\begin{equation}
	\begin{split}
	\label{eq:eta differences5}
	\sup_{0\leq \tau \leq t}\abs{\eta_1(\tau) - \eta_2(\tau)}^2
	&\leq C \del{ 
		\sup_{0\leq \tau \leq t}\int_0^\infty\mu^2(\tau,x)\dif x + \int_0^t\!\int_0^\infty \mu^2(\tau,x) \dif x\dif \tau + \int_0^t\!\int_0^\infty \mu_x^2(\tau,x) \dif x\dif \tau}\\
	& \quad
	+C\int_0^t\int_0^\infty \Big[\abs{u_{1x}(\tau,x)-u_{2x}(\tau,x)}^2\big(m_1(\tau,x)+m_2(\tau,x)\big)\Big] \dif x\dif \tau.
	\end{split}
	\end{equation}
	Putting \eqref{eq:eta differences5} into \eqref{eq:mu2 estimate2}, assuming $\varepsilon$ is small enough, and applying Gronwall's Lemma again, we deduce for each time $t$ that
	\begin{equation*}
	\begin{split} 
	&(1-C\varepsilon)\del{\int_0^\infty \mu^2(t,x)\dif x + \int_0^t\int_0^\infty \mu_x^2(\tau,x)\dif x \dif \tau}\\
	&\leq 
	C\int_0^t\int_0^\infty \abs{u_{1x}(\tau,x) - u_{2x}(\tau,x)}^2 \del{m_1(\tau,x)+m_2(\tau,x)}\dif x \dif \tau,
	\end{split}
	\end{equation*}
	which then directly ensures that
	\begin{equation} \label{eq:mu2 estimate3}
	\begin{split}
	&\sup_{t \in [0,T]}\int_0^\infty \mu^2(t,x)\dif x + \int_0^T\int_0^\infty \mu_x^2(\tau,x)\dif x \dif \tau\\
	&\leq \frac{C}{1-C\varepsilon}\int_0^T\int_0^\infty \abs{u_{1x}(\tau,x) - u_{2x}(\tau,x)}^2 \del{m_1(\tau,x)+m_2(\tau,x)}\dif x \dif \tau.
	\end{split}
	\end{equation}
	Finally, it remains to estimate $\int_0^T\int_0^\infty \abs{u_{1x}(\tau,x) - u_{2x}(\tau,x)}^2 \del{m_1(\tau,x)+m_2(\tau,x)}\dif x \dif \tau$. To that end, {\color{red}using here the notation $ \pi_i:={\color{red}\pi}[u_{ix},m_i]$, $i=1,2$,} we can apply \eqref{eq:energy differences}  to see that 
	{\small
		\begin{equation} \label{eq:energy differences1}
		\begin{split}
		&\int_0^T\!\! \int_0^\infty\! e^{-rt}\del{H(t,u_{2x},{\color{red}\pi_2},\eta_2)-H(t,u_{2x},{\color{red}\pi_1},\eta_1)}m_1(t,x) \dif x \dif t\\
		& \quad 
		+ \int_0^T \! \!\!\int_0^\infty\! e^{-rt}\del{H(t,u_{1x},{\color{red}\pi_1},\eta_1)-H(t,u_{1x},{\color{red}\pi_2},\eta_2)}m_2(t,x) \dif x \dif t
		\\ 
		&=\!\int_0^T\!\!\! \int_0^\infty\! \! e^{-rt}\del{H(t,u_{2x},{\color{red}\pi_1},\eta_1)\!-\!H(t,u_{1x},{\color{red}\pi_1},\eta_1)\!-\!D_\xi H(t,u_{1x},{\color{red}\pi_1},\eta_1) (u_{2x}\!-\!u_{1x})}m_1(t,x) \dif x \dif t\\
		& \quad
		\!+ \int_0^T\! \!\!\int_0^\infty\!\!\! e^{-rt}\del{H(t,u_{1x},{\color{red}\pi_2},\eta_2)\!-\!H(t,u_{2x},{\color{red}\pi_2},\eta_2)  \!-\!D_\xi H(t,u_{2x},{\color{red}\pi_2},\eta_2) (u_{1x}\!-\!u_{2x})}m_2(t,x) \dif x \dif t.
		\end{split}
		\end{equation}}
	Now, as by  Assumption \ref{as:smoothness} we have that the Hamiltonian is uniformly convex in $\xi$, we deduce that
	\begin{equation} \label{eq:energy differences2}
	\begin{split}
	&\int_0^T \int_0^\infty e^{-rt}\del{H(t,u_{2x},{\color{red}\pi}[u_{2x},m_2],\eta_2)(t,x)-H(t,u_{2x},{\color{red}\pi}[u_{1x},m_1],\eta_1)(t,x)}m_1(t,x) \dif x \dif t\\
	& \quad + \int_0^T \int_0^\infty e^{-rt}\del{H(t,u_{1x},{\color{red}\pi}[u_{1x},m_1],\eta_1)(t,x)-H(t,u_{1x},{\color{red}\pi}[u_{2x},m_2],\eta_2)(t,x)}m_2(t,x) \dif x \dif t
	\\ 
	&\geq \frac{1}{C}\int_0^T \int_0^\infty e^{-rt}\abs{u_{2x}(t,x)-u_{1x}(t,x)}^2 (m_1(t,x)+m_2(t,x)) \dif x \dif t.
	\end{split}
	\end{equation}
	Moreover, by applying a Taylor expansion and using \eqref{eq:epsilon Lipschitz} applied to $D_\upsilon H$ and $D_\eta H$, we have pointwise (i.e.\ for any $t,x$) that
	\begin{equation*}
	\begin{split}
	H(t,u_{2x},{\color{red}\pi}[u_{2x},m_2],\eta_2)-H(t,u_{2x},{\color{red}\pi}[u_{1x},m_1],\eta_1)
	&\leq 
	C\varepsilon\del{\abs{{\color{red}\pi}[u_{2x},m_2](t,x)-{\color{red}\pi}[u_{1x},m_1]}^2 + \abs{\eta_1-\eta_2}^2}\\
	& \quad + D_\upsilon H(t,u_{2x},{\color{red}\pi}[u_{2x},m_2],\eta_2)\del{{\color{red}\pi}[u_{2x},m_2]-{\color{red}\pi}[u_{1x},m_1]}\\
	& \quad 
	+ D_\eta H(t,u_{2x},{\color{red}\pi}[u_{2x},m_2],\eta_2)(\eta_1-\eta_2)
	\end{split}
	\end{equation*}
	as well as
	\begin{equation*}
	\begin{split}
	H(t,u_{1x},{\color{red}\pi}[u_{1x},m_1],\eta_1)-H(t,u_{1x},{\color{red}\pi}[u_{2x},m_2],\eta_2)
	&\leq C\varepsilon\del{\abs{{\color{red}\pi}[u_{2x},m_2]-{\color{red}\pi}[u_{1x},m_1]}^2 + \abs{\eta_1-\eta_2}^2}\\
	& \quad 
	- D_\upsilon H(t,u_{1x},{\color{red}\pi}[u_{1x},m_1],\eta_1)\del{{\color{red}\pi}[u_{2x},m_2]-{\color{red}\pi}[u_{1x},m_1]}\\
	& \quad
	- D_\eta H(t,u_{1x},{\color{red}\pi}[u_{1x},m_1],\eta_1)(\eta_1-\eta_2).
	\end{split}
	\end{equation*}
	Therefore, using that
	$0\leq\int_0^\infty m_i(t,x)\,dx\leq 1$, $i:=1,2,$ we see that
	\begin{equation*} 
	\begin{split}
	&\int_0^T \int_0^\infty e^{-rt}\del{H(t,u_{2x},{\color{red}\pi}[u_{2x},m_2],\eta_2)(t,x)-H(t,u_{2x},{\color{red}\pi}[u_{1x},m_1],\eta_1)(t,x)}m_1(t,x) \dif x \dif t\\
	& \quad + \int_0^T \int_0^\infty e^{-rt}\del{H(t,u_{1x},{\color{red}\pi}[u_{1x},m_1],\eta_1)(t,x)-H(t,u_{1x},{\color{red}\pi}[u_{2x},m_2],\eta_2)(t,x)}m_2(t,x) \dif x \dif t
	\\ 
	&\leq C\varepsilon\int_0^T e^{-rt}\del{\abs{{\color{red}\pi}[u_{2x},m_2](t)-{\color{red}\pi}[u_{1x},m_1](t)}^2 + \abs{\eta_1(t)-\eta_2(t)}^2}\dif t
	\\
	& \quad
	+ \int_0^T \int_0^\infty e^{-rt}D_\upsilon H(t,u_{2x},{\color{red}\pi}[u_{2x},m_2],\eta_2)\del{{\color{red}\pi}[u_{2x},m_2]-{\color{red}\pi}[u_{1x},m_1]}m_1(t,x) \dif x \dif t\\
	& \quad
	+ \int_0^T \int_0^\infty e^{-rt}D_\eta H(t,u_{2x},{\color{red}\pi}[u_{2x},m_2],\eta_2)(\eta_1-\eta_2)m_1(t,x) \dif x \dif t\\
	&\quad
	- \int_0^T \int_0^\infty e^{-rt}D_\upsilon H(t,u_{1x},{\color{red}\pi}[u_{1x},m_1],\eta_1)\del{{\color{red}\pi}[u_{2x},m_2]-{\color{red}\pi}[u_{1x},m_1]}m_2 \dif x \dif t\\
	& \quad 
	- \int_0^T \int_0^\infty e^{-rt}D_\eta H(t,u_{1x},{\color{red}\pi}[u_{1x},m_1],\eta_1)(\eta_1-\eta_2)m_2(t,x) \dif x \dif t.
	\end{split}
	\end{equation*}
	This, together with a repeated use of Assumption~\ref{as:small Lipschitz} hence ensures that 
	\begin{equation}
	\begin{split}
	\label{eq:energy differences3}
	&	\int_0^T \int_0^\infty e^{-rt}\del{H(t,u_{2x},{\color{red}\pi}[u_{2x},m_2],\eta_2)(t,x)-H(t,u_{2x},{\color{red}\pi}[u_{1x},m_1],\eta_1)(t,x)}m_1(t,x) \dif x \dif t\\
	&\quad +
	\int_0^T \int_0^\infty e^{-rt}\del{H(t,u_{1x},{\color{red}\pi}[u_{1x},m_1],\eta_1)(t,x)-H(t,u_{1x},{\color{red}\pi}[u_{2x},m_2],\eta_2)(t,x)}m_2(t,x) \dif x \dif t\\
	&\leq
	C\varepsilon\int_0^T e^{-rt}\del{\abs{{\color{red}\pi}[u_{2x},m_2](t)-{\color{red}\pi}[u_{1x},m_1](t)}^2 + \abs{\eta_1(t)-\eta_2(t)}^2}\dif t\\
	& \quad
	+ C\varepsilon\int_0^T\!\! \int_0^\infty\!\! e^{-rt}\abs{u_{1x}(t,x)-u_{2x}(t,x)}\Big(\abs{{\color{red}\pi}[u_{2x},m_2]-{\color{red}\pi}[u_{1x},m_1]}(t)+\abs{\eta_1-\eta_2}(t)\Big) m_1(t,x) \dif x \dif t\\
	& \quad
	+ \int_0^T\!\! \int_0^\infty\!\! e^{-rt}D_\upsilon H(t,u_{1x},{\color{red}\pi}[u_{1x},m_1],\eta_1)\del{{\color{red}\pi}[u_{2x},m_2]-{\color{red}\pi}[u_{1x},m_1]}(m_1-m_2)(t,x) \dif x \dif t\\
	& \quad
	+ \int_0^T \int_0^\infty e^{-rt}D_\eta H(t,u_{1x},{\color{red}\pi}[u_{1x},m_1],\eta_1)(\eta_1-\eta_2)(m_1-m_2)(t,x) \dif x \dif t\\
	&\leq C\varepsilon\int_0^T e^{-rt}\del{\abs{{\color{red}\pi}[u_{2x},m_2](t)-{\color{red}\pi}[u_{1x},m_1](t)}^2 + \abs{\eta_1(t)-\eta_2(t)}^2}\dif t\\
	& \quad 
	+ C\varepsilon\int_0^T \int_0^\infty e^{-rt}\del{\abs{u_{1x}(t,x)-u_{2x}(t,x)}^2 m_1(t,x) + \abs{m_1(t,x)-m_2(t,x)}^2} \dif x \dif t.
	\end{split}
	\end{equation}
	Therefore, taking  \eqref{eq:[ux,m]differences} into account, \eqref{eq:energy differences3} is reduced to
	\begin{equation} \label{eq:energy differences4}
	\begin{split}
	&\int_0^T \int_0^\infty e^{-rt}\del{H(t,u_{2x},{\color{red}\pi}[u_{2x},m_2],\eta_2)-H(t,u_{2x},{\color{red}\pi}[u_{1x},m_1],\eta_1)}m_1(t,x) \dif x \dif t\\
	& \quad 
	+ \int_0^T \int_0^\infty e^{-rt}\del{H(t,u_{1x},{\color{red}\pi}[u_{1x},m_1],\eta_1)-H(t,u_{1x},{\color{red}\pi}[u_{2x},m_2],\eta_2)}m_2(t,x) \dif x \dif t
	\\ 
	&\leq C\varepsilon\int_0^T e^{-rt}\abs{\eta_1(t)-\eta_2(t)}^2\dif t\\
	& \quad + C\varepsilon\int_0^T \int_0^\infty \del{\abs{u_{1x}(t,x)-u_{2x}(t,x)}^2(m_1(t,x)+m_2(t,x)) + \abs{m_1(t,x)-m_2(t,x)}^2} \dif x \dif t.
	\end{split}
	\end{equation}
	Combining \eqref{eq:energy differences4} with \eqref{eq:energy differences2} hence ensures that
	\begin{equation} \label{eq:energy differences5}
	\begin{split}
	&\del{\frac{1}{C}-C\varepsilon}\int_0^T \int_0^\infty e^{-rt}\abs{u_{2x}(t,x)-u_{1x}(t,x)}^2 (m_1(t,x)+m_2(t,x)) \dif x \dif t\\
	&\leq C\varepsilon\int_0^T e^{-rt}\abs{\eta_1(t)-\eta_2(t)}^2\dif t
	+ C\varepsilon\int_0^T \int_0^\infty \abs{m_1(t,x)-m_2(t,x)}^2 \dif x \dif t.
	\end{split}
	\end{equation}
	Thus, Combining \eqref{eq:energy differences5}, \eqref{eq:eta differences5}, and \eqref{eq:mu2 estimate3}, we get, after adjusting the value of $C$, that
	\begin{equation} \label{eq:mu2 estimate4}
	\begin{split}
	&\sup_{t \in [0,T]}\int_0^\infty \mu^2(t,x)\dif x + \int_0^T\int_0^\infty \mu_x^2(\tau,x)\dif x \dif \tau\\
	&\leq
	\frac{C\varepsilon}{1-C\varepsilon}\del{\sup_{t \in [0,T]}\int_0^\infty \mu^2(t,x)\dif x + \int_0^T\int_0^\infty \mu_x^2(\tau,x)\dif x \dif \tau}.
	\end{split}
	\end{equation}
	Now, we recall that $C$ depends only on the data.
	Hence there exists $\varepsilon_0>0$ small enough such that whenever $0<\varepsilon\leq \varepsilon_0$ we can
	deduce from \eqref{eq:mu2 estimate4} that $\mu \equiv 0$, i.e.~$m_1 = m_2$. This  in turn immediately implies that $\eta_1 = \eta_2$.
	Moreover, \eqref{eq:energy differences5} ensures that $u_{2x}=u_{1x}$ on the support of $m_1$.  
	By appealing to Assumption \ref{as:coupling continuity} we hence obtain that ${\color{red}\pi}[u_{1x},m_1] = {\color{red}\pi}[u_{2x},m_2]$.
	Finally, we may now conclude that $u_1,u_2$ satisfy the same parabolic equation for given coefficients determined by $m_1$ and ${\color{red}\pi}[u_{1x},m_1]$. Thus by considering the corresponding parabolic equation for  $u_1-u_2$  and using Lemma \ref{lem:existence Holder}, we conclude that $u_1=u_2$.
\end{proof}
%
%





\bibliography{mybibfile}

\end{document}